\documentclass[twoside,11pt]{article}

\usepackage[preprint]{jmlr2e}

\usepackage{url}
\usepackage{enumitem}
\usepackage{makecell}
\usepackage{diagbox}
\usepackage{subfigure}
\usepackage{amsmath, bm}
\usepackage{dsfont}
\usepackage{amssymb}
\usepackage{multirow}
\usepackage{varwidth}
\usepackage{pifont}
\usepackage{makecell}
\usepackage{wrapfig}
\usepackage{multicol}
\usepackage{graphicx} 
\usepackage{comment}
\usepackage{color}
\usepackage{xcolor}
\usepackage{colortbl,booktabs}
\usepackage{cleveref}
\usepackage{braket}

\definecolor{Tianlong_color}{rgb}{0.858, 0.188, 0.478}

\definecolor{Xiaohan_color}{rgb}{0.098, 0.643, 0.071}

 \makeatletter
\renewcommand\normalsize{%
\@setfontsize\normalsize{11}{12.8}}
\makeatother
\normalsize

\usepackage[font=small]{caption}

\newcommand{\cT}{\mathcal{T}}

\newcommand{\mE}{\mathbb{E}}

\newcommand{\x}{\mathbf{x}}

\newcommand{\z}{\mathbf{z}}

\newcommand{\RR}{\mathds{R}} 

\DeclareMathOperator*{\argmin}{arg\,min}

\jmlrheading{1}{2021}{1-48}{4/00}{10/00}{meila00a}{($\alpha$-$\beta$) T. Chen, X. Chen, W. Chen, H. Heaton, J. Liu, Z. Wang and W. Yin}

\ShortHeadings{Learning to Optimize: A Primer and A Benchmark}{($\alpha$-$\beta$) T. Chen, X. Chen, W. Chen, H. Heaton, J. Liu, Z. Wang and W. Yin}
\firstpageno{1}

\begin{document}
\setcitestyle{square,numbers}

\title{Learning to Optimize: A Primer and A Benchmark}

\author{
\center
\name Tianlong Chen, Xiaohan Chen, Wuyang Chen, Zhangyang Wang\thanks{All authors are listed in alphabetic order. Correspondence shall be addressed to Z. Wang.} \email \\ \{tianlong.chen,xiaohan.chen,wuyang.chen,atlaswang\}@utexas.edu \\
       \addr Department of Electrical and Computer and Engineering\\
      The University of Texas at Austin, 
       Austin, TX 78712, USA
       \AND
\name Howard Heaton \email \\ hheaton@ucla.edu \\
       \addr Department of Mathematics, 
       University of California Los Angeles, 
       Los Angeles, CA 90095, USA
       \AND
       \name Jialin Liu, Wotao Yin \email \\ \{jialin.liu,wotao.yin\}@alibaba-inc.com \\
       \addr Alibaba US, Damo Academy, Decision Intelligence Lab, Bellevue, WA 98004, USA \\
       }

\editor{~}

\maketitle

\begin{abstract}

Learning to optimize (\textbf{L2O}) is an emerging approach that leverages machine learning to develop optimization methods, aiming at reducing the laborious iterations of hand engineering. It automates the design of an optimization method based on its performance on a set of training problems. This data-driven procedure generates methods that can efficiently solve problems similar to those in the training. In sharp contrast, the typical and traditional designs of optimization methods are theory-driven, so they obtain performance guarantees over the classes of problems specified by the theory.
The difference makes L2O suitable for repeatedly solving a certain type of optimization problems over a specific distribution of data, while it typically fails on out-of-distribution problems.
The practicality of L2O depends on the type of target optimization, the chosen architecture of the method to learn, and the training procedure. This new paradigm has motivated a community of researchers to explore L2O and report their findings.

This article is poised to be the first comprehensive survey and benchmark of L2O for continuous optimization. We set up taxonomies, categorize existing works and research directions, present insights, and identify open challenges. We also benchmarked many existing L2O approaches on a few but representative optimization problems. 
For reproducible research and fair benchmarking purposes, we released our software implementation and data in the package \textbf{Open-L2O} at \url{https://github.com/VITA-Group/Open-L2O}.





\end{abstract}





\section{Introduction}

\subsection{Background and Motivation}
Classic optimization methods are typically hand-built by optimization experts based on  theories and their experience. As a paradigm shift from this conventional design, \textit{learning to optimize} (\textbf{L2O}) uses machine learning to improve an optimization method or even generate a completely new method. This paper provides a timely and up-to-date review of the rapid growing body of L2O results, with a focus on continuous optimization. 

\begin{figure}
    \centering
    \subfigure[\footnotesize{Classic Optimizer}]{\includegraphics[height=1.1in]{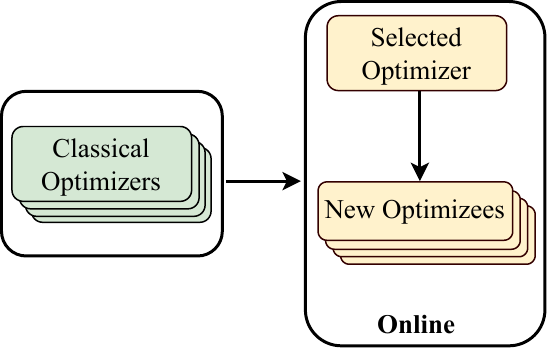}} \hspace{2.5em}
    \subfigure[\footnotesize{Learned Optimizer by L2O}]{\includegraphics[height=1.1in]{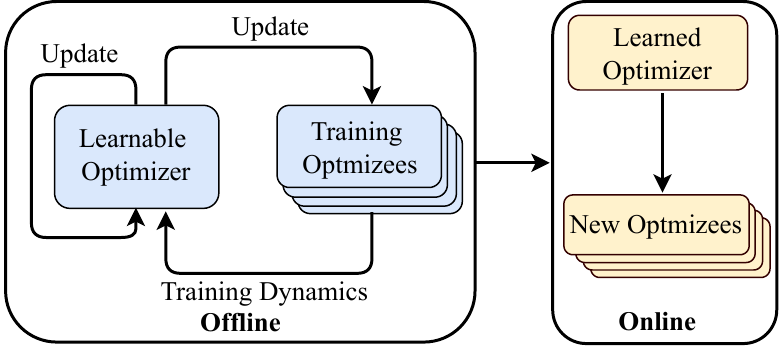}}
    \caption{(a) Classic optimizers are manually designed; they usually have few or no tuning parameters; (b) Learned optimizers are trained in an L2O framework over a set of similar optimizees (called a task distribution) and designed to solve unseen optimizees from the same distribution.
    }  
    \label{fig: optimizer-comp}
\end{figure}

Classic optimization methods are built upon components that are basic methods --- such as gradient descent, conjugate gradient, Newton steps, Simplex basis update, and stochastic sampling --- in a theoretically justified manner. Most conventional optimization methods can be written in a few lines, and their performance is guaranteed by their theories. To solve an optimization problem in practice, we can select a method that supports the problem type at hand and expect the method to return a solution no worse than its guarantee. 

L2O is an alternative paradigm that develops an optimization method by training, i.e., learning from its performance on sample problems. 
The method may lack a solid theoretical basis but improves its performance during the training process. 
The training process often occurs offline and is time consuming. However, the online application of the method is (aimed to be) time saving. When it comes to problems where the target solutions are difficult to obtain, such as nonconvex optimization and inverse-problem applications, the solution of a well-trained L2O method can have better qualities than those of classic methods. Let us call the optimization method (either hand-engineered or trained by L2O) the \emph{optimizer} and call the optimization problem solvable by the method the \emph{optimizee}. Figure \ref{fig: optimizer-comp} compares classic optimizers and L2O optimizers and illustrates how they are applied to optimizees (yellow boxes).  

In many applications of optimization, the task is to perform a certain type of optimization over a specific distribution of data repeatedly. Each time, the input data that define the optimization are new but similar to the past. We say such an application has a narrow \emph{task distribution}. Conventional optimizers may be tuned for a task distribution, but the underlying methods are design for a theory-specified class of optimization problems. 
We often describe a conventional optimizer by the formulation (and its math properties), not the task distribution. For example, we say an optimizer can minimize a \emph{smooth-and-convex} objective function subject to \emph{linear} constraints. In L2O, the training process shapes the optimizer according to both the formulation \textit{and} the task distribution. When the distribution is concentrated, the learned optimizer can ``overfit'' to the tasks and may discover ``short cuts'' that classic optimizers do not take.

\begin{figure}[htb] 
\center
    \centering 
    \includegraphics[width=0.49\linewidth]{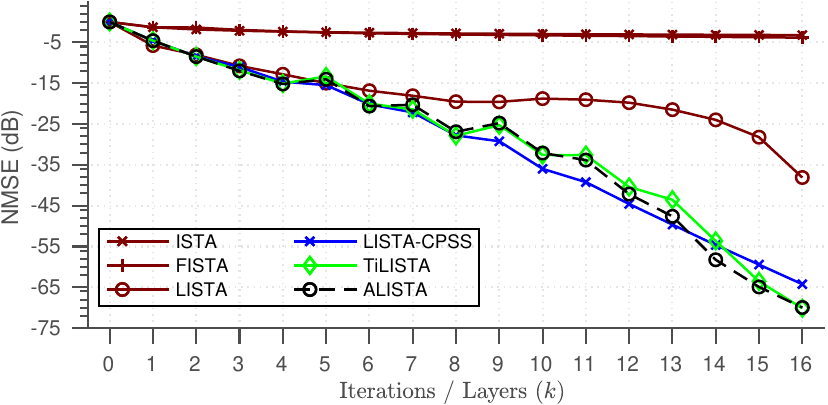}
    \includegraphics[width=0.49\linewidth]{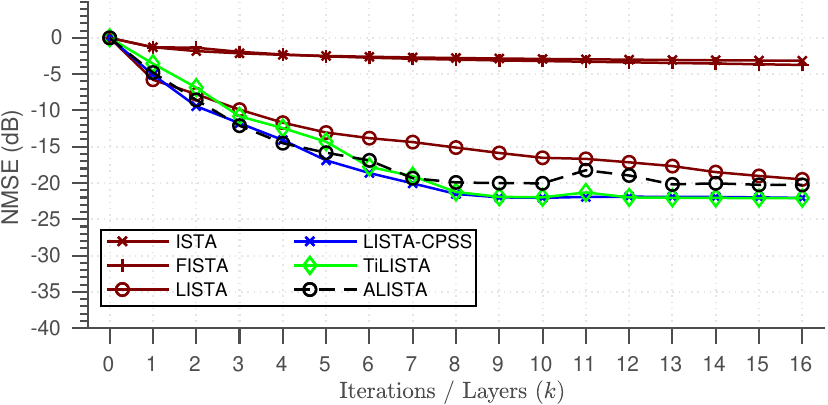}
    \vspace{-1em}
    \caption{\footnotesize L2O for compressive sensing. In our prior work \cite{Chen_Liu_Wang_Yin_2018,LiuChenWangYin2019_alista}, we successfully demonstrated that L2O optimizers (LISTA, LISTA-CPSS, TiLISTA and ALISTA) converge much faster than the two popular iterative solvers (ISTA, FISTA). Left: noiseless case. Right: noisy case (SNR = 20). X-axis is the number of iterations; Y-axis is normalized mean squared error (lower is better). \textit{See more details from the image source: Figure \ref{fig: optimizer-comp} of  \cite{LiuChenWangYin2019_alista}}.
    }
    \label{fig:lista}
  \end{figure}

L2O aims to generate optimizers with the following strengths: 
\begin{enumerate}
  \item An optimizer learned by L2O is expected to complete a set of optimizees from the same task distribution at a much faster speed than classic methods. In particular, such an L2O method can run even faster than so-called optimal methods such as Nesterov’s faster gradient descent (FGD) \cite{nesterov2005smooth} on a set of problems well suited for these optimal methods.
  \item The learned optimizer may also return a higher-quality solution to a difficult task than classic methods given a similar amount of computing budget. For instance, an L2O method can recover a sparse signal that is much more faithful (meanwhile much faster) than LASSO\footnote{In this example, the task is to recover a sparse signal, not necessarily solving the LASSO minimization.}.
\end{enumerate}
Therefore, when it comes to a small set of optimization tasks, L2O presents a potential approach to break the limits of analytic methods.

L2O has started to demonstrate great potential in some areas of optimization and applications. Examples 
include  convex $\ell_1$-minimization \cite{gregor2010learning}, neural network training \cite{andrychowicz2016learning}, black-box optimization \cite{chen2017learning} and combinatorial optimization \cite{khalil2017learning}. For example, a class of L2O methods, in the form of algorithm unrolling (which we review below), 
has state-of-the-art results in compressed sensing and inverse problems. The L2O methods in~\cite{Chen_Liu_Wang_Yin_2018,LiuChenWangYin2019_alista} converge much faster than the classic ISTA/FISTA optimizers~(see \Cref{fig:lista}), using an order-of-magnitude fewer iterations on unseen optimizees from the same distribution to reach the same accuracy. Application areas benefiting from 
L2O methods include computer vision \cite{Wang_Liu_Chang_Ling_Yang_Huang_2016,zhang2018ista,Corbineau_Bertocchi_Chouzenoux_Prato_Pesquet_2019}, medical imaging \cite{liang2020deep,yin2021end}, signal processing and communication \cite{Borgerding_Schniter_Rangan_2017,balatsoukas2019deep}, policy learning \cite{marino2020iterative}, game theory \cite{vadori2020calibration,vadori2021consensus}, computational biology \cite{cao2019learning,chen2019rna}, and even software engineering \cite{agrawal2020better}. For example, the applications of algorithm unrolling, a popular form of model-based L2O, has been specifically reviewed in Section \ref{app}. As another example, in the practical application of domain adaptation, model-free L2O has demonstrated the ability to adapt a pre-trained classifier to a new domain at less computational resource costs,  achieving the same accuracy or generalization performance~\cite{li_2020_halo,chen2020automated}. 

\subsection{Preliminaries}
\label{pre}

L2O starts with an architecture of the optimizer with free parameters to learn, so one must prepare a set of training samples that represent the task distribution, as well as choosing a training approach. 
The training process (Figure \ref{fig: optimizer-comp} (b) left) executes the approach, which iteratively applies the current optimizer to the sample optimizees and uses the observed performance to update the parameters. Over the time, the optimizer adapts to the training samples. This process trains an optimizer to take a series of actions on the optimizees, analogous of classic learning where we train ML models to make predictions or decisions.



 

We now formalize our notation. Consider an optimization problem $\min_\x f(\x)$ where $\x \in \RR^d$. A classic optimizer often iteratively updates $\x$ based on a handcrafted rule. For example, the first-order gradient descent algorithm is to take an update at each iteration $t$ based on the local landscape information at the instantaneous point $\x_t$ : $\x_{t+1}=\x_t-\alpha \nabla f(\x_t)$, where $\alpha$ is the step size. 



L2O has lots of freedom to use the available information. 
The information $\z_t$ available at time $t$ can include the existing iterates $\x_0,\ldots,\x_t$ as well as their gradients $\nabla f(x_0),\ldots,\nabla f(\x_t)$, etc.
Then the L2O approach models an update rule by a function $g$ of $\z_t$: $\x_{t+1}=\x_t- g(\z_t, \phi)$, where the mapping function $g$ is parameterized by $\phi$. 

Finding an optimal update rule can be formulated mathematically as searching for a good $\phi$ over the parameter space of $g$. 
In order to find a desired $\phi$ associated with a fast optimizer, \cite{andrychowicz2016learning} proposed to minimize the weighted sum of the objective function $f(\x_t)$ over a time span (which is called the unrolling length) $T$ given by:
\begin{align}\label{eq:l2o_obj}
\min_{\phi}  \mE_{f \in \cT} \left[\sum\limits_{t=1}^T w_t f(\x_t)\right], \quad \text{with}\quad \x_{t+1}=\x_t- g(\z_t, \phi),~t=1,\dots,T-1,
\end{align}
where $w_1,\dots,w_T$ 
are the weights and $f$ represents an optimizee from an ensemble $\cT$ of optimizees that represent the target task distribution. Note that the parameter $\phi$ determines the objective value through determining the iterates $\x_t$. L2O solves the training problem \cref{eq:l2o_obj} 
for a desirable $\phi$ and the update rule $g(\z_t, \phi)$. 


In practice, the choices of $w_t$ vary by case and depend on empirical settings. For example, many L2O models for sparse coding unroll to a fixed length $T$ for all optimizees, and then only minimize the step-$T$ functional value \cite{zhou2018sc2net,heaton2020safeguarded}, i.e., $w_T = 1$ and $w_1=\cdots=w_{T-1} = 0$. The work in \cite{andrychowicz2016learning}, on the other hand, assigns nonzero $w_t$ values and report more efficient Backpropagation Through Time (BPTT).

In L2O, one specifies the architecture of optimizer, which consists of both learnable components $\phi$ and fixed components.
We address both of them in this paper. The update rule $g$ is often parameterized by multi-layer neural networks  or recurrent neural networks. In theory, neural networks are universal approximators, so L2O may discover completely new, optimal update rules without referring to any existing updates. Since this kind of L2O architecture does not need a model,
we refer to it as \textit{model-free} L2O. 

The shortcomings of model-free L2O include: lacking convergence guarantees and requiring a high number of training samples. 
On tasks where classic operations --- such as projection, normalization, and decaying step sizes --- are critical to good performance, 
model-free L2O either cannot achieve good performance or require largely many training problems to discover classic operations from the scratch. 
To avoid these shortcomings, we consider incorporating the existing methods as base or starting points for learning, which reduce the search to fewer parameters and a smaller space of algorithms. We call this alternative approach \textit{model-based} L2O.

  \subsection{Broader Contexts}
  \label{context}


 
 When the tasks are general machine learning tasks, e.g., determining the parameters of a prediction model by minimizing its training loss, L2O overlaps with meta-learning \cite{vilalta2002perspective,hospedales2020meta,chen2020mate}. 
It is worth noting that ``meta-learning'' in different communities has different meanings. Here, meta-learning refers to using a method to improve learning algorithm(s), also called ``learning to learn''~\cite{andrychowicz2016learning,chen2017learning}. 
Recent results in this line of research contributes a significant portion of L2O development \cite{li2017learning}. The goal of L2O captures two main aspects of meta learning: rapid learning within tasks, and transferable learning across many tasks from some same distribution. However, L2O is not entirely meta-learning, because it takes into account domain knowledge of optimization and applies to many non-learning optimization tasks such as solving inverse problems.
 
 
  
L2O shares many similarities with AutoML \cite{yao2018taking}. The term ``AutoML" broadly refers to the automation of any step(s) in the ML lifecycle.
Traditionally, AutoML research focuses on model selection, algorithm selection, and hyperparameter optimization. These methods accelerate the design iterations of many types of ML algorithms, such as random forests, gradient boosting, and neural networks. AutoML recently draws (back) the mainstream attention because of its significant success in enhancing the performance of deep learning \cite{elsken2018neural}. Among the topics under AutoML, those most relevant to L2O are algorithm configuration \cite{hutter2011sequential} and hyperparameter optimization (HPO) \cite{feurer2019hyperparameter}. Algorithm configuration determines a high-performing configuration of some algorithm across a given set of problem instances. HPO tunes a set of hyperparameters specifically for an ML algorithm, 
via Bayesian optimization \cite{klein2017fast} or classification-based optimization \cite{yu2016derivative}. \cite{thornton2013auto} has combined algorithm selection and hyperparameter optimization, also known as CASH. 
The main difference between these works and L2O lies in that L2O can discover \textit{new optimization methods} from a parameterized algorithmic space of optimizers, rather than only \textit{selecting} from or \textit{tuning} a few existing optimizers. Yet, the boundary between HPO and L2O is often blurry since certain L2O methods~\cite{LiuChenWangYin2019_alista,xu2019learning,wang2020guarantees} configure or predict hyperparameters for existing optimizers only.


L2O is closely related to the new frontier of \textit{learning-augmented algorithms} \cite{MIT}. Classic algorithms are designed with worst-case performance guarantees in mind and do not adapt to actual inputs. On the other hand, ML algorithms often achieve competitive performance by adapting to inputs, but their worst-case performance on (unseen) inputs degrades significantly. Learning-augmented algorithms combine the best of both worlds, using ML to improve the performance of classic algorithms, by adapting their behaviors to the input distribution. 
Examples include learning-augmented data structures \cite{kraska2018case,mitzenmacher2018model}, streaming and sketching algorithms \cite{jiang2019learning,hsu2019learning}, online algorithms \cite{foster2018online}, error-correcting codes \cite{kim2018communication,kim2020deepcode}, scheduling algorithms \cite{mitzenmacher2020scheduling}, approximation algorithms \cite{indyk2019learning}, and safeguarded learned algorithms \cite{heaton2020safeguarded}. L2O can be counted as a subset of those learning-augmented algorithms. 
A comprehensive list of relevant materials can be found in \cite{MIT}.


\vspace{-0.5em}
\paragraph{Previous review articles and differences}  There have been a number of review articles on meta learning \cite{vilalta2002perspective,hospedales2020meta} and AutoML \cite{yao2018taking,hutter2019automated,he2021automl}. 

The work in \cite{monga2019algorithm} surveys algorithm unrolling. 
Work \cite{shlezinger2020model} reviews 
model-based deep learning. 
These articles have overlapping scopes with ours as algorithm unrolling is a main (but not the only) technique of model-based L2O. This article features a more comprehensive coverage including both model-based and model-free L2O approaches.

Lastly, 
let us mention recent works that leverage ML for combinatorial and discrete optimization \cite{LiChenQifengKoltun2018_Combinatorial,DaiKhalilZhangDilkinaSong2018_learning,bertsimas2019online,bertsimas2021voice,cauligi2020learning}. \cite{BengioLodiProuvost2018_Machine} provides a comprehensive review on ML for combinatorial optimization.

\subsection{Paper Scope and Organization}

      We draw distinctions between model-free and model based L2O approaches and review many L2O techniques. 
      Emphases are given to recurrent network-based L2O methods, algorithm unrolling, and plug-and-play. We discuss how to train them effectively and how their designs can benefit from the structures of the optimizees and classic optimization methods. 
      
    We benchmark existing model-based and model-based L2O approaches on a few  representative optimization problems. We have released our software implementation and test cases as the \textbf{Open-L2O} package at (\url{https://github.com/VITA-Group/Open-L2O}), in  the  spirit of reproducible research and fair benchmarking.

The rest of the article is organized as follows. Section 2 and Section 3 review and categorize the existing works in model-free and model-based L2O, respectively. Section 4 describes the Open-L2O testbeds, our comparison experiments, and our result analyses. The article is concluded by Section 5 with a discussion on research gaps and future work.

\section{Model-Free L2O Approaches}

A model-free L2O approach in general aims to learn a parameterized update rule of optimization  
without taking the form of any analytic update  
(except the updates being iterative). The recent mainstream works in this vein \cite{andrychowicz2016learning,chen2017learning,lv2017learning,wichrowska2017learned,cao2019learning} leverage recurrent neural networks (RNNs), most of which use the long short-term memory (LSTM) architecture. An LSTM is unrolled to perform iterative updates and trained to find short optimization trajectories. One set of parameters are shared across all the unrolled steps. At each step, the LSTM takes as input the optimizee's local states (such as zero-order and first-order information) and returns the next iterate. 


Model-free L2O shares the same high-level workflow depicted in Figure \ref{fig: optimizer-comp} (b).
It is divided into two stages: an offline \textit{L2O training} stage that learns the optimizer with a set of optimizees sampled from the task distribution $\mathcal{T}$; then, an online \textit{L2O testing} stage that applies the learned optimizer to new optimizees, assuming they are samples from the same task distribution. Table~\ref{tab:summary} compares a few recent, representative works in this area, from multiple lens including what optimizer architecture to select, what objectives and metrics that meta training and meta testing each use, and what input feature and other techniques are adopted. There are a great variety. 


\begin{table}[ht]
\caption{Summary and comparison of representative model-free L2O methods.}
\label{tab:summary}
\centering
\resizebox{1\textwidth}{!}{
\begin{tabular}{l|c|c|c|c|c}
\toprule
Paper & Optimizer Architecture & Input Feature & Meta Training Objective & Additional Technique & Evaluation Metric \\ \midrule
\cite{andrychowicz2016learning} & LSTM & Gradient & Meta Loss  & \begin{tabular}[c]{@{}c@{}} Transform input gradient $\nabla$ \\ into $\mathrm{log(\nabla)}$ and $\mathrm{sign(\nabla)}$\end{tabular} & Training Loss\\
\midrule
\cite{chen2017learning} & LSTM & Objective Value & Objective Value  & N/A & Objective Value \\
\midrule
\cite{lv2017learning} & LSTM & Gradient & Meta Loss  & \begin{tabular}[c]{@{}c@{}}Random Scaling\\ Combination with Convex Functions\end{tabular}  & Training Loss\\
\midrule
\cite{wichrowska2017learned} & Hierarchical RNNs & \begin{tabular}[c]{@{}c@{}}Scaled averaged gradients, \\relative log gradient magnitudes, \\relative log learning rate \end{tabular} & Log Meta Loss  & \begin{tabular}[c]{@{}c@{}}Gradient History Attention \\ Nesterov Momentum \end{tabular}  & Training Loss\\
\midrule
\cite{pmlr-v97-metz19a} & MLP & Gradient & Meta Loss  & Unbiased Gradient Estimators  & \begin{tabular}[c]{@{}c@{}} Training Loss \\ Testing Loss \end{tabular}\\
\midrule
\cite{li2016learning} & RNN Controller & Loss, Gradient & Meta Loss & Coordinate Groups  & Training Loss\\
\midrule
\cite{pmlr-v70-bello17a} & \begin{tabular}[c]{@{}c@{}} Searched Mathematical Rule \\ by Primitive Functions \end{tabular} & Scaled averaged gradients & Meta Loss & N/A & Testing Accuracy \\
\midrule
\cite{cao2019learning} & Multiple LSTMs & \begin{tabular}[c]{@{}c@{}} Gradient, momentum, \\ particle's velocity and attraction \end{tabular} & \begin{tabular}[c]{@{}c@{}} Meta Loss and \\ Entropy Regularizer
\end{tabular} & Sample- and Feature- Attention & Training Loss\\
\midrule
\cite{jiang2018learning} & RNN & Input Images, Input Gradient & Meta Loss  & N/A & \begin{tabular}[c]{@{}c@{}} Standard and Robust \\ Test Accuracies \end{tabular}\\
\midrule
\cite{xiong2020improved} & LSTM & Input Gradient & Meta Loss & N/A & \begin{tabular}[c]{@{}c@{}} Training Loss and \\ Robust Test Accuracy \end{tabular}\\
\bottomrule
\end{tabular}}
\end{table}

\subsection{LSTM Optimizer for Continuous Minimization: Basic Idea and Variants} \label{sec:model_free_basics}

As the optimization process can be regarded as a trajectory of iterative updates, RNNs are one natural option with a good inductive bias to learn the update rule. The first pioneering model-free L2O method \cite{andrychowicz2016learning}  proposed to model the update rule implicitly by gradient descent. This first LSTM method for L2O has often been referred to as \underline{L2O-DM} in literature\footnote{Methods that are underscored in this section are also later evaluated and benchmarked in Section.~\ref{sec:benchmark}.}.

At high-level, the learned optimizer (modeled by LSTM) is updated by the gradients from minimizing the loss induced by the optimizee, and optimizees are updated by the update rule predicted by the optimizer:
\begin{equation}
    \mathcal{L}(\phi)=\mathbb{E}_{(\theta_0, f) \in \cT}\left[\sum_{t=1}^{\mathtt{T}} w_{t} 
    f\left(\theta_{t}\right)\right] \quad \mathrm { where } \quad 
    \begin{array}{r}\theta_{t+1}=\theta_{t}+g_{t} \\
    {\left[\begin{array}{l}g_{t} \\
    h_{t+1}\end{array}\right]=m\left(\nabla_{t}, h_{t}, \phi\right)}\end{array}.
    \label{eq:andrychowicz2016learning}
\end{equation}
In Eq.~\ref{eq:andrychowicz2016learning}, $f$ represents an optimizee (e.g. a neural network with its loss function), with $\theta_0$ denoting the initialization of $f$. A set of optimizees and their initializations are sampled from the task distribution $\cT$. 
$\phi$ is the parameters of the L2O optimizer. $\nabla_{t}=\nabla_{\theta} f\left(\theta_{t}\right)$ is the gradient of the objective function with respect to the optimizee's parameters. $m$ is the LSTM optimizer, and $g$ and $h$ are update and hidden state produced by $m$, respectively. $\mathtt{T}$ is the maximal unrolling length for LTSM, often set due to the memory limitation; and $w_t$ = 1 was used here for every $t$.

Two ad-hoc training techniques were utilized here. Firstly, each optimization variable (or called coordinate) of the optimizee's parameter $\theta_t$ shared the same LSTM weights, but with different hidden states, in order to reduce the memory overhead faced when scaling up the number of optimization coordinates; this trick has been inherited by many LSTM-type L2O works. Secondly, to stabilize the optimizer's training, the authors proposed to reduce the dynamic range of the optimizee's gradient magnitudes, by preprocessing $\nabla_{t}$ into $(\mathrm{log}(\nabla_{t}), \mathrm{sgn}(\nabla_{t}))$ as the input into the optimizer.



 The authors of \cite{andrychowicz2016learning} conducted a few proof-of-concept studies on small-scale tasks such as MNIST classification, where they showed the learned optimizer $m$ can converge faster than some stochastic gradient descent based optimizers such as SGD, RMSprop, and Adam. Despite this initial success, the scalability and generalizability of \cite{andrychowicz2016learning} remain underwhelming. Two specific questionable shortcomings are:
\begin{itemize}
    \item \textbf{Generalizability of learned optimizer to unseen and potentially more complicated optimizees.} The L2O training set contains sampled optimizees from the target distribution. Just like typical machine learning models, we might expect a learned optimizer to both interpolate and extrapolate from the seen (meta-)training set, the latter being more challenging. Taking training deep networks for example, during meta-testing, we might expect the learned optimizer to generalize (extrapolate) to the training of deeper or wider networks beyond instances seen in the meta training set. This generalizability is demanded also due to the memory bottleneck during meta training. To update the L2O optimizer $m(\phi)$ using back-propagation, we need to keep the gradients and the computation graph of the optimizee in memory. Therefore, the memory bottleneck arises when the unroll length of the $m(\phi)$ becomes large, or the dimension of optimizee's parameter $\theta$ is high. Other types of generalization, such as training a network with unseen different activation functions or loss functions, are also found to be desirable yet challenging \cite{andrychowicz2016learning}.
    
    Moreover, applying the L2O optimizers to train more complex neural networks may meet more sophisticated loss landscapes, which cannot be well learned from observing simple training dynamics from small networks. For example, it is well known that training over-parameterized networks will exhibit sharply distinct behaviors at small and large learning rates, and the two regimes usually co-exist in practical training and are separated by a phase transition \cite{leclerc2020two,lewkowycz2020large,wu2020direction}.

    \item \textbf{Generalizability of learned optimizer to longer training iterations.} Training larger models naturally takes  more iterations, which calls for longer-term L2O modeling. An LSTM model can in principle characterize longer-term dependency by unfolding more time steps (i.e., longer unrolling length $\mathtt{T}$). However, that often results in L2O training instability due to the optimization artifacts of LSTM such as gradient explosion or vanishing, in addition to the practical memory bottleneck. Therefore, most LSTM-based L2O methods \cite{andrychowicz2016learning,chen2017learning,lv2017learning,wichrowska2017learned,pmlr-v97-metz19a,cao2019learning,chen2019rna,chen2020learning2stop} are forced to limit their maximal unroll lengths and to truncate the unrolled optimization (e.g., up to 20 steps). As a result, the entire optimization trajectory is divided into consecutive shorter pieces, where each piece is optimized by applying a truncated LSTM. 
    
    However, choosing the unrolling (or division) length faces a well-known dilemma \cite{lv2017learning}: on one hand, a short-truncated LSTM can result in premature termination of the iterative solution. Naively unfolding an LSTM-based L2O optimizer (trained with a small $\mathtt{T}$) to more time stamps at meta-testing usually yields unsatisfactory results. The resulting ``truncation bias" causes learned optimizers to exhibit instability and yield poor-quality solutions when applied to training optimizees.

    
\end{itemize}
The above research gaps motivate several lines of follow-up works, as summarized below:

\paragraph{Debiasing LSTM Truncation} As explained above, selecting the unrolling length $\mathtt{T}$ reflects a fundamental dilemma in LSTM-type L2O models: a longer $\mathtt{T}$ causes optimization difficulty at L2O training, while a smaller $\mathtt{T}$ incurs larger ``truncation bias" and hampers the generalization at meta-testing.

\begin{figure}[!ht]
    \centering
    \includegraphics[width=0.99\textwidth]{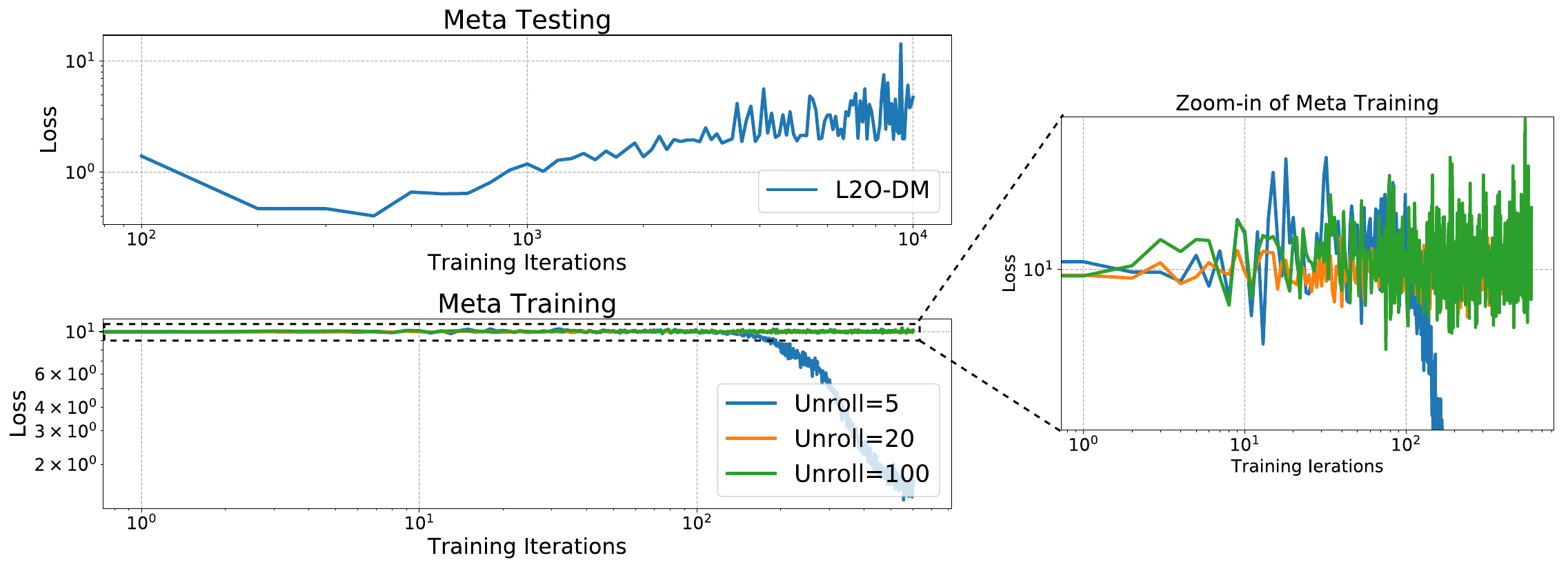}
    \caption{\small{ The dilemma of longer unrolling in both L2O testing and L2O training: (1) \textit{generalization failure} at L2O testing: the upper figure shows the meta-testing loss of the vanilla L2O method \cite{andrychowicz2016learning}} to quickly diverge, as we increase the number of steps at the meta-testing. This failure case was also observed in Figure \ref{fig:lista} of \cite{wichrowska2017learned} and in Figure \ref{fig:lista} of \cite{lv2017learning}. (2) \textit{optimization failure} at L2O training: the bottom figure shows the L2O training loss to also not decrease when it adopts longer unrolling lengths, due to the optimization artifacts of LSTM such as gradient explosion or vanishing.}
    \label{fig:meta_fail}
\end{figure}

Many LSTM-type works focus on addressing this truncation bias. \cite{wichrowska2017learned} drew the unroll length of their truncated optimization from a heavy tailed distribution. \cite{pmlr-v97-metz19a} proposed to replace LSTM with just MLP for the optimizer. The authors of \cite{pmlr-v97-metz19a} proposed to smooth the loss landscape by using two unbiased gradient estimators with dynamic re-weighting:
\begin{align}
    g_{\mathrm{rp}} &= \frac{1}{S}\sum_{s=1}^S\nabla_\theta f(\theta+n_s), \quad n_1,\cdots,n_S \sim N(0, \sigma^2I) \quad \mathrm{i.i.d.} \\
    g_{\mathrm{es}} &= \frac{1}{S}\sum_{s=1}^S f(\Tilde{\theta}_s)\nabla_{\theta}[\mathrm{log}(\mathds{P}(\Tilde{\theta}_s))],  \quad \Tilde{\theta}_1,\cdots,\Tilde{\theta}_S \sim N(\theta,\sigma^2I) \quad \mathrm{i.i.d.}\\
    g_{\mathrm{merged}} &= \frac{g_\mathrm{rp}\sigma_\mathrm{rp}^{-2}+g_\mathrm{es}\sigma_\mathrm{es}^{-2}}{\sigma_\mathrm{rp}^{-2}+\sigma_\mathrm{es}^{-2}}
\end{align}
$g_\mathrm{rp}$ and $g_\mathrm{es}$ are graident estimated by ``reparameterization trick'' \cite{kingma2013auto} and ``evolutionary strategies'' \cite{wierstra2008natural}, and $\sigma_\mathrm{rp}$ and $\sigma_\mathrm{es}$ are empirical estimates of the variances of $g_\mathrm{rp}$ and $g_\mathrm{es}$, respectively.
In this way, they demonstrated L2O can train convolutional network networks faster in wall-clock time compared to tuned ﬁrst-order methods, with reduced test losses. Their later work \cite{metz2019using} also found that the meta-learned optimizer can train image classification models such that they are robust to unseen image corruptions.

\paragraph{Stronger LSTM Architectures}
Another direction explored is to introduce stronger RNN/LSTM architecture to L2O. Instead of using a single RNN layer, Wichrowska et al. \cite{wichrowska2017learned} leveraged three RNN layers to learn the optimization in different scales. They organized three RNN layers in a hierachical fashion, here we simply denote them as ``bottom RNN'', ``middle RNN'', ``upper RNN''. The ``bottom RNN'' directly takes the scaled gradients from the optimizee as input. Given a specific training step, the ``middle RNN'' receives the average all hidden states from ``bottom RNN'' across the optimizee's parameter coordinates, and outputs a bias term to the ``bottom RNN''. Further, the ``upper RNN'' takes the averaged hidden states from ``middle RNN'' across a certain window of optimizee's training steps. By using smaller hidden states in ``bottom RNN'', this hierarchical design of L2O optimizer achieved lower memory and compute overhead, while achieving better generalization. We refer \cite{wichrowska2017learned} as \underline{L2O-Scale} in this article.

\paragraph{Improved L2O Training Techniques}
In order to improve the generalization of the learned optimizer to both longer unrolling (i.e. longer optimization iterations) and unseen functions, \cite{lv2017learning} proposed two training tricks. The first one, called \textit{random scaling}, could be viewed as a special ``data augmentation" for L2O: a coordinate-wise scale factor $\mathbf{c}$ was randomly generated at each iteration to scale the parameters of the optimizee during L2O training: $f_{\mathbf{c}}(\theta)=f(\mathbf{c} \theta)$. It was motivated by the observation that in many optimization problems such as the quadratic function $f(\theta)=\lambda\|\theta\|_{2}^{2}$, the ideal update rule should achieve the same minima under varying $\lambda$. The second trick was to add a convex term during L2O training, as inspired by the proximal algorithms \cite{parikh2014proximal}; and avoided large random updates when the L2O optimizer is under-trained.
We refer \cite{lv2017learning} as \underline{L2O-RNNprop} in this article.


Lately, the authors of \cite{chen2020training} took a deeper dive into improved training techniques for L2O models. The authors first presented a progressive training scheme, which gradually increased the optimizer unroll length to mitigate the L2O dilemma of truncation bias (shorter unrolling) versus gradient explosion (longer unrolling). Furthermore, they presented an off-policy imitation learning approach to guide the L2O training, by forcing the L2O optimizer to mimic the update rules generated by analytic optimizers. The authors of \cite{chen2020training} evaluated their improved training techniques with a variety of state-of-the-art L2O models \cite{andrychowicz2016learning,lv2017learning,wichrowska2017learned}, and achieved boosted performance (lower training losses of unseen optimizees) without changing the original L2O RNN architecture in each method.
We refer \cite{chen2020training} as \underline{L2O-enhanced} in this article.


\subsection{Other Common Implementations for Mode-Free L2O}

While LSTM is so far the mainstream model, other optimizer models have also been explored. We describe two alternatives: \textit{reinforcement learning} (RL), and \textit{neural symbolics}. 

\cite{li2016learning} proposed to learn an RL policy $\pi$ to predict the update rule, as the learned optimizer. The policy $\pi$ samples the update steps from a Gaussian distribution. The mean and variance of the Gaussian distribution are learnable parameters of the L2O policy, updated by reinforcement learning. The observation of the policy is composed of the current value of objective function (i.e. loss), the recent gradients, and changes of the objective values up to $H$ steps ($H = 25$ in their experiments). The policy receives the decrease of training loss as the reward.
Logistic regression functions and two-layer neural nets were leveraged as the testbeds. Further on, \cite{li2017learning} proposed to group coordinates under permutation invariance (e.g., weight matrix or a bias vector) into a coordinate group. This formulation reduced the computation cost of expensive optimization problems, making the proposed method extensible to wider neural networks. Generalizability of the learned policy was also studied by training the policy on shallow and narrow networks, and test on wider layers.
\cite{almeida2021generalizable} learned to update optimizer
hyperparameters instead of model parameters, directly using novel features, actions, and a reward function to feed RL. They demonstrated promising scalability to large-scale real problems. 

It is worthy to mention a special L2O work \cite{pmlr-v70-bello17a}, that explored L2O from the neural symbolic prospective. The authors also leveraged an RL controller, but avoided directly modeling the update rules. Instead, they designed a search space of operands (gradient $g$, gradient's running exponential moving average $\hat{m}$, etc.), unary functions of input $x$ ($e^x$, $\mathrm{log}(x)$, $\sqrt{|x|}$, etc.), and binary functions (mapping $(x, y)$ to $x + y$, $x - y$, $x * y$, etc.). The RL controller was learned to select a sequence of elements from the search space, formulate a function to process the input, and output the update rule.
Their searched best L2O optimizer (the RL controller) shows strong transferability, and improved training performance (lower training losses) on different tasks and architectures, including ImageNet classification and Google's neural machine translation system.
Their idea was further developed in \cite{real2020automl}, to automatically discover complete machine learning algorithms from raw  mathematical operations as building blocks, which concerns a more general problem than L2O.

\subsection{More Optimization Tasks for Model-Free L2O} \label{sec:model_free_more}
\paragraph{Black-box Optimization} \cite{chen2017learning} pioneered to extend the LSTM L2O framework \cite{andrychowicz2016learning} to derivative-free or black-box function optimization. Due to the absence of the optimizee's gradients as input, the authors of \cite{chen2017learning} instead treated the optimizee's input-output pair as the observation of the optimizer, and formulated the optimization as an exploration-exploitation trade-off problem. They updated the optimizer's hidden state $\mathbf{h}_{t}$ by the observation from the last step ($\mathbf{x}_{t-1}, y_{t-1}$) and then chose a new query point $\mathbf{x}_{t}$ to explore
\begin{align}
    \mathbf{h}_{t}, \mathbf{x}_{t} &=\operatorname{RNN}_{\theta}\left(\mathbf{h}_{t-1}, \mathbf{x}_{t-1}, y_{t-1}\right)\\
    y_{t} &\sim p\left(y \mid \mathbf{x}_{t}, \cdots, \mathbf{x}_{1}\right),
    \label{eq:chen2017learning}
\end{align}
where the function value $y_t$ was incrementally sampled from a Gaussian Process ($p$) at each query point $\mathbf{x}_{t}$ in the experiments. During L2O training, it was assumed that the derivatives of function value $y_t$ can be computed with respect to the input $\mathbf{x}_{t}$, which means the errors will be backpropagated for L2O training, but not needed at L2O testing time. \cite{chen2017learning} demonstrated that their learned RNN optimizers are competitive with state-of-the-art Bayesian optimization packages (Spearmint \cite{snoek2014input}, SMAC \cite{hutter2011sequential}, and TPE \cite{bergstra2011algorithms}).


\paragraph{Particle Swarm Optimization}
Current L2O methods mostly learn in the space of continuous optimization algorithms that are point-based and uncertainty unaware. Inspired by population-based algorithms (e.g. swarm optimization), \cite{cao2019learning} estimated the posterior directly over the global optimum and used an uncertainty measure to help guide the learning process. The authors designed a novel architecture where a population of LSTMs jointly learned iterative update formula for a population of samples (or a swarm of particles). The model can take as input both point-based input features, such as gradient momentum; and population-based features, such as particle's velocity and attraction from swarm algorithms. To balance exploration and exploitation in search, the authors directly estimated the posterior over the optimum
and included in the meta-loss function the differential entropy of the posterior. Furthermore, they learn feature- and sample-level importance reweighting (often called ``attention" in deep learning) 
in the L2O model, for more interpretable learned optimization rules. Their empirical results over non-convex test functions and the protein-docking application demonstrated that this new L2O largely outperforms the off-the-shelf Particle Swarm Optimization (PSO) algorithm \cite{moal2010swarmdock} and the vanilla LSTM-based L2O methods that are not uncertainty-aware \cite{andrychowicz2016learning}. We refer \cite{cao2019learning} as \underline{L2O-Swarm} in this article.

\paragraph{Minimax Optimization} One more challenging testbed for model-free L2O is to solve continuous minimax optimization, that is of extensive practical interest \cite{goodfellow2014generative,madry2018towards,wu2018towards}, yet notoriously unstable and difficult. Three prior works \cite{jiang2018learning,ruan2019learning,xiong2020improved} tried to plug in L2O into a specific application of minimax optimization called adversarial training \cite{madry2018towards}:
\begin{equation}
    \min _{\theta} \mathbb{E}_{(\boldsymbol{x}, \boldsymbol{y}) \sim D}\left\{\max _{\boldsymbol{x}^{\prime} \in \mathbb{B}(\boldsymbol{x}, \epsilon)} \mathcal{L}\left(f\left(\boldsymbol{x}^{\prime}\right), \boldsymbol{y}\right)\right\},
\end{equation}
where $D$ is the empirical distribution of input data, and the inner maximization is defined as the worst-case loss within a small neighborhood $\mathbb{B}(\bm{x}, \epsilon)$. In both works, the L2O only predicted update directions for the inner maximization problem, while the outer minimization was still solved by classic optimizer. They empirically showed that L2O can improve the solution quality to the inner maximization optimization, hence also leading to a better minimax solution and a more robustly trained model.

A latest work \cite{minimax} extended L2O to solve minimax optimization from end to end, for the first time. The authors proposed Twin L2O, consisting of two LSTMs for updating min and max variables, respectively. This decoupled design was shown by ablation experiments to facilitate learning,
particularly when the min and max variables are highly non-symmetric. Several enhanced variants and training techniques were also discussed. The authored benchmarked their L2O algorithm
on several relatively basic and low-dimensional test problems, and on which L2O compared favorably against state-of-the-art minimax solvers. How to scale up L2O for fully solving minimax problems of practical interest, such as adversarial training or GANs, remained to be an open challenge. 

\paragraph{Game Theory} RL-based model-free L2O has very recently found interest from the game theory field. \cite{vadori2020calibration} proposed to train multi-agent systems (MAS) to achieve symmetric pure Nash equilibria. Such equilibria needs to satisfy certain constraints so that MAS are calibrated for practical use. The authors adopted a novel dual-RL-based algorithm to fit emergent behaviors of agents in a shared equilibrium to external targets. They used parameter sharing with decentralized execution, to train multiple agents using a single policy network, while each agent can be conditioned on agent-specific information. The methodology shared similarities to \cite{cao2019learning}. Another work \cite{vadori2021consensus} extended consensus optimization to the constrained case. They introduced a new framework for online learning in non zero-sum games, where the update rule’s gradient and Hessian coefficients along a trajectory are learned by an RL policy, conditioned on the game signature. The authors mentioned that the same problem might potentially be solved by using the LSTM-based approach too \cite{minimax}.

\paragraph{Few-shot Learning} Another application of the LSTM L2O framework \cite{andrychowicz2016learning} explores the application of model-free L2O in the small data regime.
\cite{ravi2016optimization} first adopted a model-free L2O approach to learn a few-shot classifier, accessing only very few labeled examples per class at training.
The authors of \cite{ravi2016optimization} proposed to learn an LSTM-based meta-learner to optimize each task-specific classifier using the cell states in the LSTM, inspired by their observation that the gradient descent update on the parameters in the classifier resembles the cell state update in an LSTM. \cite{li2017meta} simplified \cite{ravi2016optimization} by constraining the optimization on the learner classifier to one step of gradient descent but with learnable initialization and step size.

\section{Model-Based L2O Approaches} 

%


We now overview model-based L2O approaches, which mark the recent movement to fuse traditional model-based optimization algorithms with powerful deep learning architectures \cite{monga2019algorithm}. Rather than use general-purpose LSTMs, these methods model their iterative update rules through a learnable architecture that is inspired by analytic optimization algorithms. Often, these learned methods can approximate problem solutions with tens of iterations whereas their classic counterparts make require hundreds or thousands of iterations \cite{LiuChenWangYin2019_alista}.

At a high level, model-based L2O can be viewed as a ``semi-parameterized" option that takes advantage of both model-based structures/priors and  data-driven learning capacity\footnote{We slightly abused the term of \textit{semi-parametric model} from statistical learning, to represent different meanings: here it refers to blending pre-defined structures/priors into black-box learning models.}. The growing popularity of this approach lies in its demonstrated effectiveness in developing compact, data-efficient, interpretable and high-performance architectures, when the underlying optimization model is assumed available or can be partially inferred. There is a large design space to flexibly balance between the two. Most model-based L2O methods take one of the two following mainstream approaches.

The first approach is known as \textbf{plug and play (PnP)}, and the key idea here is to \textit{plug} a previously trained neural network (NN) into part of the update for an optimization algorithm (i.e., in place of an analytic expression), and then \textit{play} by immediately applying the modified algorithm to problems of interest (without additional training).  
We illustrate this with the common alternating direction method of multipliers (ADMM) formulation of PnP (e.g., see  \cite{venkatakrishnan2013plug}). Here the underlying task of the original optimization algorithm is to minimize the sum of two functions $f$ and $g$, using successive proximal operations. The PnP formulation replaces the proximal for $g$ with an operator $H_\theta$ to obtain the iteration:
\begin{subequations}
    \begin{align}
        x^{k+1} & = H_\theta(y^k - u^k) \label{eq:pnp-admm-1} \\
        y^{k+1} & = \mathrm{prox}_{\alpha f}(x^{k+1}+u^k) \label{eq:pnp-admm-2} \\
        u^{k+1} & = u^k + x^{k+1} - y^{k+1}.\label{eq:pnp-admm-3}
    \end{align}
\end{subequations}

Before inserting $H_\theta$ into the iteration above, the parameters
$\theta$  
are learned independently as the solution to a training problem, i.e.,
\begin{equation}
    \label{eq:model-unroll-learn}
    \theta \in \argmin_{\tilde{\theta}} \mathcal{L}(\tilde{\theta}).
\end{equation}
The loss function $\mathcal{L}$ is designed by an independent goal (e.g., to learn a natural image denoising operator).
For example, one might model an image recovery problem by using total variation (TV) as a regularizer \cite{rudin1992nonlinear}. In the optimization scheme chosen (e.g., ADMM) to recover the image, one step of the process could be to perform a proximal operation with TV, which effectively acts as a denoiser. The PnP framework proposes replacement of the TV proximal with an existing denoiser (e.g., a neural network \cite{ryu2019plug} or BM3D \cite{dabov2007image}).

The second approach is known as \textbf{algorithm unrolling}, whichs unrolls a truncated optimization algorithm into the structure of a neural network \cite{monga2019algorithm}. Updates take the form
\begin{equation}
    \label{eq:model-unroll}
    x^{k+1} = T(x^k; \theta^k), \quad k = 0,1,2,\cdots,K.
\end{equation}
Upon establishing this form, the parameters are learned by an end-to-end approach:
\begin{equation}
    \label{eq:model-unroll-learn_2}
   \min_{\{\theta^k\}_{k=0}^{K-1}} \mathcal{L}\big(x^K(\{\theta^k\}_{k=0}^{K-1})\big),
\end{equation}
where $\mathcal{L}$ is the loss function we use in training. 
We emphasize the distinction that the parameters in unrolled schemes are trained end-to-end using the iterate $x^K$ as a function of each $\theta^k$ whereas training occurs separately for PnP.
We will introduce how to design $\mathcal{L}$ later in this section.
Below we discuss each of these approaches, their typical features, and which approach is suitable for various applications.

\begin{figure}
    \centering
    \includegraphics[width=5.5in]{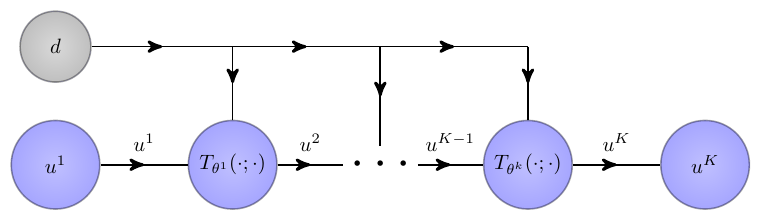}
    \caption{A common approach in various L2O schemes is to form feed forward  networks by unrolling an iterative algorithm, truncated to $K$ iterations, and tuning parameters at each layer/iteration $k$. This generalizes the update formula $u^{k+1} = T(u^k; d)$ to include a dependence on weights $\theta^k$, denoted by a subscript.}
    \label{fig:unrolling}
\end{figure}

\begin{table}[!ht]
    \centering
    \begin{tabular}{c|c|c|c}
         & Opt   & IP   & PnP  \\\hline
        \rowcolor{black!20!white}
        Tunable Optimization Model &  & \checkmark & \checkmark \\\hline
        Tunable Update Formulas & \checkmark & \checkmark & \checkmark  \\\hline
        \rowcolor{black!20!white}
        Training Loss tied to Model & \checkmark &  &  \\\hline
        Training Loss measure $u^\star$ error &  & \checkmark & \\\hline
        \rowcolor{black!20!white}
        Convergence guarantee to ideal $u^\star$ & S & S &  \\\hline
        Interpretable Updates & S & S & $\checkmark$\\\hline
    \end{tabular}
    \caption{Comparison of properties for three types of model-based L2O methods: unrolled objective-based  (Opt), inverse problem  (IP), and Plug and Play (PnP) methods. Here $\checkmark$ and S mean that the corresponding property \textit{always} and \textit{sometimes} holds, respectively.}
    \label{tab:comparison-props-Opt-IP-PnP}
\end{table}

\subsection{Plug and Play}
PnP methods date back to 2013 when they were first proposed for ADMM \cite{venkatakrishnan2013plug}.
The authors replaced one of the proximal operators in ADMM with a denoiser (e.g., K-SVD \cite{aharon2006k}, BM3D \cite{dabov2007image}, and non-local means \cite{buades2005nonlocal}) and showed empirical performance surpassing the original ADMM algorithm.
Several other PnP schemes were also introduced and their empirical success was shown in the literature \cite{heide2014flexisp,metzler2015bm3d,rond2016poisson,brifman2016turning,wang2017parameter,ono2017primal,kamilov2017plug,he2018optimizing,gupta2018cnn,yang2018proximal,ye2018deep,lyu2019iterative,zhang2019deep,yuan2020plug,ahmad2020plug,mataev2019deepred,song2020new}. 
The theoretical guarantees of convergence of PnP-ADMM was first provided in \cite{sreehari2016plug} under the assumption that the derivative of the denoiser was doubly stochastic. Then \cite{chan2016plug} proved the convergence of PnP-ADMM with the use of a more realistic ``bounded denoisers" assumption.
Some other PnP schemes were analyzed under different assumptions
\cite{teodoro2017scene,teodoro2019image,buzzard2018plug,dong2018denoising,tirer2018image,chan2019performance,sun2019online,gavaskar2020plug,xu2020provable,sun2020scalable,sun2020async}.

PnP was initially not connected to L2O, until some concurrent works \cite{meinhardt2017learning,rick2017one,zhang2017learning} introduced the concept of “learning operators” into a PnP framework. Instead of using a manual-designed denoiser, they modeled the proximal operator as a deep neural network and learned it from data. The empirical performance of such an approach largely exceeded the prior PnP works. From an L2O perspective, the learned methods were able to improve performance by either accelerating execution of a subroutine in the optimization algorithm or providing a ``better'' solution than previously obtained.

Besides learning the operators in a PnP framework, one can also learn a functional as the regularizer. \cite{bigdeli2017deep,visapp18} modeled the statistical distribution of nature images as a Denoising Autoencoder (DAE) and constructed the optimization objective by maximizing the likelihood. The DAE is differentiable; thus, the inverse problem can be solved with gradient-based optimization algorithms.
A parameterized discriminator function defined in the Wasserstein distance between two distributions was used by \cite{lunz2018adversarial}  as a   learned a functional that discriminates between a ground-truth and fake images. The authors treated this learned functional as a regularizer and used it in a gradient descent algorithm for image inverse problems.

Learning is not only able to help find the denoiser/regularizer in PnP, but is also able to be used in a more meta manner to help find good parameters in PnP. \cite{wei2020tuning} learned a policy network to tune the parameters in Plug-and-Play with reinforcement learning. By using the learned policy, the guided optimization can reach comparable results to the ones using oracle parameters tuned via the inaccessible ground truth.

The marriage of PnP with L2O also provides theoretical blessings. 
While a manual-designed regularizer may not guarantee the convergence of PnP to the desired solution, a learning-based method provides the flexibility to meet the condition of convergence.  \cite{ryu2019plug} studied the convergence of some PnP frameworks using fixed-point iterations, showing guaranteed convergence under a certain Lipschitz condition on the denoisers. They then proposed a normalization method for training deep learning-based denoisers to satisfy the proposed Lipschitz condition. 
Similarly, \cite{terris2021enhanced} also provided   an approach for building convergent PnP algorithms using monotone operator theory  and constraining the Lipschitz constant of the denoiser during   training.
This fixed point framework was   extended by \cite{cohen2020regularization} to the RED-PRO framework, which  also showed convergence to global solutions. 
Independently of RED-PRO, \cite{heaton2020projecting}  proposed a method within the RED-PRO framework for learning projection operators onto compact manifolds of true data.
Utilizing assumptions about the manifold and sufficient representation by the sampled data, the authors proved their approach to constructing a PnP operator can provably approximate  the projection (in probability) onto a low dimensional manifold of true data.



\subsection{Algorithm Unrolling}
Herein we overview L2O methods comprised of unrolling iterative optimization algorithms. 
We start by emphasizing there are two distinct goals of unrolled L2O methods: either to minimize an objective, or to recover a signal.  The distinction is that the aim for inverse problems is to tune the parameters so that they minimize the reconstruction accuracy of some ``true'' signal rather than find a minimizer of the model's objective function.  
The practicality for objective minimization is to speed up convergence. 
Keeping these categorizations in mind will help provide the reader intuitive lenses for looking at different approaches. See Table \ref{tab:comparison-props-Opt-IP-PnP} for a comparison of common qualities in these methods and PnP. 

Below we provide a comprehensive and organized review of existing works, along multiple dimensions: what problems they solve, what algorithms they unroll, what goals they pursue (objective minimization or signal recovery), and to what extent they freely parameterize their learnable update rule.

\subsubsection{Different target problems}
\label{sec:model-based-problems}

We broadly categorize the main problems tackled by model-based L2O in four categories: probabilistic graph models, sparse and low rank regression, differential equations, and quadratic optimization.
\begin{itemize}
    \item \textbf{Sparse and low rank regression}. The most investigated problems in the literature of unrolling is probably the sparse regression, inspired by the first seminal work \citep{gregor2010learning}, which unrolls the Iterative Shrinkage Thresholding Algorithm (ISTA) or its block coordinate variant as a recurrent neural network to solve LASSO for fast sparse coding.
    The unrolling philosophy is  also used in low rank regression as it shares some common nature with sparse regression.
    \begin{itemize}
        \item \textbf{LASSO:} Most of unrolling works following \citep{gregor2010learning} also  solve LASSO-type optimization problems for sparse coding \citep{sprechmann2013supervised,Moreau_Bruna_2017,Perdios_Besson_Rossinelli_Thiran_2017,Giryes_Eldar_Bronstein_Sapiro_2018,zhou2018sc2net,ablin2019learning,Hara_Chen_Washio_Wazawa_Nagai_2019,cowen2019lsalsa,Wu_Dimakis_Sanghavi_Yu_Holtmann-Rice_Storcheus_Rostamizadeh_Kumar_2019}.
        Beyond naive LASSO, A natural extension is to apply the same \textit{unrolling and truncating} methodology to solve \textit{group LASSO} \citep{Bronstein_Sprechmann_Sapiro_2012}, i.e. finding solutions with structured sparsity constraints. \cite{papyan2017convolutional,aberdam2019multi} extended solving (convolutional) LASSO in a multi-layer setting, and an adapted version of ISTA was unrolled in \cite{sulam2019multi}.
        \item \textbf{Analysis model:} Different from the sparse prior in LASSO (also known as the \textit{Synthesis} formulation), which assumes the signal is the linear combinations of a small number of atoms in a dictionary, the sparse analysis model assumes the existence of a forward transformation that sparsifies the signal. \cite{cherkaoui2020learning} applied unrolling to regression problems with total variation regularizations, which follows the analysis formulation. Problems with analysis sparse priors were also tackled using unrolling in \cite{malezieux2021dictionary} in a bi-level optimization context.
        \item \textbf{Other sparse/anti-sparse regression problems:} \cite{wang2016learning} and \citep{xin2016maximal} unroll the iterative hard thresholding (IHT) algorithm that solves \textit{$\ell_0$ minimization} problem instead of LASSO. Unrolling for $\ell_\infty$ minimization was also considered by \cite{WangYangChangLingHuang2016_learning}, leading to so-called anti-sparse representation learning \cite{studer2014democratic}. 
        \item \textbf{Low Rank Regression:} \citep{HersheyRouxWeninger2014_deep,sprechmann2014supervised,SprechmannBronsteinSapiro2015_learning,yakar2013bilevel} all extend the unrolling idea to low-rank matrix factorization. Specifically, \citep{SprechmannBronsteinSapiro2015_learning} proposes tailored pursuit architectures for both robust PCA and non-negative matrix factorization with specialized algorithms inspired by the non-convex optimization techniques to alleviate the expensive SVD involved in each iteration.
    \end{itemize}
    \item \textbf{Probabilistic Graphical Model}. The unrolling method can  also  be adopted to solve probabilistic graphical models. For example, \citep{HersheyRouxWeninger2014_deep} interpret conventional networks as mean-field inference in Markov random fields and obtains new architectures by modeling the belief propagation algorithm to solve the Markov random field problems. \citep{zheng2015conditional} formulate Conditional Random Fields with Gaussian pairwise potentials and mean-ﬁeld approximate inference as recurrent neural networks.
    \item \textbf{Differential equations}. Another line of works \citep{Chen_Pock_2017,long_pde-net_2018,long_pde-net_2019,greenfeld_learning_2019} unroll the evolution in \textit{Partial Differential Equation} (PDE) systems. While PDEs are commonly derived based on empirical observations. those recent advances offer new opportunities for data-driven discovery of (time-dependent) PDEs from observed dynamic data.  One can train feed-forward or recurrent neural networks to approximate PDEs, with applications such as fluid simulation \cite{wiewel2019latent}. Loosely related is also significant work on connecting neural networks with \textit{Ordinary Differential Equation} (ODE) systems.
    \item \textbf{Quadratic optimization}. Some recent works investigate the unrolling method in quadratic optimization problems \citep{wang2020guarantees,Chen_Zhang_Reisinger_Song_2020}, which are easier to solve compared to the problems studied above. The focus here is more on the theoretical analysis of the convergence, and/or the interplay between the unrolled algorithm and the resultant deep model's property (e.g., generalization and stability).
\end{itemize}

In some scenarios, we are not directly interested in the output of the unrolled model but use it for downstream tasks, e.g. \text{clustering} \citep{Wang_Chang_Zhou_Wang_Huang_2016} and \text{classification} \citep{wang2016learning,zhou2018sc2net,cowen2019lsalsa}. That will often lead to task-driven joint optimization and the end task output becomes the focus of evaluation, in place of the original unrolling algorithm's output fidelity.

\paragraph{Inverse Problems}
In many cases, however, optimization problems with manually designed objective functions  only provide approximations to the original signals that we are really interested. This is often due to inexact prior knowledge about the original signals. For example, sparsity and total variation regularizations  only partially reflect the complexity of natural images (approximated sparsity and smoothness) which are hardly true in real-world applications and do not depict the exact characteristics of natural images.
Therefore, many works solve the \textit{inverse problem} directly, striving to recover the original signals that we are really interested in.
We return to this task in a later subsection.

\paragraph{Constrained Optimization Problems} Most target problems mentioned above in this subsection are unconstrained optimization problems. Some exceptions such as the sparsity-constrained and low-rank regression problems, e.g., LASSO, have equivalent unconstrained formulation under proper regularizations. Some others have easy-to-implement forms of projection onto the constraint sets, including the $\ell_{0/\infty}$-constrained regression \cite{wang2016learning,xin2016maximal,WangYangChangLingHuang2016_learning}, and non-negative matrix factorization with non-negative constraints \cite{SprechmannBronsteinSapiro2015_learning}.

There have been a few efforts directly tackling more general constrained optimization problems using unrolling. \cite{liu2019frank} unrolled Frank-Wolfe algorithm to solve the structured regression with general $\ell_p$-norm constraint ($p\ge 1$), and proposed a novel closed-form nonlinear pooling unit parameterized by $p$ for the projection. 
\cite{pauwels2021hcgm} unrolled Frank-Wolfe algorithm for least square problems with affine, non-negative and $\ell_p$-norm constraints. It was also the first to apply the unrolled network to financial data processing.
\cite{Corbineau_Bertocchi_Chouzenoux_Prato_Pesquet_2019,Bertocchi_Chouzenoux_Corbineau_Pesquet_Prato_2020} investigated image restoration with various hard constraints and unrolled proximal interior point algorithms while incorporating the constraints using a logarithmic barrier.

\subsubsection{Different Algorithms Unrolled}

For the bulk of iterative optimization algorithms that have been unrolled, we classify them into three categories: forward backward splitting, primal-dual methods, and \textit{others}.
The first two categories consist entirely of first-order algorithms, which is due to their low computational complexity and the fact their resultant L2O methods are often more reliable to train. We also emphasize that although various works discussed below may use the same underlying algorithm as the base, they can vary greatly in their performance due to choices regarding what parameters are learned and what safeguard precautions are used to ensure convergence (discussed further in subsequent sections). We also emphasize that the majority of these algorithms revolve around obtaining some form of sparsity/low rank solution.\\
  
We begin with forward backward splitting (FBS).
The simplest learned scheme is \textit{Gradient Descent}, which is where the the gradient descent operation is applied (i.e., the forward operator) and the backward operator is simply the identity. This class of methods is studied in many works (e.g., see \citep{Giryes_Eldar_Bronstein_Sapiro_2018,gupta2018cnn,Diamond_Sitzmann_Heide_Wetzstein_2018,Takabe_Wadayama_2020,Chen_Zhang_Reisinger_Song_2020,adler2017solving,adler2018learned, domke2012generic,Putzky_Welling_2017,Chen_Pock_2017,lunz2018adversarial,mukherjee2020learned,Wu_Dimakis_Sanghavi_Yu_Holtmann-Rice_Storcheus_Rostamizadeh_Kumar_2019,Wadayama_Takabe_2019,kofler2020neural}).
However, the most popular focus in the literature is on problems with sparsity constraints, which are usually modeled by $\ell_1$/$\ell_0$-minimization. The former, also known as LASSO, can be solved by the \textit{iterative shrinkage-thresholding algorithm} (ISTA) \citep{blumensath2008iterative} and its variants.
This has yielded great interest in L2O schemes, as evidenced by the fact the works \citep{gregor2010learning,Wang_Chang_Zhou_Wang_Huang_2016,Bronstein_Sprechmann_Sapiro_2012,Wang_Liu_Chang_Ling_Yang_Huang_2016,zhang2018ista,Chen_Liu_Wang_Yin_2018,LiuChenWangYin2019_alista,takabe2020complex,ablin2019learning,Ito_Takabe_Wadayama_2019,Yao_Dang_Zhang_Wu_2019,Hara_Chen_Washio_Wazawa_Nagai_2019,Aberdam_Golts_Elad_2020,Behrens_Sauder_Jung_2020}, among others, provide various ways to unroll the original ISTA algorithm.
Additionally, \citep{Moreau_Bruna_2017,Perdios_Besson_Rossinelli_Thiran_2017,Tolooshams_Dey_Ba_2018,Aberdam_Golts_Elad_2020} unroll a Nesterov accelerated ISTA known as \textit{FISTA} (Fast ISTA) \citep{beck2009fast}. 
Continuing in the vein of sparsity,
\citep{Borgerding_Schniter_2016,Borgerding_Schniter_Rangan_2017,metzler2017learned} unroll another algorithm, called \textit{approximate message passing} ({AMP}), which introduces Onsager correction terms that whitens the noise in intermediate signals while \citep{he2019model,Ito_Takabe_Wadayama_2019,takabe2020complex} unroll \textit{Orthogonal AMP} \citep{ma2017orthogonal}, an extension to the original AMP.
For the $\ell_0$-minimization problems, the \textit{iterative hard-thresholding} ({IHT}) algorithm, which replaces the soft-thresholding function in ISTA with hard-thresholding, is unrolled in \citep{wang2016learning,xin2016maximal}.
Switching gears, for problems that involve minimize an objective that is the sum of several terms, the incremental proximal gradient method has also been unrolled \cite{kobler2017variational}.  
Further generalization of FBS is given by an abstract collection of flexible iterative modularization algorithms (FIMAs) \cite{liu2019convergence} that perform updates using the composition of two learned operators (provided the composition monotonically decreases an energy). A caveat of this particular approach is that the learned operator is also composed with a classic proximal-gradient operation.  \\

The next category of unrolled L2O algorithms consists of primal-dual schemes. These typically include variations of \textit{primal-dual hybrid gradient} (PDHG) and the \textit{alternating direction method of multipliers} (ADMM) (e.g., see \citep{sprechmann2013supervised,yang_Sun_Li_Xu_2016,rick2017one,adler2018learned,xie2019differentiable,cowen2019lsalsa, cheng2019model}). 
Another variation used is the \textit{Condat-Vu primal-dual hybrid gradient} \cite{jiu2020deep}. As above, many of these methods also focus on leveraging the sparse structure of data. 
\\


    

The remaining group of miscellaneous methods take various forms. These include another  first-order algorithm, \textit{Half Quadratic Splitting} (HQS), for image restoration \citep{zhang2017learning,yang2018proximal}. 
Using a substitution and soft constraint, the HQS approach solves a problem using a model that approximates variational problems (VPs).
Beyond first-order algorithms, unrolling is also applied to \textit{second-order (Newton or Quasi-Newton) methods} \citep{Xiong_De_la_Torre_2013,lilearning}, \textit{Differential Equations} \citep{long_pde-net_2018,Long_Lu_Dong_2019,Zhang_Lu_Liu_Dong_2018} and \textit{Interior Point} method \citep{Corbineau_Bertocchi_Chouzenoux_Prato_Pesquet_2019,Bertocchi_Chouzenoux_Corbineau_Pesquet_Prato_2020}, and \textit{Frank-Wolfe} algorithm \cite{liu2019frank,pauwels2021hcgm}. 
\textit{Conjugate Gradient} was proposed \cite{aggarwal2018modl} to help yield a more accurate enforcement of the data-consistency constraint at each iteration than comparable proximal-gradient methods. Lastly, \cite{khatib2020learned} propose an unfolded version of a greedy pursuit algorithm, i.e., \textit{orthogonal matching pursuit} (OMP), that directly targets at the original combinatorial sparse selection problem. Their methods called Learned Greedy Method (LGM) can accommodate a dynamic number of unfolded layers, and a stopping mechanism based on representation error, both adapted to the input. Besides, there are also L2O schemes that appear to be based entirely on heuristic combinations of classic optimization methods (e.g., see \cite{huang2020data}).
\ \\

\subsubsection{Objective-Based v.s. Inverse Problems} 
\label{sec:obj}

\paragraph{Objective Based.} 
The simplest L2O unrolling scheme is objective based.
Training is applied here to yield rapid convergence for a particular distribution of data.
The training loss can take various forms, including minimizing the expectation of an objective function \cite{ablin2019learning}, the objective's gradient norm, the  
distance to the optimal solution, or the fixed point residual of the algorithm \cite{heaton2020safeguarded} (e.g., if $T$ is the update operator, then $\|x-T(x)\|$ is the fixed point residual). 
Examples of objectives that have been extensively studied in the literature are presented in Section.~\ref{sec:model-based-problems}.

In addition, \emph{safeguarding} can be used in this situation, for guiding learned updates to ensure convergence  \citep{moeller2019controlling,heaton2020safeguarded}.
This can be accomplished in multiple ways. A typical approach is to generate a tentative update using an L2O scheme and then check whether the tentative update satisfies some form of descent inequality (e.g., yields a lesser energy value or fixed point residual). If descent is obtained, then the tentative update is used; otherwise, a classic optimization update is used. These safeguarded schemes provide the benefit of reducing computational costs via   L2O machinery while maintaining theoretical convergence guarantees.

\paragraph{Inverse Problems} 
Several L2O methods attempt to solve inverse problems (IPs) (e.g., see \cite{hammernik2017deep,huang2020data,jiu2020deep,kobler2020total,kofler2020neural,li2020nett,liu2019convergence,mukherjee2020learned}).
Here the task is to reconstruct a signal $x^\star$ from indirect noisy measurements.
The measurements are typically expressed by $d\in \mathbb{R}^m$ and are related to a forward operator $A :\mathbb{R}^n\rightarrow \mathbb{R}^m$ by
\begin{equation}
    d = A(x^\star) + \varepsilon,
\end{equation}
where $\varepsilon$ is noise. A few repeated themes arise the L2O IP literature.
The typical process is to i) set up a variational problem \eqref{eq:variatonal-problem} as a surrogate model, ii) choose a parameterized optimization algorithm that can solve\footnote{We mean to say that, for some choice of parameters, the algorithms solves \eqref{eq:variatonal-problem}.} \eqref{eq:variatonal-problem}, and iii) perform supervised learning to identify the optimal parameter settings. We expound upon these themes and their nuances below.

First, by creating a variational problem, one assumes $u^\star$ approximately solves
\begin{equation}
    \min_{x \in \mathbb{R}^n} \ell\left( d, A(x) \right) + J(x),
    \tag{VP}
    \label{eq:variatonal-problem}
\end{equation}
where $\ell:\mathbb{R}^m\times\mathbb{R}^m\rightarrow \mathbb{R}$ is a fidelity term that encourages the estimate $u$ to be consistent with the measurements $d$ and $J:\mathbb{R}^n\rightarrow \mathbb{R}$ is a regularizer. 
Learning comes into play since $u^\star$ is usually {\it not} the solution to (\ref{eq:variatonal-problem}), but instead some ``close'' variant.
Thus, upon choosing the form of $\ell$ and $J$,  one includes tunable parameters in the model and unrolls a parameterized optimization scheme for (typically) a fixed number $K$ of iterations. 
Algorithm updates can be parameterized beyond their classic form (e.g., replace a fixed matrix with a matrix that has tunable entries as done in \cite{gregor2010learning}). 
In most cases, the training loss takes the form of minimizing the expected value of the square of the Euclidean distance between the output estimate $u^K$ and the true signal $u^\star$, i.e., the parameters $\Theta$ are trained to  solve
\begin{equation}
    \min_{\Theta} \mathbb{E}_{d\sim\mathcal{D}} \left[ \|x^K(\Theta, d) - x_d^\star\|^2 \right].
\end{equation}
The primary alternative tranining loss for L2O is to use estimates of the Wasserstein-1 distance between the distribution of reconstructed signals and the distributio of true signals (e.g., see \citep{lunz2018adversarial}), which yields unsupervised training.




\subsubsection{Learned Parameter Roles}

A key aspect of unrolled L2O methods is to determine how to parameterize each update.
This subsection discusses some of the common considerations and roles for these parameters.

\paragraph{Learning parameters in the iterations}
First, the direct method of parameterization is to convert scalars/vectors/matrix/filters used in iterative algorithms into learnable parameters and learn them through data-driven training process. In this type of parameterization, learning can overcome the need to hand choose hyperparameters. For example, LISTA \citep{gregor2010learning} unrolls ISTA, which usually has update formulation
\begin{equation}
    x^{k+1} = \eta_{\lambda/L}\left(x^k - \frac{1}{L}A^T(Ax^k - d)\right), \label{eq:ISTA}
\end{equation}
where $\lambda$ is the coefficient before the $\ell_1$ regularization in LASSO, $L$ is the largest eigenvalue of $A^TA$, and $\eta_\theta(\cdot)$ is the coordinate-wise soft-thresholding function parameterized with threshold $\theta$\footnote{The soft-thresholding function takes the formula as $\eta_\theta(z)=\mathrm{sign}(z)\cdot\max(0,|z|-\theta)$}. Then LISTA parameterizes ISTA as a recurrent formulation
\begin{equation}
    x^{k+1} = \eta_{\theta}\left(W_ed + S x^k\right), \label{eq:LISTA}
\end{equation}
where taking $W_e\equiv A^T/L$, $S\equiv I-A^TA/L$ and $\theta\equiv\lambda/L$ reduces (\ref{eq:LISTA}) to ISTA.
A large amount of unrolling works follows the same methodology, but differ from each other in specific perspectives during the parameterization process, e.g. drop the recurrent structure in \citep{gregor2010learning} and use feed-forward modeling by untying the learnable parameters in different layers and update them independently \citep{HersheyRouxWeninger2014_deep,Borgerding_Schniter_2016}. This untying formulation can enlarge the capacity of the unrolled model \citep{HersheyRouxWeninger2014_deep} but can also cause the training instability due to the overwhelming parameter space and the difficulty of theoretical analysis on the convergence and other properties of the unrolled model.

To this end, there has been a line of efforts that strive to reduce the number of trainable parameters by relating and coupling different parts of parameterization \citep{xin2016maximal, Chen_Liu_Wang_Yin_2018} and analytically generate parameters that satisfy certain constraints \citep{LiuChenWangYin2019_alista}.
ALISTA \cite{LiuChenWangYin2019_alista} learns thresholding and step size parameters, reducing the number of parameters significantly -- to two scalars per layer, which stabilizes the data-driven process while achieving state-of-the-art recovery performance at the time with theoretical convergence guarantee.

Besides, recent trends also started to apply model-free L2O methods, to predict hyperparameters in classic iterative optimization algorithms. This is seen as the fusion among model-free L2O (since the algorithmic techniques are LSTM- or RL-based), model-specific L2O (as the update rule is eventually based on the classic iterates), and more traditional hyperparameter optimization. Examples include learning an adaptive learning rate schedule for training deep networks with stochastic graduate descent  \cite{xu2019learning,egidio2020learning,shu2020meta}; coordinating specifically for layer-wise training \cite{you2020l2} or domain adaptation speeds \cite{chen2020automated}; predicting the update combination weights and decay rates in Adam \cite{wang2019hyperadam}; or estimating training sample importance \cite{shu2020meta,chen2020self}, among others.

\paragraph{Deep priors learned from data}
Another popular parameterization in unrolling is to use a deep model to learn a data-dependent prior on the interested signals to replace hand-designed priors such as sparsity and total variation, which are not accurate in real-world applications and thus introduce bias.
However, previous works have various ways to use the ``learned prior''.
For example in \citep{rick2017one,metzler2017learned,gupta2018cnn,zhang2017learning}, people used data-driven training to learn a proximal or a projection operator that is iteratively used in the unrolled algorithm. The learned operator takes recovered signals contaminated by noises and artifacts as inputs and outputs a refined estimation. Sometimes the learned operator is found to fit a denoiser. The main difference of these learned prior methods from Plug-and-Play methods is that the operator is usually learned in a end-to-end way, making it overfitting to the current task or data distribution and not be able to be plugged into other algorithms. \cite{bora2017compressed,van2018compressed} supposed that the relevant vectors lie near the range of a generative model, and hence used generative adversarial networks (GANs) as a learning-based prior for image compressive sensing. 
\citep{lunz2018adversarial,mukherjee2020learned} perceive the prior as a loss function that can distinguish between coarsely recovered signals without considering any priors, and real on-domain signals with good quality. The prior as a loss function is adversarially trained to output high losses for coarse recoveries and low losses for real signals. Then we use the (sub-)gradient generative by the learned network via back-propagation in the unrolled algorithms e.g. gradient descent.

\paragraph{Others}
Besides the above two major ways of parameterization, there are also works that learn black-box agents that directly generate next-step estimations given the current iterate and historical information \citep{adler2017solving,Xiong_De_la_Torre_2013,Putzky_Welling_2017}. For instance, \citep{Behrens_Sauder_Jung_2020} learns an LSTM network that generates the step sizes and threshold parameters within each layer in ALISTA, instead of training them using back-propagation.


\subsection{Applications} 
\label{app}
Unrolled L2O schemes have found many applications. Below we identify several of these and note that there is some overlap among the three approaches discussed in this section.  


\begin{itemize}
    \item \textbf{Image Restoration and Reconstruction.} The model-based L2O methods, including both Plug-and-Play and unrolling methods, are widely used for various tasks in image restoration, enhancement, and reconstruction. Popular application examples include 
    denoising \citep{Diamond_Sitzmann_Heide_Wetzstein_2018,Putzky_Welling_2017,zhang2017learning,Chen_Pock_2017,sreter2018learned,LiuChenWangYin2019_alista,xie2019differentiable,lunz2018adversarial,mardani2019degrees},
    deblurring \citep{Diamond_Sitzmann_Heide_Wetzstein_2018,zhang2017learning,meinhardt2017learning,liu2018bridging,Corbineau_Bertocchi_Chouzenoux_Prato_Pesquet_2019,mardani2019degrees,mukherjee2020learned,zhang2019deep,li2020efficient},
    super-resolution \citep{rick2017one,Giryes_Eldar_Bronstein_Sapiro_2018,Putzky_Welling_2017,zhang2017learning,Chen_Pock_2017,liu2018bridging,zhang2019deep},
    inpainting \citep{rick2017one,Putzky_Welling_2017,sreter2018learned,Aberdam_Golts_Elad_2020}, and compressive sensing \citep{rick2017one,metzler2017learned,Chen_Liu_Wang_Yin_2018,Diamond_Sitzmann_Heide_Wetzstein_2018,Perdios_Besson_Rossinelli_Thiran_2017,mardani2018neural,zhang2018ista,Ito_Takabe_Wadayama_2019,mardani2019degrees}, JPEG artifacts reduction \citep{Wang_Liu_Chang_Ling_Yang_Huang_2016,Chen_Pock_2017,fu2019jpeg}, demosaicking \citep{meinhardt2017learning}, dehazing \citep{yang2018proximal} and deraining \cite{wang2020model}. Note that not all those works identically stick to the algorithm's original architecture; instead many only follow loosely the idea, and replace various components with convolutions or other modern deep learning building blocks. 
    \item \textbf{Medical and Biological Imaging.} We specifically separate Medical and Biology Imaging applications from the previous Image Restoration part because the former has its own special scenario and challenges. Imaging techniques in medical such as MRI and CT require accurate reconstructions of images with as few measurements as possible that result in minimal discomfort or side-effect. It is also challenging to extend the model-based methods in natural images to properly deal with the complex-valued inputs \citep{yang_Sun_Li_Xu_2016}. Other work that applies model-based L2O to medical and biology imaging includes \citep{adler2018learned,tamir2019unsupervised,adler2017solving,Diamond_Sitzmann_Heide_Wetzstein_2018,lunz2018adversarial,dardikman2020learned,mukherjee2020learned,Hara_Chen_Washio_Wazawa_Nagai_2019,solomon2019deep}.
    %
    \item \textbf{Wireless Communication.} Tasks in wireless communication systems can also be solved by unrolling methods, e.g. resource management \citep{sun2017learning,takabe2020complex,chowdhury2020unfolding}, channel estimation \citep{he2019model}, signal detection \citep{he2019model} and LDPC coding \citep{Wadayama_Takabe_2019}.For example, MIMO detection, which can be formulated as a sparse recovery problem, was shown to benefit from L2O methods based on ``deep unfolding" of an iterative algorithm added with trainable parameters \cite{balatsoukas2019deep}, such as LISTA. The model-based L2O approaches have already exhibited superior robustness and stability to low signal-to-noise (SNR), channel correlation, modulation symbol and MIMO configuration mismatches \cite{he2019model}. We refer the readers to a seminal survey about unfolding methods in communication systems \citep{balatsoukas2019deep}.
    \item \textbf{Seismic Imaging.} Another important application of model-based L2O is seismic imaging \citep{wang2018velocity,yang2019deep,yang2020deep,zhang2020can}.  Most of them adopt a more plug-and-play manner to learn CNN-based projectors that are first trained using data and then integrated into classic iterative updates, due to the desirable emphasis on the physical modeling.
    \item \textbf{Miscellaneous Applications:} such as \text{clustering} \citep{Wang_Chang_Zhou_Wang_Huang_2016, peng2016deep, peng2018k} and \text{classification} \citep{wang2016learning,zhou2018sc2net,cowen2019lsalsa}, phase retrieval \citep{icsil2019deep}, RNA second structure prediction \cite{chen2019rna}, speech recognition and source separation \cite{HersheyRouxWeninger2014_deep,wang2018end}, remote sensing \cite{lohit2019unrolled}, smart grid \cite{zhang2019real}, graph recovery \cite{Shrivastava2020GLAD:}, and photometric stereo estimation \cite{xin2016maximal}.
\end{itemize}

\subsection{Theoretical Efforts}  
\label{theory}

Although successful empirical results show significant potential of L2O, limited theory exists due to black-box training pertaining to such learned optimizers. That important gap often makes the broader usability of L2O questionable. We note that, the specific optimization problem and algorithm structure in model-based L2O often offer more opportunities for their theoretical analysis, by bridging us to the wealth of classic optimization tools, compared to model-free L2O. Therefore, model-based L2O in a few specific problems has been the main focus of existing theoretical efforts so far.

To explain why and how model-based L2O methods outperform traditional optimization algorithms, there are several essential questions at the heart of model-based L2O:
\begin{itemize}
    \item (Capacity). Given the model of L2O, do there exist parameters in the model that makes L2O provably outperform traditional optimization algorithms, over the task distribution (i.e., ``in distribution")? Is there any ``safeguard" mechanism available on L2O to ensure them at least as good as traditional algorithms even on examples out of the task distribution (i.e., ``out of distribution" or OoD)?
    \item (Trainability). Given the existence of such parameters, what training method should we use to obtain those parameters? Do guarantees exist that the training method converges to the ideal parameters?
    \item (Generalization). Does the trained model generalize, say, to testing instances from the same source of training instances (i.e., ``interpolation”)? Can the trained models ``extrapolate", e.g., on testing instances more complicated than any training one?
    \item (Interpretability). How can we explain what the L2O models have learned?
\end{itemize}
We list brief answers to the four questions here. Capacity and interpretability have been partially solved on signal/image-processing related problems by some recent works, as to be detailed in the subsections below. Generalization gains increasing attention recently and some works provide bounds of generalization gap of some specific L2O models. Lastly, to the best of our knowledge, there has been very limited theoretical work on the trainability of L2O, due to the high nonconvexity of the training objectives (see Section \ref{sec:obj}). The main exception here is the recent work \cite{wang2020guarantees} discussing the local minima and gradient explosion/vanishing in L2O training for quadratic minimization. This was also partially addressed by \cite{heaton2020projecting} where the authors used sorting activation functions \cite{anil2019sorting} and Lipschitz networks to encourage gradient norm preservation.

\subsubsection{Capacity}
\label{subsubsec:capacity}
To our best knowledge, \cite{sun2017learning} is the first effort on theories of L2O. It approximates a traditional method WMMSE (weighted minimum mean square error) with a fully-connected neural network and proves that the output of the neural network can be arbitrarily close to the result of WMMSE as long as its number of layers and units is large enough. In other words, this work estimates the approximation capacity of the neural network.

Approximation capacity is a generic notion in machine learning. Specifically for L2O, \textbf{convergence} can be used to describe the model capacity:
\emph{Do there exist parameters $\{\theta^k\}_k$ in (\ref{eq:model-unroll}) that make $\{x^k\}$ converge better than classic optimization algorithms as $k\to \infty$?}
The work \cite{Chen_Liu_Wang_Yin_2018} adopts this way to describe the convergence of L2O on the sparse recovery problem and give a convergence rate of LISTA which is better than that of classic algorithms ISTA/FISTA. \cite{LiuChenWangYin2019_alista,wu2019sparse,Aberdam_Golts_Elad_2020,ablin2019learning,yang2020learning,zarka2019deep} improve the theoretical result of LISTA by designing delicate models, in another word, designing operator $T$ in (\ref{eq:model-unroll}).  \cite{xie2019differentiable} analyzes the convergence rate of differentiable linearized ADMM. 

Instead of studying parameters $\theta^k$, another approach to establish convergence is to propose mathematical conditions on the operators (for example, operator $T$ in (\ref{eq:model-unroll}) and operator $H$ in (\ref{eq:pnp-admm-1})) in the models. These mathematical conditions should not only guarantee the convergence but also can be satisfied practically. \cite{ryu2019plug} proposes a continuity condition that can be satisfied by specifically designed training method. \cite{chen2020learnable} assumes smoothness of the operator that is satisfied by choosing smooth activation functions in a neural network (e.g., the sigmoid function). \cite{mukherjee2020learned} assumes convexity of the regularizer in their math proof and proposes to parameterize the regularizer by a convex neural network. \cite{heaton2020projecting} proves the convergence under assumptions about the manifold and sufficient representation by the sampled data, which are usually satisfied in practice.

While the above efforts focus on studying the convergence and acceleration effect of L2O over the target task distribution, a parallel and same important topic is to bound or characterize the convergence of L2O under OoD inputs: how much can the L2O convergence degrade when applied to optimizees deviating from the task distribution? Seemingly daunting at the first glance, that goal may be fulfilled by L2O with a safeguard mechanism, that can provide a way to establish convergence independent of the parameters and data. In this sense, the capacity of the original L2O models can be considered as enlarged as the convergence is attained on more OoD optimizees.
The most common approach is to i) compute a tentative L2O update using any method under the sun, ii) check if the tentative update yields a reduction in the value of an energy, and iii) accept the L2O update if the energy is less than some relative bound (e.g., the energy at the current iterate).  If the tentative L2O update is rejected, then a fallback scheme is to apply the update of a classic optimization algorithm. 
Because the energy is monotonically decreasing, it can be shown that the overall algorithm converges to a minimizer of the energy.
The energy is typically defined to be the objective function in the variational problem (VP) or its differential  \cite{liu2018bridging,liu2018learning,liu2019convergence,moeller2019controlling}.  
An alternative approach is define the energy to measure the residual between updates (e.g., \cite{heaton2020safeguarded}).
That is, if $T$ is the update operator for a classic algorithm, then at a tentative update $u^k$ we check if the residual $\|u^k-T(u^k)\|$ is less than some relevant bound.

\subsubsection{Interpretability}
Interpretability is significant to a learning model, now even more than ever. Some efforts have been made on the interpretability of L2O, mainly about linking or reducing their behaviors to those of some analytic, better understood optimization algorithms.
\cite{xin2016maximal} studies unrolled iterative hard-thresholding (IHT) and points out that unrolled IHT adopts better dictionary and a wider range of RIP condition than IHT.
\cite{Moreau_Bruna_2017} demonstrates that the mechanism of LISTA is related to a specific matrix factorization of the Gram matrix of the dictionary.
\cite{Giryes_Eldar_Bronstein_Sapiro_2018} explains the success of LISTA with a tradeoff between convergence speed and reconstruction accuracy.
\cite{LiuChenWangYin2019_alista} shows that the weights in LISTA has low coherence with the dictionary and proposes an analytic approach to calculate the weights.
\cite{Shrivastava2020GLAD:} analyzes alternating minimization algorithm quantitatively on a graph recovery problem and reveals the analogy between the learned L2O model and  solving a sequence of adaptive convex programs iteratively.
\cite{Takabe_Wadayama_2020} points out the analogy between deep-unfolded gradient descent and gradient descent with Chebyshev step-size and shows that the learned step size of deep-unfolded gradient descent can be qualitatively reproduced by Chebyshev step-size.

\subsubsection{Generalization}

Generalization is an important topic for L2O, just like for any other ML domain. Recent work \cite{metz2020tasks,almeida2021generalizable} described a ``generalization-first" perspective for L2O. The relationship between generalization gap and the number of training instances provides us an estimate on how many samples we should use. In the recent literature,  \cite{Chen_Zhang_Reisinger_Song_2020},  \cite{behboodi2020generalization} and \cite{van2020interpretable} studied  the Rademacher complexity, an upper bound of generalization gap.  \cite{Chen_Zhang_Reisinger_Song_2020} estimates the Rademacher complexity of gradient
descent and Nesterov’s accelerated gradients on a parameterized quadratic optimization problem; \cite{behboodi2020generalization} estimates the Rademacher complexity of a deep thresholding network on sparse linear representation problem.
\cite{van2020interpretable} proposes a reweighted RNN for the signal reconstruction problem and provides its Rademacher complexity.
As its definition suggests, generalization measures how the trained L2O models perform on unseen samples from the same distribution seen in training, but not on samples deviating from that distribution. Please refer to the previous discussion in Section.~\ref{subsubsec:capacity} on how methods such as Safeguard can mitigate the challenges of OoD optimizees.

Stability is highly related with generalization \cite{Chen_Zhang_Reisinger_Song_2020}, some L2O works analyze the stability of their methods.
For example, \cite{kobler2020total} estimates Lipschitz constants of their method to provide stability analysis results both with respect to input signals and network parameters. This yields a growth condition inequality.
Another measure of stability is to ensure stability of the model as the weighting of the regularizer term in (VP) tends to zero \cite{kofler2020neural,li2020nett}.
A third approach is that of $\Phi$-regularization used with null space networks \cite{schwab2019deep}, which learns a regularizer based on the null space of $A$.
Each of these methods yields different insights about the stability of L2O schemes for inverse problems.

Similar to the bias-variance trade-off in ML, L2O models are also subject to the trade-off between the capacity and the generalization. 
\cite{Chen_Zhang_Reisinger_Song_2020} provided insights on such trade-off in an unrolled model by proving a generalization bound as a function of model depth, and related to the properties of the algorithm unrolled.


\section{The Open-L2O Benchmark}
\label{sec:benchmark}

The community has made diverse attempts to explore L2O and generated a rich literature for solving different optimization tasks on different data or problems, as well as different software implementations on various platforms. Each L2O method has its own training recipe and hyperparameters. Because a common benchmark for the field has not yet been established, comparisons among different L2O methods are inconsistent and sometimes unfair. In this section, we present our efforts toward creating a comprehensive benchmark that enables fair comparisons. To the best of our knowledge, this is the first time of such an attempt. 

\textbf{Testbed problems}: We choose some popular and representative test problems that have been used in the existing L2O literature: (i) convex sparse optimization, including both sparse inverse problems and LASSO minimization; (ii) minimizing the nonconvex Rastrigin function, and (iii) training neural networks (NNs), which is a more challenging nonconvex minimization problem. 

\textbf{Task distributions}: Inspired by both research practice and real-world demands, we define task distributions in the following problem-specific way: (i) for sparse optimization, the optimizees during training and testing are optimization problems with the same objective function and decision variables but  \textit{different data}; (ii) in the Rastrigin-function test, during both training and testing, the optimizees have different decision variables and use random initializations; (iii) in the NN training test, the training optimizees use the same network architecture and dataset but random initializations; however, testing samples optimizees from a \textit{different distribution}, that is, the L2O optimizer is applied to train a network of an unseen architecture and on a different dataset. 

\textbf{Compared Methods and Evaluation Settings}: For each problem, we choose applicable approaches that include both model-free and/or model-based ones, implement them in the TensorFlow framework, ensure identical training/testing data, and evaluate them in the same but problem-specific metrics. 
After presenting the results, we draw observations from   these benchmark experiments. 

Our datasets and software are available as the \textbf{Open-L2O} package at: \url{https://github.com/VITA-Group/Open-L2O}. We hope that Open-L2O will foster reproducible research, fair benchmarking, and coordinated development efforts in the L2O field. 

\subsection{Test 1: Convex Sparse Optimization}

\subsubsection{Learning to perform sparse optimization and sparse-signal recovery}
\label{sec:exp_sparse_recovery}

\paragraph{Problem definition} The sparse recovery problem has been widely studied in the model-based L2O literature \cite{gregor2010learning,Chen_Liu_Wang_Yin_2018,LiuChenWangYin2019_alista}. The task is to recover a sparse vector from its noisy linear measurements:
\begin{align}
    b_q = A x^\ast_q + \varepsilon_q,
    \label{eq:sparse_recovery_forward}
\end{align}
\noindent
where $x_q^\ast\in\mathbb{R}^{n}$ is a sparse vector, $\varepsilon_q\in\mathbb{R}^{m}$ is an additive noise, and $q$ indexes an optimizee instance.
While $x^\ast_q,b_q,\varepsilon_q$ change for each $q$, the measurement matrix $A\in\mathbb{R}^{m\times n}$ is fixed across all training and testing. In practice, there is a matrix $A$ associated with each sensor or sensing procedure.

\paragraph{Data generation}
We generate $51,200$ samples as the training set and $1,024$ pairs as the validation and testing sets, following the i.i.d. sampling procedure in~\cite{Chen_Liu_Wang_Yin_2018}. We sample sparse vectors $x^*_q$ with components drawn i.i.d. drawn from the distribution $\mathrm{Ber}(0.1)\cdot N(0,1)$, yielding an average sparsity of $\sim10\%$. We run numerical experiments in four settings:
\begin{itemize}
    \item \textbf{Noiseless}. We take $(m,n)=(256,512)$ and sample $A$ with $A_{ij}\sim N(0,1/m)$ and then normalize its columns to have the unit $\ell_2$ norm. The noise $\varepsilon$ is always zero.
    \item \textbf{Noisy}. The same as above except for Gaussian measurement noises $\epsilon_q$ with a signal-to-noise ratio (SNR) of 20dB.
    \item \textbf{Coherent}. A Gaussian random matrix is highly incoherent, making sparse recovery relatively easy. To increase the challenge, we compute a dictionary $D\in\mathbb{R}^{256\times 512}$ from 400 natural images in the BSD500 dataset \cite{martin2001database} using the block proximal gradient method \cite{xu2013block} and, then, use it as the measurement matrix $A$, which has a high coherence. Other settings remain unchanged from the first setting above. 
    \item \textbf{Larger scale}. We scale up the noiseless case to $(m,n)=(512,1024)$; the other settings stay the same.
\end{itemize}

\paragraph{Model and training settings}
All our model-based L2O approaches take measurements $b_q$ as input and return  estimates $\hat{x}_q \approx x^*_q$. They are trained to minimize the mean squared error
$\mathbb{E}_{q\sim Q} \|\hat{x}_q - x^*_q\|_2^2$.
We adopt a progressive training scheme following~ \cite{Borgerding_Schniter_Rangan_2017,Chen_Liu_Wang_Yin_2018,wu2019sparse}. We use a batch size of 128 for training and a learning rate of $5\times 10^{-4}$. Other hyperparameters follow the default suggestions in their original papers.
After training, the learned models are evaluated on the test set in NMSE (Normalized Mean Squared Error) in the decibel
(dB) unit:
\[
\mathrm{NMSE}_\mathrm{dB}(\hat{x}_q,x^*_q)=
10 \log_{10}
\left(\|\hat{x}_q-x^*_q\|^2/\|x^*_q\|^2\right).\]

We compare the following model-based L2O approaches, all of which are unrolled to 16 iterations:
$(1)$ a feed-forward version \textbf{LISTA} as proposed in \cite{gregor2010learning}, with untied weights across layers;
$(2)$ \textbf{LISTA-CP}, a variant of LISTA with weight coupling \cite{xin2016maximal,Chen_Liu_Wang_Yin_2018};
$(3)$ \textbf{LISTA-CPSS}, a variant of LISTA with weight coupling and support selection techniques \cite{Chen_Liu_Wang_Yin_2018,LiuChenWangYin2019_alista};
$(4)$ \textbf{ALISTA} \cite{LiuChenWangYin2019_alista}, the variant of LISTA-CPSS with minimal learnable parameters;
$(5)$ \textbf{LAMP} \cite{Borgerding_Schniter_Rangan_2017}, an L2O network unrolled from the AMP algorithm;
$(6)$ \textbf{LFISTA} \cite{MoreauBruna2017_understanding}, an L2O network unrolled from the FISTA algorithm;
%
$(7)$ \textbf{GLISTA} \cite{wu2019sparse}, an improved LISTA model by adding gating mechanisms.
%

\begin{figure}[t]
\centering
\includegraphics[width=0.95\textwidth]{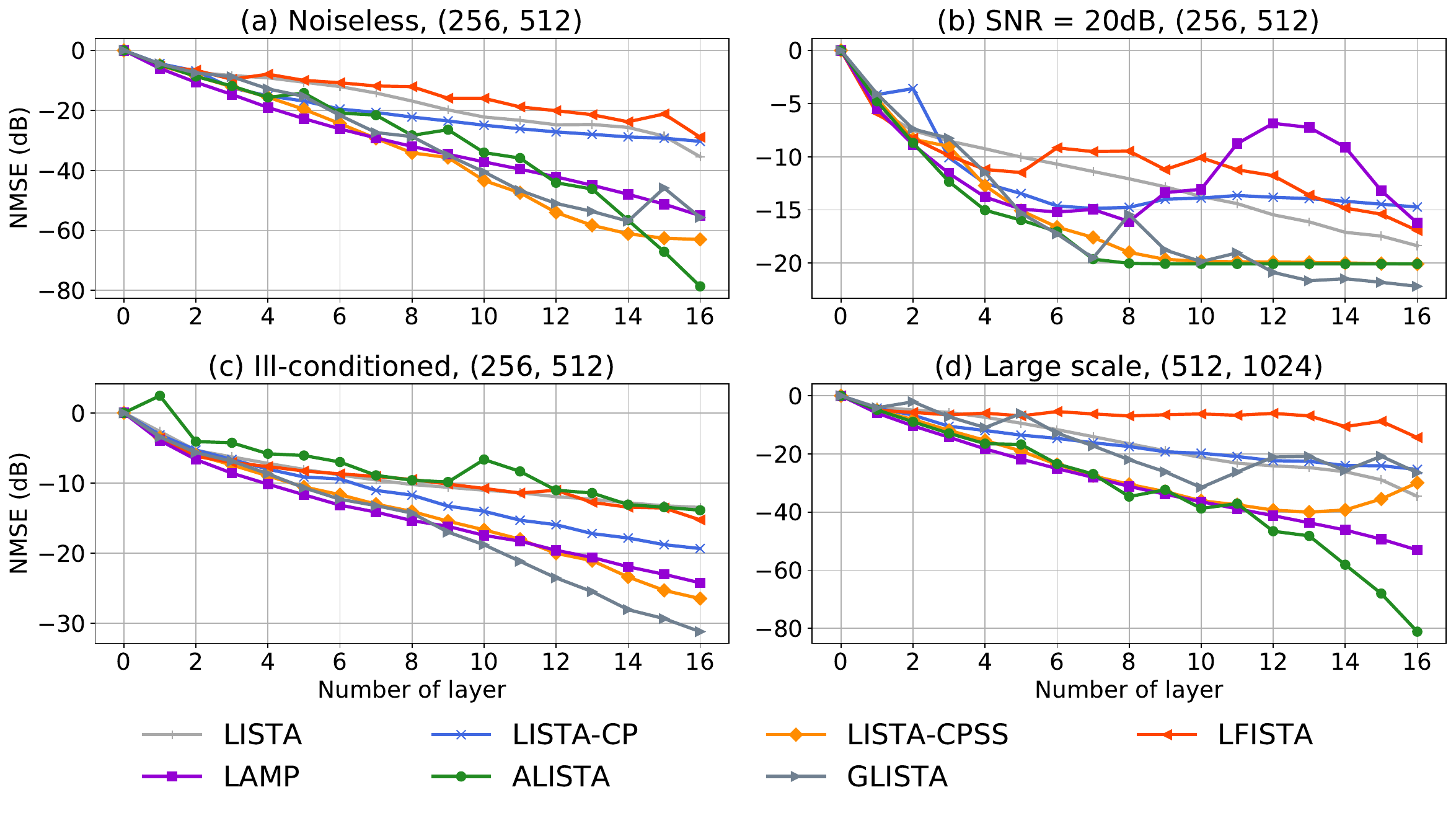}
\caption{Results of sparse recovery in four different settings: (a) noiseless with $(m,n)=(256,512)$; (b) additive Gaussian measurement noises of SNR=20dB and $(m,n)=(256,512)$; (c) coherent dictionary with $(m,n)=(256,512)$; (d) larger scale with $(m,n)=(512,1024)$. The x-axis counts the layers and the y-axis is the NMSE of the recovery.}
\label{fig:exp_sparse_recovery}
\end{figure}

\paragraph{Results} Figure~\ref{fig:exp_sparse_recovery} summarizes the comparisons; its four subfigures present the NMSEs of sparse recovery under noiseless, noisy, coherent dictionary, and large scale settings, respectively. We make the following observations
:
\begin{itemize}
    \item In the noiseless setting, ALISTA yields the best final NMSE, and the support selection technique in LISTA-CPSS contributes to good performance. LAMP and GLISTA achieve slightly worse NMSEs than LISTA-CPSS but are much better than LISTA. Surprisingly, LFISTA fails to outperform LISTA, which we attribute to its heavier parameterization in \cite{MoreauBruna2017_understanding} than other models and consequently harder training. 
    \item In the noisy setting of 20dB SNR, the performance of LAMP degrades severely while LISTA-CPSS and ALISTA still perform robustly. GLISTA outperforms all the others in the setting.
    \item In setting of a coherent dictionary, ALISTA suffers significantly from the coherence while the other methods that learn the weight matrices from data cope this issue better. In particular, GLISTA has the best performance thanks to the gate mechanisms.
    \item When it comes to a larger scale setting $(m,n)=(512,1024)$, LISTA, LISTA-CP, LFISTA and GLISTA have performance degradation due to the higher parameterization burden. In comparison, ALISTA becomes the best among all since it has the fewest parameters to learn: only the layer-wise thresholds and step-sizes. Therefore, ALISTA is least impacted by the problem scale and produces nearly the same performance as it does the smaller setting of $(m,n)=(256,512)$. 
\end{itemize}


\subsubsection{Learning to minimize the Lasso model} \label{sec:exp_lasso}

\paragraph{Problem definition} Instead of sparse recovery, which aims to recover the original sparse vectors, our goal here is to minimize the LASSO objective even though its solution is often different from the true sparse vector:
\begin{equation} \label{equ:lasso}
x^\mathrm{Lasso}_q = \argmin_x f_q(x),~~~\mbox{where}~f_q(x)= \frac{1}{2}\|A x-b_q\|_2^2+\lambda\|x\|_1,
\end{equation}
\noindent where $A\in \mathbb{R}^{m\times n}$ is a known, fixed, and normalized dictionary matrix, whose elements are sampled i.i.d. from a Gaussian distribution. An optimizee instance with index $q$ 
is characterized by a sparse vector $x^*_q$ and $b_q\in\mathbb{R}^{m\times 1}$ is the observed measurement under $A$ from $x^*_q$ following the same linear generation model in the previous subsection. $\lambda$ is a hyperparameter usually selected manually or by cross-validation, and is chosen to be $0.005$ by default in all our experiments.

Uniquely, we compare \textbf{both model-based and model-free} L2O methods in this section. For the former, we can adopt similar algorithm unrolling recipes as in the previous subsection. For the latter, we treat the problem as generic minimization and apply LSTM-based L2O methods. To our best knowledge, this is the first comparison between the two distinct L2O mainstreams. We hope the results provide quantitative understanding how much we can gain from incorporating problem-specific structures (when available) into L2O.

\paragraph{Data generation}
We run the experiments with $(m,n)=(5,10)$, as well as $(m,n)=(25,50)$. We did not go larger due to the high memory cost of LSTM-based model-free L2O methods. 
We sample $12,800$ pairs of $x^*_q$ and $b_q$ for training and $1,280$ pairs for validation and testing. The samples are noiseless. During the testing, we set $1,000$ iterations for both model-free L2O methods and classic optimizers. We run model-based L2O ones with a smaller fixed number of iterations (layers), which are standard for them. 

For each sample $b_q$, we let $f^*_q$ denote the optimal Lasso objective value, $f_q(x^\mathrm{Lasso})$. The optimization performance is measured with a modified relative loss:
\begin{equation}
    R_{f,\mathcal{Q}}(x) = \frac{\mathbb{E}_{q\sim\mathcal{Q}} [f_q(x) - f^*_q] }{\mathbb{E}_{q\sim\mathcal{Q}} [f^*_q] },
    \label{eq:lasso-metric}
\end{equation}
where the optimal solution is generated by $2,000$ iterations of FISTA and the expectations are taken over the entire testing set. 


We run four categories of methods to solve the LASSO optimization problem:
\begin{itemize}
    \item Three sub-gradient descent algorithms: GD, ADAM and RMSProp, which  serve as ``natural baselines" with neither learning nor problem-specific structure. All these algorithms use a step size of $10^{-3}$. We verified that changing step sizes did not notably alter their performance. (Even though GD and its generalizations do not handle nonsmooth objectives without smoothing or using proximal operators, they are popular optimizers that people try on almost everything, so we choose to test them anyway.)
    
    \item Two problem-specific analytic optimizers: ISTA \cite{blumensath2008iterative}, a forward backward splitting algorithm using a soft-thresholding backward operator and its Nesterov-accelerated version, FISTA \cite{beck2009fast}. For both methods, we use a step size of $1/L$, where $L$ is the largest eigenvalue of $A^T A$, and a threshold of $\lambda/L$ for the soft-thresholding function.
    
    \item Two model-based L2O models: vanilla LISTA \cite{gregor2010learning} and ALISTA \cite{LiuChenWangYin2019_alista}, both unrolled to 16 layers. We follow the same training setting in the last subsection except that we replace the training loss from the mean squared error (w.r.t. the true solution) with the LASSO loss. Since other compared methods take far more iterations, we also tried to expand LISTA/ALISTA to more iterations during testing by appending FISTA iterations. 
    
    \item Three model-free L2O methods: L2O-DM \cite{andrychowicz2016learning}, L2O-RNNprop \cite{lv2017learning} and L2O-enhanced optimizers \cite{chen2020training} (Section \ref{sec:model_free_basics}). All L2O optimizers are trained with the Lasso objective function as the reward for $100$ epochs and $1,000$ iterations per epoch. The reported evaluation is the average performance over $10$ random starting points.
\end{itemize}


\begin{figure}[t]
\subfigure[($m$,$n$)=(5,10)]{
\begin{minipage}[ht]{0.48\linewidth}
    \centering
    \includegraphics[width=1\textwidth]{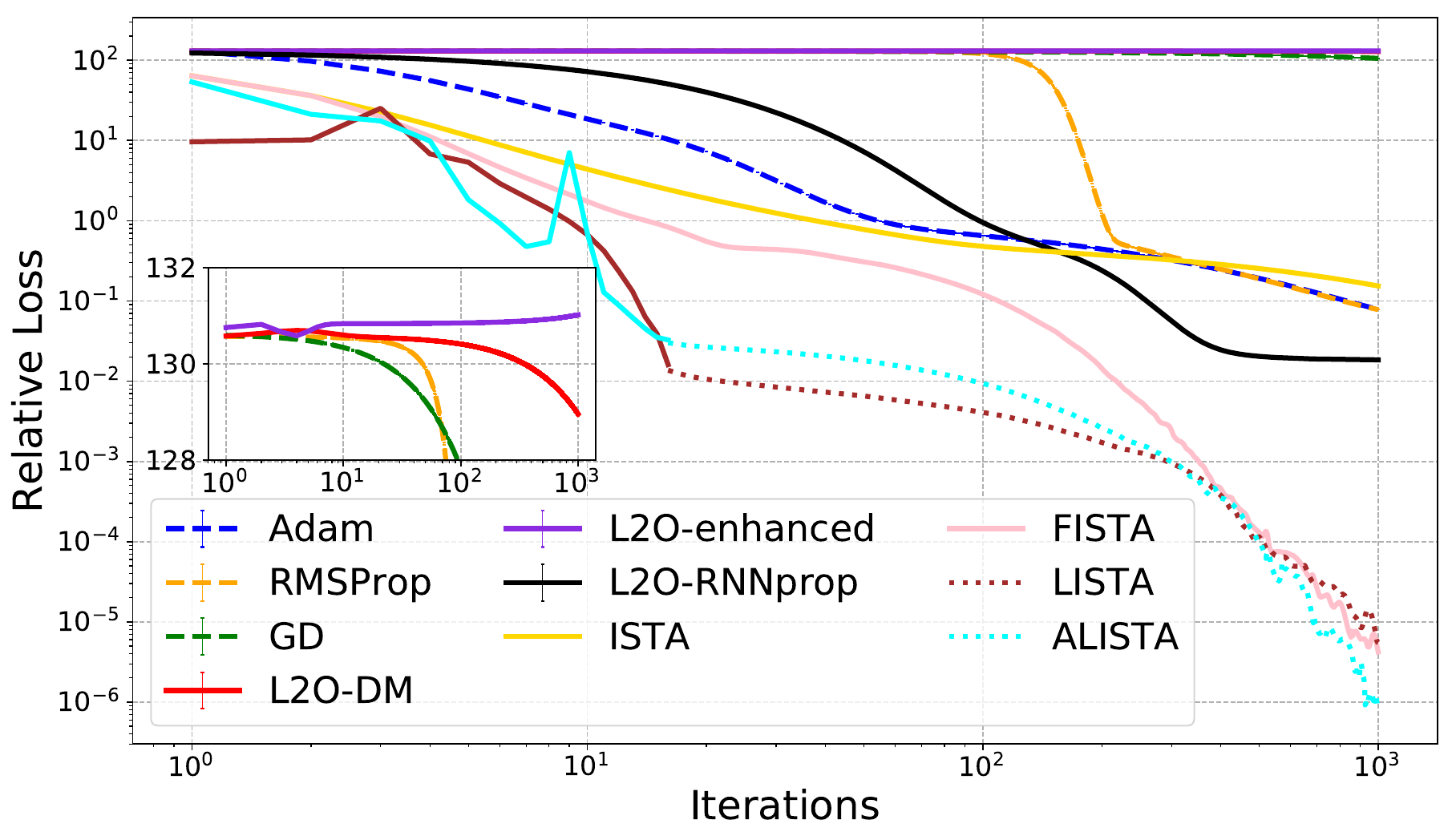}
\end{minipage}
}
\quad
\subfigure[($m$,$n$)=(25,50)]{
\begin{minipage}[!ht]{0.48\linewidth}
    \centering
    \includegraphics[width=1\textwidth]{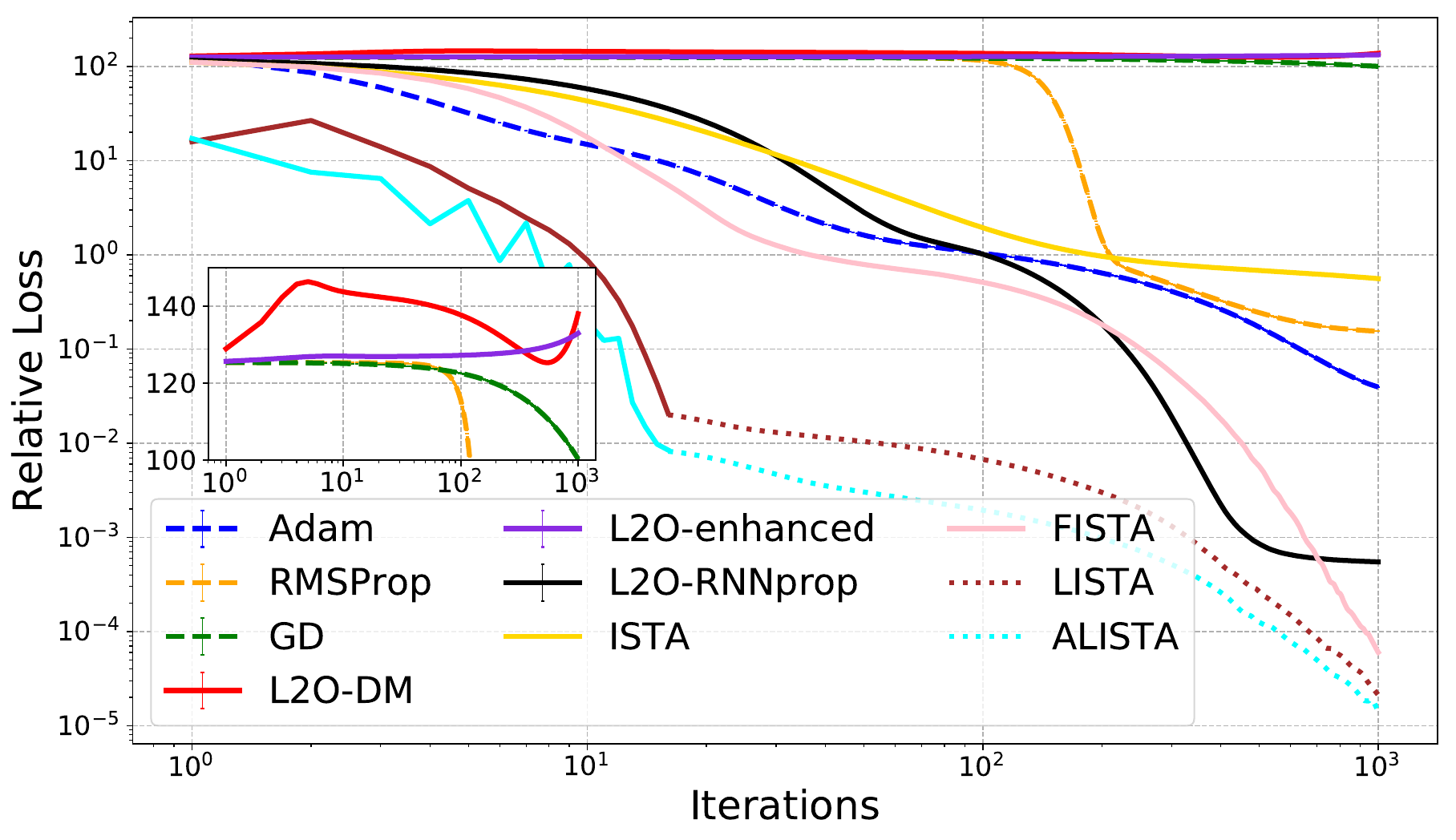}
\end{minipage}
}
\caption{Evaluation comparisons among analytic, model-based L2O, and model-free L2O optimizers on Lasso. y-axis represents the modified relative loss (\ref{eq:lasso-metric}), and x-axis denotes the number of iterations, both in the logarithmic scale.}
\label{fig:lassologloss}
\end{figure}

\paragraph{Results} From Figure \ref{fig:lassologloss},  we make the following observations:
\begin{itemize}
    \item At the small problem size $(5, 10)$, both ISTA and FISTA converge fast in tens of iterations with FISTA being slightly ahead. Both model-based L2O models, ISTA and ALISTA of 16 iterations (layers), converge to solutions of precision compared to what FISTA can achieve after hundreds of iterations and better than what ISTA can do at 1,000 iterations. The advantage of ISTA/ALISTA can sustain beyond 16 iterations using FISTA updates.
    \item In comparison, at the small problem size $(5, 10)$, model-free L2O methods exhibit far worse performance. L2O-RNNprop has slower convergence than ISTA/FISTA  and only produces reasonably good solutions after 1,000 iterations -- though still much better than analytic optimizers Adam and RMSProp. L2O-DM, L2O-enhanced, and vanilla GD completely fail to decrease the objective value.
    \item At the larger problem size $(25, 50)$, LISTA and ALISTA still converge to high-precision solutions with only 16 iterations with sustained advantages from the FISTA extension. Interestingly, ISTA now becomes much slower than FISTA. All the sub-gradient descent and model-free L2O methods remain to perform poorly; only L2O-RNNprop can converge faster than ISTA and comparable to FISTA, though reaching lower precision.
    \item Our experiments clearly demonstrate the dominant advantage of incorporating problem-specific structures to the optimizers when it comes to {both analytic and learned} optimizers. 
\end{itemize}



\subsection{Test 2: Minimization of non-convex function Rastrigin} \label{sec:exp_rastrigin}
We now turn to non-convex minimization. One popular non-convex test function is called the Rastrigin function:
\begin{equation} \label{equ:ras_test}
f(\bm{x}) =  \frac{1}{2}\sum_{i=1}^{n}x_i^2 - \sum_{i=1}^{n} \alpha\cos{(2\pi x_i)} + \alpha n,
\end{equation}
where $\alpha=10$. We consider a broad family of similar functions $f_q(\bm{x})$ that generalizes Rastrigin function:
\begin{equation} \label{equ:ras_train}
 f_q(\bm{x}) = \frac{1}{2}\|\bm{A}_q\bm{x}-\bm{b}_q\|_2^2 - \alpha\bm{c}_q\cos(2\pi\bm{x}) + \alpha n, 
\end{equation}
where $\bm{A}_q \in \mathbb{R}^{n \times n}$, $\bm{b}_q \in \mathbb{R}^{n \times 1}$ and  $\bm{c}_q \in \mathbb{R}^{n \times 1}$ are parameters whose elements are sampled i.i.d. from $\mathcal{N}(0,1)$. Obviously, the function (\ref{equ:ras_test}) is a special case in this family with $\bm{A}=\bm{I}, \bm{b}=\{0,0,\ldots,0\}^{\mathrm{T}}, \bm{c}=\{1,1,\ldots,1\}^{\mathrm{T}}$. For training, we sample $1,280$ triplets of \{$A_q$, $b_q$, $c_q$\} om two problem scales: $n=2$ and $n=10$. For evaluation, we sample another $128$ combinations of $A_q$, $b_q$ and $c_q$ and report the average performance over $10$ random starting points. The number of steps for evaluation is $1,000$.

\begin{figure}[ht]
\subfigure[$n$=2]{
\begin{minipage}[ht]{0.48\linewidth}
    \centering
    \includegraphics[width=1\textwidth]{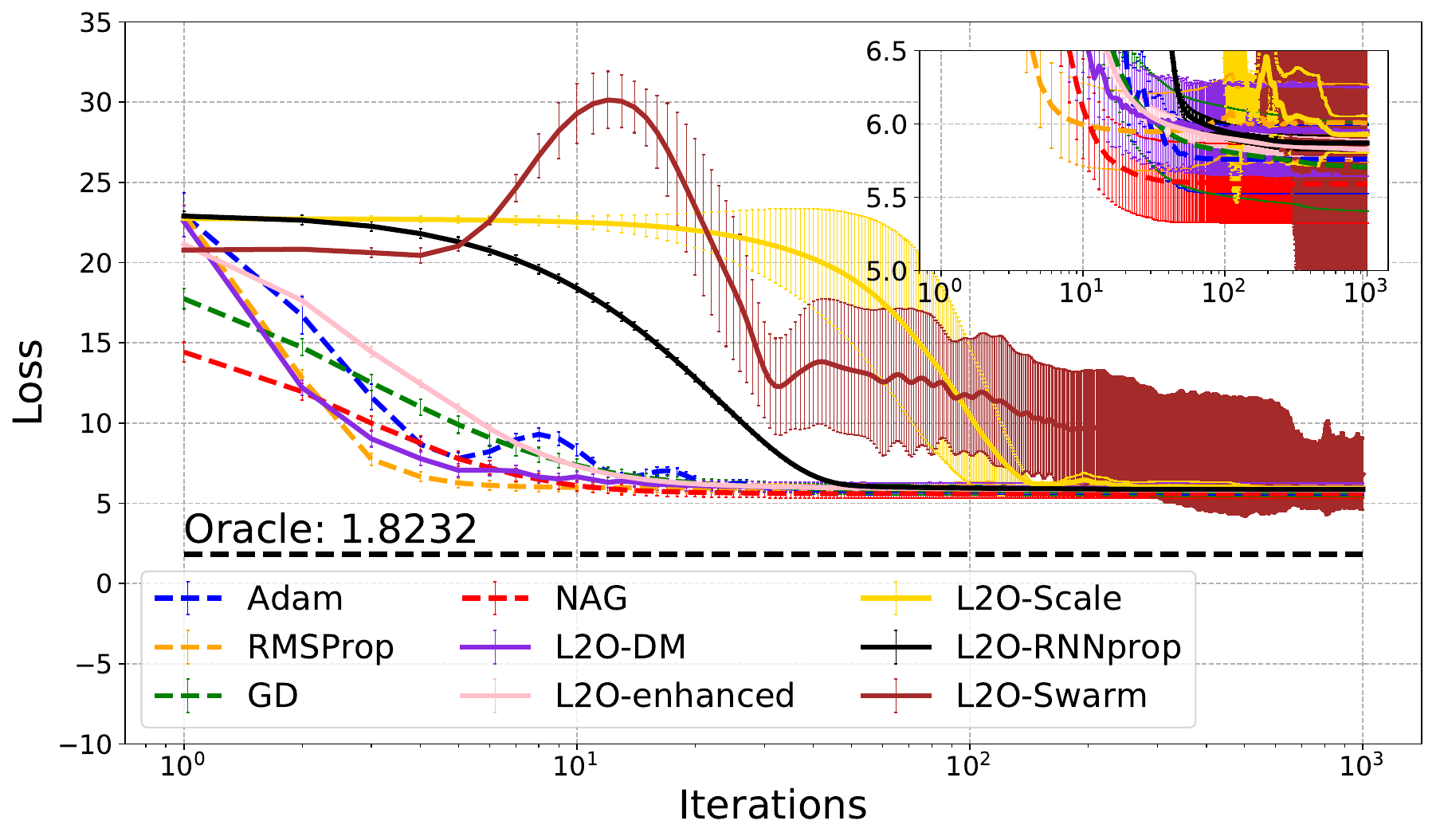}
\end{minipage}
}
\quad
\subfigure[$n$=10]{
\begin{minipage}[!ht]{0.48\linewidth}
    \centering
    \includegraphics[width=1\textwidth]{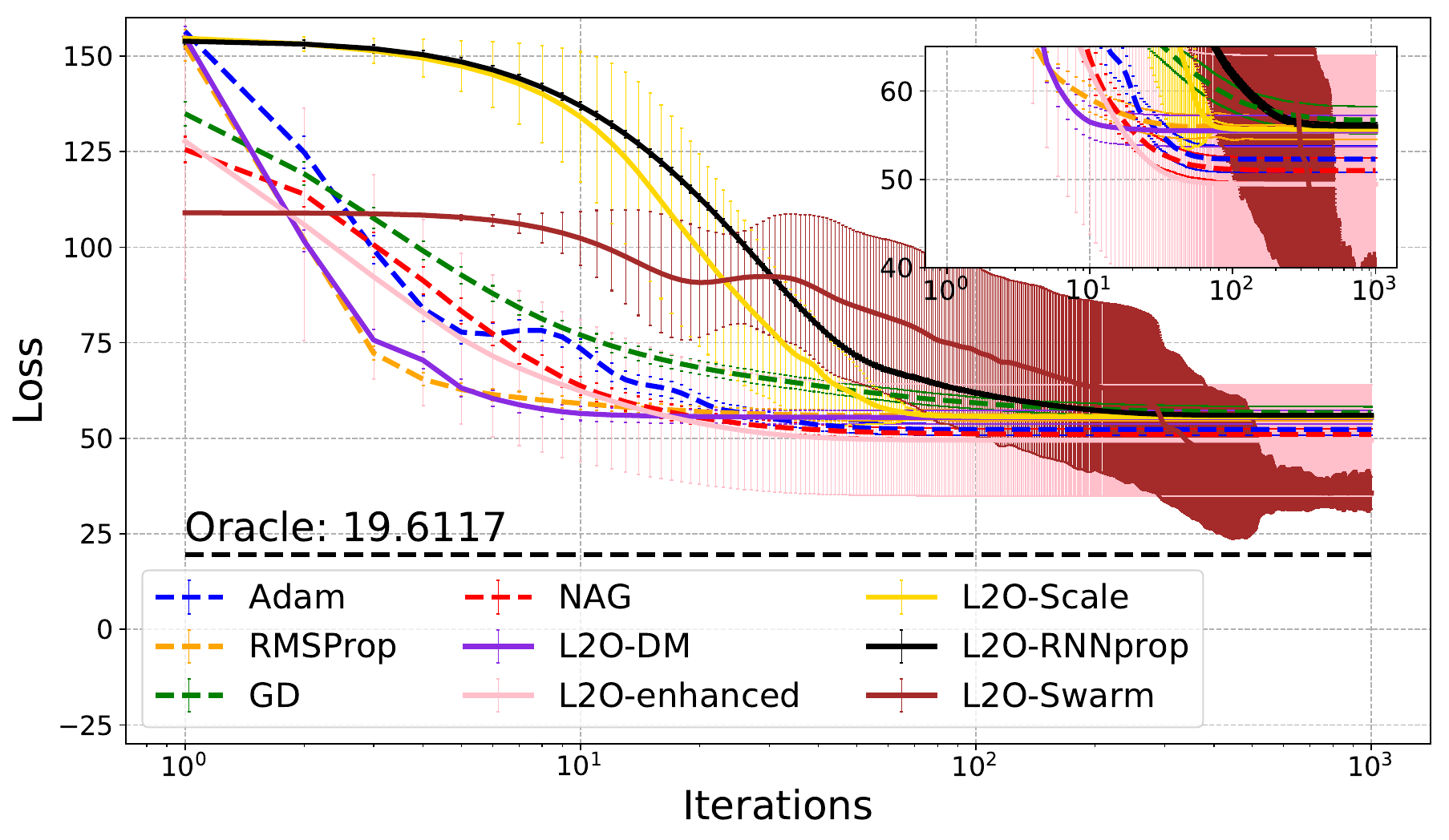}
\end{minipage}
}
\caption{Evaluation and comparison across model-free L2O and analytic optimizers on the generalized family of Rastrigin functions. y-axis represents the average loss of function values, and x-axis denotes the number of iterations in the logarithm scale.}
\label{fig:raslogloss}
\end{figure}

Two groups of methods are compared: (1) Four traditional gradient descent algorithms, including ADAM with a $10^{-1}$ step size, RMSProp with a $3\times 10^{-1}$ step size, GD with the line-searched step size started from $10^{-1}$, and NAG (Nesterov Accelerated Gradient) with the line-searched step size started from $10^{-1}$. All other hyperparamters are tuned by careful grid search. (2) Five model-free L2O methods, including L2O-DM \cite{andrychowicz2016learning}, L2O-enhanced \cite{chen2020training}, L2O-Scale \cite{wichrowska2017learned}, L2O-RNNprop \cite{lv2017learning} (Section \ref{sec:model_free_basics}), and L2O-Swarm \cite{cao2019learning} (Section \ref{sec:model_free_more}). All L2O optimizers are trained for $100$ epochs and $1000$ iterations per epoch. At the testing stage, we evaluate and report the logarithmic loss of unseen functions from the same family,  which are plotted in Figure~\ref{fig:raslogloss}.

\paragraph{Results} From Figure~\ref{fig:raslogloss} we draw the following observations:
\begin{itemize}
     \item ADAM and RMSProp converge slightly more quickly than GD and NAG with line search, especially in the early stage, for both $n=2$ and $n=10$. All analytic optimizers converge to local minima of similar qualities for $n=2$, and NAG finds a slightly better solution for $n=10$.
    \item L2O-DM, L2O-enhanced, L2O-RNNProp, and L2O-Scale perform similarly to analytic optimizers regarding both solution quality and convergence speed, showing no obvious advantage over analytic optimizers. For $n=10$, L2O-enhanced finds a solution of better quality than the other two. But the three model-free L2Os have larger error bars at $n=10$, which indicates model instability.
    \item Although not converging faster than others, L2O-Swarm locates higher quality solutions (lower losses) for both $n$ values, especially $n=10$. Since L2O-Swarm is the only method that leverages a population of LSTM-based L2O ``particles", it explores a larger landscape than the other methods, so it is not surprising that its higher search cost leads to better solutions.
    \item The oracle objectives (dashed black lines) in Figure~\ref{fig:raslogloss} are generated by the Run-and-Inspect method~\cite{chen2019run}, which can provably find a global minimum if the objective function can be decomposed into a smooth strongly convex function (e.g., a quadratic function) plus a restricted function (e.g., sinusoidal oscillations). For $n=2$, almost all methods reach a similar loss lied in [$5.5,6.0$] in the end. For $n=10$, L2O-Swarm performs significantly better than other methods in terms of the achieved final loss, though it converges more slowly than most of the other approaches.
\end{itemize}

\subsection{Test 3: Neural Network Training} \label{sec:exp_nn}

Our last test is training multi-layer neural networks (NNs), one of the most common tasks of L2O since its beginning. This has been the playground for model-free L2O methods. There are few problem-specific structures to explore. Common optimizers are stochastic gradient descent and their variants. We hope this test to address two questions:
\begin{itemize}
    \item Can model-free L2O optimizers outperform analytic ones on neural network training? If so, how much is the margin?
    \item Can model-free L2O optimizers generalize to unseen network architectures and data?
\end{itemize}
To fairly compare different methods and answer these questions, we train L2O optimizers on the same neural network used in~\cite{andrychowicz2016learning}: a simple Multi-layer Perceptron (MLP) with one $20$-dimension hidden layer and the sigmoid activation function, trained on the MNIST dataset to minimize the cross-entropy loss. Therefore, the task distribution becomes optimizing the same MLP model with different random initializations.

We probe the (out of distribution) generalizability of the learned optimizers by testing them on two unseen optimization tasks, following the practice in~ \cite{andrychowicz2016learning,chen2020training}:
\vspace{-0.5em}
\begin{enumerate}[leftmargin=*]
    \item Train another MLP with one $20$-dimension hidden layer, but using the ReLU activation function, on the MNIST dataset.
    \item Train another ConvNet on the MNIST dataset, consisting two convolution layers, two max pooling layers, and one last fully connected layer. The first convolution layer uses $16$ $3\times 3$ filters with stride $2$. The second convolution layer uses $32$ $5\times 5$ filters with stride $2$. The max pooling layers have size $2\times 2$ with stride $2$.
\end{enumerate}
We compare four model-free L2O optimizers: L2O-DM \cite{andrychowicz2016learning}, L2O-enhanced \cite{chen2020training}, L2O-Scale \cite{wichrowska2017learned}, and L2O-RNNprop \cite{lv2017learning} (Section \ref{sec:model_free_basics}), all following the hyperparameters suggested in their original papers. SGD, Adam, and RMSProp are analytical optimizers served as the baselines. 
All L2O optimizers are trained with the single model from 10,000 different random initializations drawn from $\mathcal{N}(0,01)$. On each optimizee, now corresponding to a random initialization, the optimizers run for $100$ iterations. The training uses a batch size of $128$. During each run of testing, we evaluate learned optimizers on an unseen testing optimizee for $10,000$ steps, which is much more than the training iteration number. We then report the training loss of the optimizee. We perform 10 independent runs and report the error bars.

It is worth mentioning that current L2O-optimizers can hardly scale to training large neural networks with over $1\times 10^6$ parameters. This is because to update the L2O optimizer using back-propagation, we need to keep the gradients during unrolling and the computation graph of the optimizee in memory, therefore, consuming too much memory. This remains an open challenge.

\begin{figure}[ht]
\subfigure[MLP]{
\begin{minipage}[ht]{0.48\linewidth}
    \centering
    \includegraphics[width=1\textwidth]{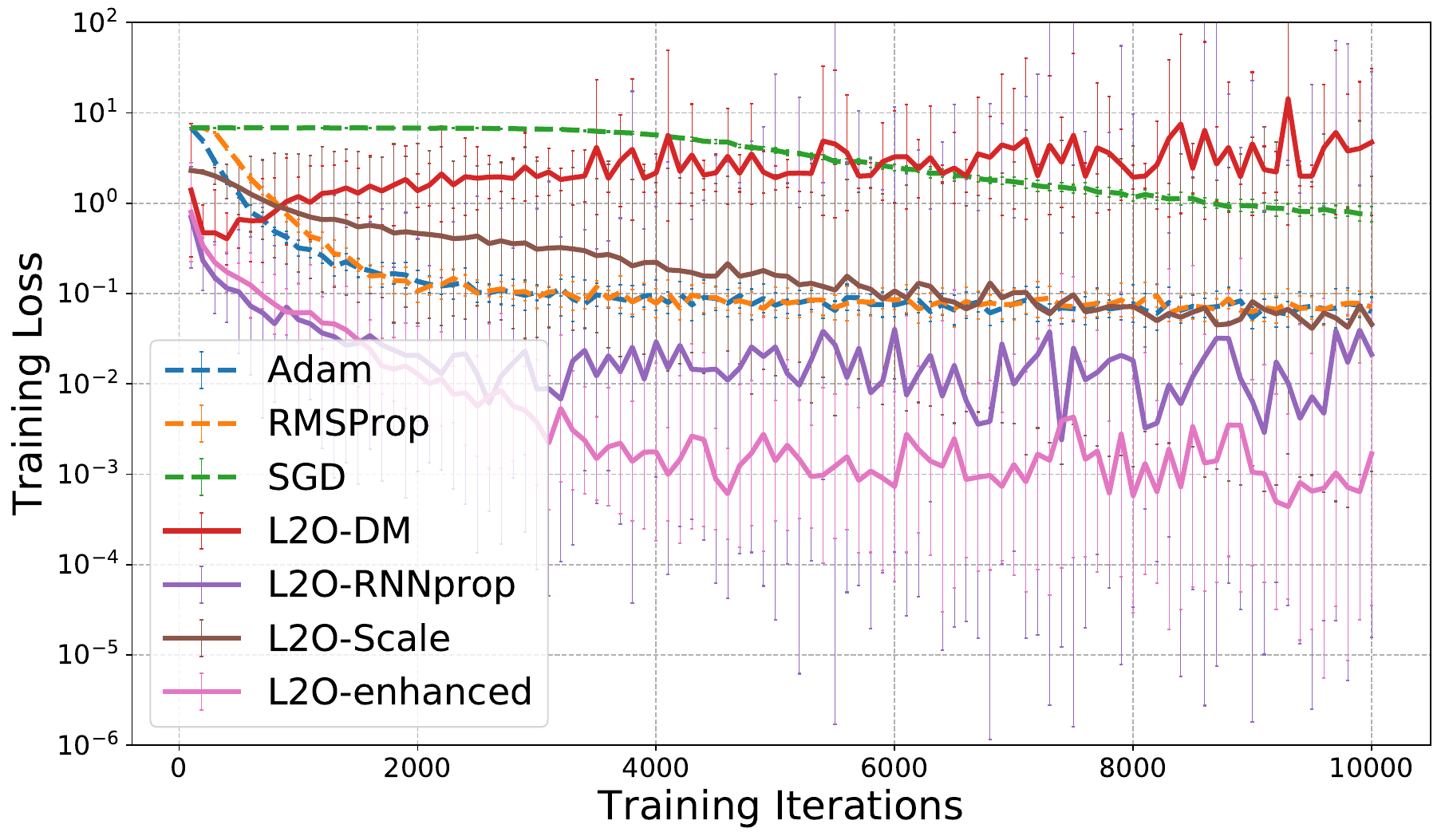}
\end{minipage}
}
\quad
\subfigure[ConvNet]{
\begin{minipage}[!ht]{0.48\linewidth}
    \centering
    \includegraphics[width=1\textwidth]{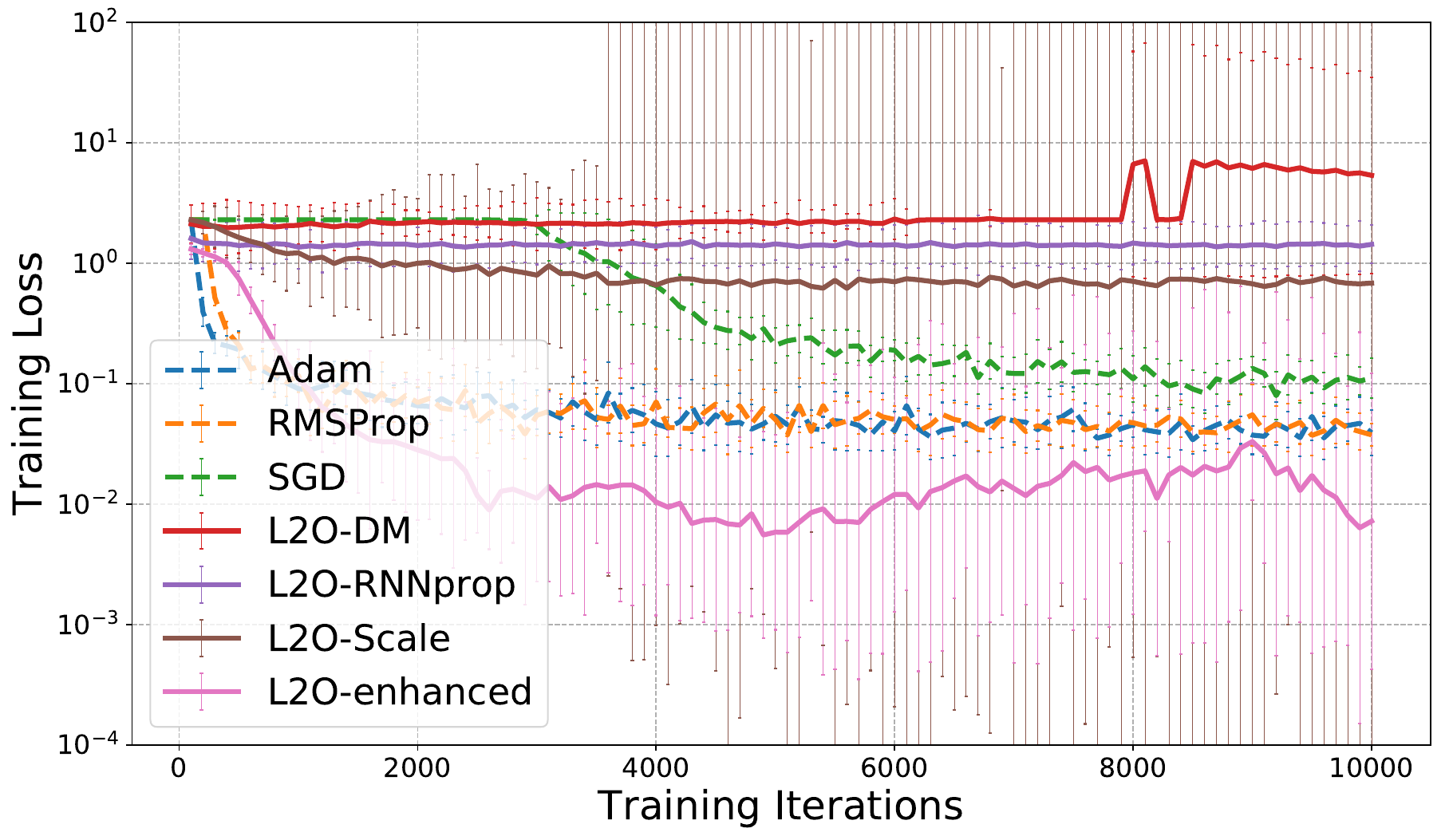}
\end{minipage}
}
\caption{Evaluation and comparison across model-free L2O and analytic optimizers on neural network training. y-axis represents the average training loss values, and x-axis denotes the number of iterations. Each curve is the average of ten runs, and error bars are plotted.}
\label{fig:nn}
\end{figure}

\paragraph{Results}
From Figure~\ref{fig:nn}, we make the following observations:
\begin{itemize}
    \item The relative performance of analytic optimizers is consistent: both Adam and RMSProp achieve similar training losses on both ReLU-MLP and ConvNets, and better than SGD. However, the learned optimizers exhibit different levels of transferability. L2O-DM completely diverges on ConvNet, but works on ReLU-MLP when the number of iterations is small, showing poor generalization. 
    \item L2O-Scale and L2O-RNNprop can generalize to ReLU-MLP as they achieve similar or lower training losses than analytic optimizer), but they cannot effectively optimize ConvNet. L2O-RNNprop only works for a small number of iterations less than 4,000. L2O-Scale generalizes better to larger iteration numbers, thanks to its random scaling during training.
    \item L2O-enhanced, which adopts increasing unroll length and off-line imitation strategy, achieves the best performance on both unseen networks after extended (optimizee) training iterations.
    \item The issue of lacking stability during testing exists for all L2O optimizers. On both ReLU-MLP and ConvNet, all L2O optimizers suffer from larger loss variances as the training iteration grows, causing the errors to fluctuate more severely that those of the analytic optimizers. 
\end{itemize}

\subsection{Take-home messages}
\begin{itemize}
    \item[$\star$] When appropriate model or problem structures are available to be leveraged, model-based L2O methods clearly outperform model-free ones and analytic optimizers. On optimizees sampled from the same task distribution, model-based L2O optimizers have solid, consistent, and reliable performance.
    \item[$\star$] Although in many cases we can also observe performance gains of model-free L2O approaches over analytic ones, the gains are not consistent.
    The benefits of current model-free L2O approaches come to questions when they are applied to larger optimizees or applied with more iterations during testing, let alone to ``out-of-distribution" optimizees.
    We find no clear winner in existing model-free L2O approaches that can consistently outperforms others across most tests. 
\end{itemize}

\section{Concluding Remarks}

This article provided the first ``panoramic" review of the state-of-the-arts in the emerging L2O field, accompanied with the first-of-its-kind benchmark. The article reveals that, despite its promise, the L2O research is still in its infancy, facing many open challenges as well as research opportunities from practice to theory.

On the theory side, Section \ref{theory} has listed many open theoretical  questions for us to understand why and how model-based L2O methods outperform traditional optimizers. Besides those, for model-free L2O methods, the theoretical foundation has been scarce if any at all. For example, although the training of L2O optimizers is often successful empirically, there has been almost no theoretical result established for the convergence performance of such a L2O training process, putting the general feasibility of obtaining L2O models in question. Also, for both model-based and mode-free L2O approaches, the generalization or adaption guarantees of trained L2O to optimizees out of the task distribution are under-studied yet highly demanded. 

On the practical side, the scalability of model-free L2O methods is perhaps the biggest hurdle for them to become more practical. That includes scaling up to both larger and more complicated models, and to more iterations during meta-testing. For model-specific L2O, the current empirical success is still limited to a few special instances in inverse problems and sparse optimization, and relying on case-by-case modeling. An exploration towards broader applications, and perhaps building a more general framework, is demanded. 



Furthermore, there is no absolute border between mode-based and model-free L2O methods, and the spectrum between the two extremes can hint many new research opportunities. \cite{monga2019algorithm} suggested a good perspective that unrolled networks (as one example of model-based L2O) is an intermediate state between generic networks and analytic algorithms, and might be more data-efficient to learn. The view was supported by \cite{meng2021a}, which further advocated that the unrolled model might be considered as a robust starting point for subsequently data-driven model search. We also believe that end-to-end learning approaches can be improved with current continuous optimization algorithms to benefit from the theoretical guarantees and state-of-the-art algorithms already available.

So, to conclude this article, let us quote Sir Winston Churchill: \textit{``Now this is not the end. It is not even the beginning of the end. But it is, perhaps, the end of the beginning."} Although most approaches we discussed in this paper are still at an exploratory level of deployment, and are apparently not yet ready as general-purpose or commercial solvers, we are strongly confident that machine learning has just began to feed the classic optimization field, and the blowout of L2O research progress has yet to start.

\newpage

\appendix

\section*{Appendix A. List of Abbreviations}

\begin{table}[h!]
\centering
{\footnotesize
\begin{tabular}{ccc}
\toprule
Abbreviation & Full Name & Description \\ \midrule
Adam & Adaptive Moment Estimation & {\scriptsize \begin{tabular}{@{}l@{}}A first-order gradient-based optimization algorithm\\of stochastic objective functions, based on adaptive\\estimates of lower-order moments.\end{tabular} } \\ \hline
ADMM & {\footnotesize \begin{tabular}{@{}c@{}}Alternating Direction\\Method of Multipliers\end{tabular} } & {\scriptsize \begin{tabular}{@{}l@{}}An optimization algorithm that solves convex\\problems by breaking them into smaller pieces, each\\of which are then easier to handle.\end{tabular} } \\ \hline
AutoML & Automated Machine Learning & {\scriptsize \begin{tabular}{@{}l@{}}The process of automating the process of applying\\machine learning to real-world problems.\end{tabular} }\\ \hline
CASH & {\scriptsize \begin{tabular}{@{}c@{}}Combined Algorithm Selection\\and Hyperparameter Optimization\end{tabular} } & {\scriptsize \begin{tabular}{@{}l@{}}A model selection strategy that consider the\\choosing of data preparation, learning algorithm,\\and algorithm hyperparameters as one large global\\optimization problem.\end{tabular} } \\ \hline
GAN & Generative Adversarial Network & {\scriptsize \begin{tabular}{@{}l@{}}A class of machine learning frameworks \cite{goodfellow2014generative},\\where two neural networks contest with each other\\in a zero-sum game, to learn to generate new data\\with the same statistics as a given training set.\end{tabular}} \\ \hline
HPO & Hyperparameter Optimization & {\scriptsize \begin{tabular}{@{}l@{}}The problem of choosing a set of optimal\\hyperparameters for a machine learning algorithm.\end{tabular} }\\ \hline
L2O & Learning to Optimize & {\scriptsize \begin{tabular}{@{}l@{}}Learnable optimizers to predict update rules (for\\optimizees) fit from data.\end{tabular} }\\ \hline
LASSO & {\scriptsize \begin{tabular}{@{}c@{}}Least Absolute Shrinkage\\and Selection Operator\end{tabular} } & {\scriptsize \begin{tabular}{@{}l@{}}A linear regression method that uses shrinkage to\\encourages simple, sparse models.\end{tabular} } \\ \hline
LSTM & Long-short Term Memory & {\scriptsize \begin{tabular}{@{}l@{}}A variant of artificial recurrent neural network\\(RNN) architecture, typically used to process\\sequence data.\end{tabular} }\\ \hline
MLP & Multi-layer Perceptron & {\scriptsize \begin{tabular}{@{}l@{}}A class of feedforward artificial neural network\\(ANN) composed of multiple layers of perceptrons \\ (with threshold activation).\end{tabular} }\\ \hline
PGD & Projected Gradient Descent & {\scriptsize \begin{tabular}{@{}l@{}}Projected gradient descent minimizes a function\\subject to a constraint. At each step we move in\\the direction of the negative gradient, and then\\``project'' onto the feasible set.\end{tabular} } \\ \hline
ReLU & Rectified Linear Unit & {\scriptsize \begin{tabular}{@{}l@{}}An activation function used in neural networks,\\defined as $f(x) = \mathrm{max}(x, 0)$.\end{tabular} } \\ \hline
RL & Reinforcement Learning & {\scriptsize \begin{tabular}{@{}l@{}}RL is an area of machine learning regarding how to\\make intelligent agents take actions in an environment\\to maximize the cumulative reward.\end{tabular} } \\ \hline
RMSProp & Root Mean Square Propagation & {\scriptsize \begin{tabular}{@{}l@{}}RMSprop gradient descent maintains a moving\\(discounted) average of the square of gradients, and\\then divide the gradient by the root of this average.\end{tabular} } \\
\bottomrule
\end{tabular}
}
\end{table}



\vskip 0.2in

\bibliography{L2O}

\begin{thebibliography}{247}
\providecommand{\natexlab}[1]{#1}
\providecommand{\url}[1]{\texttt{#1}}
\expandafter\ifx\csname urlstyle\endcsname\relax
  \providecommand{\doi}[1]{doi: #1}\else
  \providecommand{\doi}{doi: \begingroup \urlstyle{rm}\Url}\fi

\bibitem[Chen et~al.(2018)Chen, Liu, Wang, and Yin]{Chen_Liu_Wang_Yin_2018}
Xiaohan Chen, Jialin Liu, Zhangyang Wang, and Wotao Yin.
\newblock Theoretical linear convergence of unfolded ista and its practical
  weights and thresholds.
\newblock In \emph{Advances in Neural Information Processing Systems}, pages
  9061--9071, 2018.

\bibitem[Liu et~al.(2019{\natexlab{a}})Liu, Chen, Wang, and
  Yin]{LiuChenWangYin2019_alista}
Jialin Liu, Xiaohan Chen, Zhangyang Wang, and Wotao Yin.
\newblock {ALISTA}: Analytic weights are as good as learned weights in {LISTA}.
\newblock In \emph{International Conference on Learning Representations},
  2019{\natexlab{a}}.

\bibitem[Nesterov(2005)]{nesterov2005smooth}
Yurii Nesterov.
\newblock Smooth minimization of non-smooth functions.
\newblock \emph{Mathematical programming}, 103\penalty0 (1):\penalty0 127--152,
  2005.

\bibitem[Gregor and LeCun(2010)]{gregor2010learning}
Karol Gregor and Yann LeCun.
\newblock Learning fast approximations of sparse coding.
\newblock In \emph{Proceedings of the 27th international conference on
  international conference on machine learning}, pages 399--406, 2010.

\bibitem[Andrychowicz et~al.(2016)Andrychowicz, Denil, Gomez, Hoffman, Pfau,
  Schaul, Shillingford, and De~Freitas]{andrychowicz2016learning}
Marcin Andrychowicz, Misha Denil, Sergio Gomez, Matthew~W Hoffman, David Pfau,
  Tom Schaul, Brendan Shillingford, and Nando De~Freitas.
\newblock Learning to learn by gradient descent by gradient descent.
\newblock In \emph{Advances in neural information processing systems}, pages
  3981--3989, 2016.

\bibitem[Chen et~al.(2017)Chen, Hoffman, Colmenarejo, Denil, Lillicrap,
  Botvinick, and Freitas]{chen2017learning}
Yutian Chen, Matthew~W Hoffman, Sergio~G{\'o}mez Colmenarejo, Misha Denil,
  Timothy~P Lillicrap, Matt Botvinick, and Nando Freitas.
\newblock Learning to learn without gradient descent by gradient descent.
\newblock In \emph{International Conference on Machine Learning}, pages
  748--756, 2017.

\bibitem[Khalil et~al.(2017)Khalil, Dai, Zhang, Dilkina, and
  Song]{khalil2017learning}
Elias Khalil, Hanjun Dai, Yuyu Zhang, Bistra Dilkina, and Le~Song.
\newblock Learning combinatorial optimization algorithms over graphs.
\newblock \emph{Advances in neural information processing systems},
  30:\penalty0 6348--6358, 2017.

\bibitem[Wang et~al.(2016{\natexlab{a}})Wang, Liu, Chang, Ling, Yang, and
  Huang]{Wang_Liu_Chang_Ling_Yang_Huang_2016}
Zhangyang Wang, Ding Liu, Shiyu Chang, Qing Ling, Yingzhen Yang, and Thomas~S
  Huang.
\newblock D3: Deep dual-domain based fast restoration of jpeg-compressed
  images.
\newblock In \emph{Proceedings of the IEEE Conference on Computer Vision and
  Pattern Recognition}, pages 2764--2772, 2016{\natexlab{a}}.

\bibitem[Zhang and Ghanem(2018)]{zhang2018ista}
Jian Zhang and Bernard Ghanem.
\newblock Ista-net: Interpretable optimization-inspired deep network for image
  compressive sensing.
\newblock In \emph{Proceedings of the IEEE conference on computer vision and
  pattern recognition}, pages 1828--1837, 2018.

\bibitem[Corbineau et~al.(2019)Corbineau, Bertocchi, Chouzenoux, Prato, and
  Pesquet]{Corbineau_Bertocchi_Chouzenoux_Prato_Pesquet_2019}
M.-C. Corbineau, C.~Bertocchi, E.~Chouzenoux, M.~Prato, and J.-C. Pesquet.
\newblock Learned image deblurring by unfolding a proximal interior point
  algorithm.
\newblock In \emph{2019 IEEE International Conference on Image Processing
  (ICIP)}, page 4664–4668. IEEE, Sep 2019.
\newblock ISBN 978-1-5386-6249-6.
\newblock \doi{10.1109/ICIP.2019.8803438}.

\bibitem[Liang et~al.(2020)Liang, Cheng, Ke, and Ying]{liang2020deep}
Dong Liang, Jing Cheng, Ziwen Ke, and Leslie Ying.
\newblock Deep magnetic resonance image reconstruction: Inverse problems meet
  neural networks.
\newblock \emph{IEEE Signal Processing Magazine}, 37\penalty0 (1):\penalty0
  141--151, 2020.

\bibitem[Yin et~al.(2021)Yin, Wu, Sun, Dalca, Yue, and Bouman]{yin2021end}
Tianwei Yin, Zihui Wu, He~Sun, Adrian~V Dalca, Yisong Yue, and Katherine~L
  Bouman.
\newblock End-to-end sequential sampling and reconstruction for mr imaging.
\newblock \emph{arXiv preprint arXiv:2105.06460}, 2021.

\bibitem[Borgerding et~al.(2017)Borgerding, Schniter, and
  Rangan]{Borgerding_Schniter_Rangan_2017}
Mark Borgerding, Philip Schniter, and Sundeep Rangan.
\newblock Amp-inspired deep networks for sparse linear inverse problems.
\newblock \emph{IEEE Transactions on Signal Processing}, 65\penalty0
  (16):\penalty0 4293–4308, Aug 2017.
\newblock ISSN 1941-0476.
\newblock \doi{10.1109/TSP.2017.2708040}.

\bibitem[Balatsoukas-Stimming and Studer(2019)]{balatsoukas2019deep}
Alexios Balatsoukas-Stimming and Christoph Studer.
\newblock Deep unfolding for communications systems: A survey and some new
  directions.
\newblock In \emph{2019 IEEE International Workshop on Signal Processing
  Systems (SiPS)}, pages 266--271. IEEE, 2019.

\bibitem[Marino et~al.(2020)Marino, Pich{\'e}, Ialongo, and
  Yue]{marino2020iterative}
Joseph Marino, Alexandre Pich{\'e}, Alessandro~Davide Ialongo, and Yisong Yue.
\newblock Iterative amortized policy optimization.
\newblock \emph{arXiv preprint arXiv:2010.10670}, 2020.

\bibitem[Vadori et~al.(2020)Vadori, Ganesh, Reddy, and
  Veloso]{vadori2020calibration}
Nelson Vadori, Sumitra Ganesh, Prashant Reddy, and Manuela Veloso.
\newblock Calibration of shared equilibria in general sum partially observable
  markov games.
\newblock \emph{Advances in Neural Information Processing Systems}, 33, 2020.

\bibitem[Vadori et~al.(2021)Vadori, Savani, Spooner, and
  Ganesh]{vadori2021consensus}
Nelson Vadori, Rahul Savani, Thomas Spooner, and Sumitra Ganesh.
\newblock Consensus multiplicative weights update: Learning to learn using
  projector-based game signatures.
\newblock \emph{arXiv preprint arXiv:2106.02615}, 2021.

\bibitem[Cao et~al.(2019)Cao, Chen, Wang, and Shen]{cao2019learning}
Yue Cao, Tianlong Chen, Zhangyang Wang, and Yang Shen.
\newblock Learning to optimize in swarms.
\newblock In \emph{Advances in Neural Information Processing Systems}, pages
  15018--15028, 2019.

\bibitem[Chen et~al.(2019{\natexlab{a}})Chen, Li, Umarov, Gao, and
  Song]{chen2019rna}
Xinshi Chen, Yu~Li, Ramzan Umarov, Xin Gao, and Le~Song.
\newblock Rna secondary structure prediction by learning unrolled algorithms.
\newblock In \emph{International Conference on Learning Representations},
  2019{\natexlab{a}}.

\bibitem[Agrawal et~al.(2020)Agrawal, Menzies, Minku, Wagner, and
  Yu]{agrawal2020better}
Amritanshu Agrawal, Tim Menzies, Leandro~L Minku, Markus Wagner, and Zhe Yu.
\newblock Better software analytics via “duo”: Data mining algorithms
  using/used-by optimizers.
\newblock \emph{Empirical Software Engineering}, 25\penalty0 (3):\penalty0
  2099--2136, 2020.

\bibitem[Li et~al.(2020{\natexlab{a}})Li, Chen, You, Wang, and
  Lin]{li_2020_halo}
Chaojian Li, Tianlong Chen, Haoran You, Zhangyang Wang, and Yingyan Lin.
\newblock Halo: Hardware-aware learning to optimize.
\newblock In \emph{Proceedings of the European Conference on Computer Vision
  (ECCV)}, September 2020{\natexlab{a}}.

\bibitem[Chen et~al.(2020{\natexlab{a}})Chen, Yu, Wang, and
  Anandkumar]{chen2020automated}
Wuyang Chen, Zhiding Yu, Zhangyang Wang, and Anima Anandkumar.
\newblock Automated synthetic-to-real generalization.
\newblock \emph{International Conference on Machine Learning (ICML)},
  2020{\natexlab{a}}.

\bibitem[Zhou et~al.(2018)Zhou, Di, Du, Peng, Yang, Pan, Tsang, Liu, Qin, and
  Goh]{zhou2018sc2net}
Joey~Tianyi Zhou, Kai Di, Jiawei Du, Xi~Peng, Hao Yang, Sinno~Jialin Pan,
  Ivor~W Tsang, Yong Liu, Zheng Qin, and Rick Siow~Mong Goh.
\newblock Sc2net: Sparse lstms for sparse coding.
\newblock In \emph{AAAI}, pages 4588--4595, 2018.

\bibitem[Heaton et~al.(2020{\natexlab{a}})Heaton, Chen, Wang, and
  Yin]{heaton2020safeguarded}
Howard Heaton, Xiaohan Chen, Zhangyang Wang, and Wotao Yin.
\newblock Safeguarded learned convex optimization.
\newblock \emph{arXiv preprint arXiv:2003.01880}, 2020{\natexlab{a}}.

\bibitem[Vilalta and Drissi(2002)]{vilalta2002perspective}
Ricardo Vilalta and Youssef Drissi.
\newblock A perspective view and survey of meta-learning.
\newblock \emph{Artificial intelligence review}, 18\penalty0 (2):\penalty0
  77--95, 2002.

\bibitem[Hospedales et~al.(2020)Hospedales, Antoniou, Micaelli, and
  Storkey]{hospedales2020meta}
Timothy Hospedales, Antreas Antoniou, Paul Micaelli, and Amos Storkey.
\newblock Meta-learning in neural networks: A survey.
\newblock \emph{arXiv preprint arXiv:2004.05439}, 2020.

\bibitem[Chen et~al.(2020{\natexlab{b}})Chen, Tang, Muandet,
  et~al.]{chen2020mate}
Xiaohan Chen, Siyu Tang, Krikamol Muandet, et~al.
\newblock Mate: Plugging in model awareness to task embedding for meta
  learning.
\newblock In \emph{Neural Information Processing Systems 2020},
  2020{\natexlab{b}}.

\bibitem[Li and Malik(2017{\natexlab{a}})]{li2017learning}
Ke~Li and Jitendra Malik.
\newblock Learning to optimize neural nets.
\newblock \emph{arXiv preprint arXiv:1703.00441}, 2017{\natexlab{a}}.

\bibitem[Yao et~al.(2018)Yao, Wang, Chen, Dai, Yi-Qi, Yu-Feng, Wei-Wei, Qiang,
  and Yang]{yao2018taking}
Quanming Yao, Mengshuo Wang, Yuqiang Chen, Wenyuan Dai, Hu~Yi-Qi, Li~Yu-Feng,
  Tu~Wei-Wei, Yang Qiang, and Yu~Yang.
\newblock Taking human out of learning applications: A survey on automated
  machine learning.
\newblock \emph{arXiv preprint arXiv:1810.13306}, 2018.

\bibitem[Elsken et~al.(2018)Elsken, Metzen, and Hutter]{elsken2018neural}
Thomas Elsken, Jan~Hendrik Metzen, and Frank Hutter.
\newblock Neural architecture search: A survey.
\newblock \emph{arXiv preprint arXiv:1808.05377}, 2018.

\bibitem[Hutter et~al.(2011)Hutter, Hoos, and
  Leyton-Brown]{hutter2011sequential}
Frank Hutter, Holger~H Hoos, and Kevin Leyton-Brown.
\newblock Sequential model-based optimization for general algorithm
  configuration.
\newblock In \emph{International conference on learning and intelligent
  optimization}, pages 507--523. Springer, 2011.

\bibitem[Feurer and Hutter(2019)]{feurer2019hyperparameter}
Matthias Feurer and Frank Hutter.
\newblock Hyperparameter optimization.
\newblock In \emph{Automated Machine Learning}, pages 3--33. Springer, Cham,
  2019.

\bibitem[Klein et~al.(2017)Klein, Falkner, Bartels, Hennig, and
  Hutter]{klein2017fast}
Aaron Klein, Stefan Falkner, Simon Bartels, Philipp Hennig, and Frank Hutter.
\newblock Fast bayesian optimization of machine learning hyperparameters on
  large datasets.
\newblock In \emph{Artificial Intelligence and Statistics}, pages 528--536.
  PMLR, 2017.

\bibitem[Yu et~al.(2016)Yu, Qian, and Hu]{yu2016derivative}
Yang Yu, Hong Qian, and Yi-Qi Hu.
\newblock Derivative-free optimization via classification.
\newblock In \emph{AAAI}, volume~16, pages 2286--2292, 2016.

\bibitem[Thornton et~al.(2013)Thornton, Hutter, Hoos, and
  Leyton-Brown]{thornton2013auto}
Chris Thornton, Frank Hutter, Holger~H Hoos, and Kevin Leyton-Brown.
\newblock Auto-weka: Combined selection and hyperparameter optimization of
  classification algorithms.
\newblock In \emph{Proceedings of the 19th ACM SIGKDD international conference
  on Knowledge discovery and data mining}, pages 847--855, 2013.

\bibitem[Xu et~al.(2019)Xu, Dai, Kemp, and Metz]{xu2019learning}
Zhen Xu, Andrew~M Dai, Jonas Kemp, and Luke Metz.
\newblock Learning an adaptive learning rate schedule.
\newblock \emph{arXiv preprint arXiv:1909.09712}, 2019.

\bibitem[Wang et~al.(2020{\natexlab{a}})Wang, Yuan, Wu, and
  Ge]{wang2020guarantees}
Xiang Wang, Shuai Yuan, Chenwei Wu, and Rong Ge.
\newblock Guarantees for tuning the step size using a learning-to-learn
  approach.
\newblock \emph{arXiv preprint arXiv:2006.16495}, 2020{\natexlab{a}}.

\bibitem[MIT()]{MIT}
Mit eecs 6.890 learning-augmented algorithms.
\newblock
  \url{https://www.eecs.mit.edu/academics-admissions/academic-information/subject-updates-spring-2019/6890}.
\newblock Spring 2019.

\bibitem[Kraska et~al.(2018)Kraska, Beutel, Chi, Dean, and
  Polyzotis]{kraska2018case}
Tim Kraska, Alex Beutel, Ed~H Chi, Jeffrey Dean, and Neoklis Polyzotis.
\newblock The case for learned index structures.
\newblock In \emph{Proceedings of the 2018 International Conference on
  Management of Data}, pages 489--504, 2018.

\bibitem[Mitzenmacher(2018)]{mitzenmacher2018model}
Michael Mitzenmacher.
\newblock A model for learned bloom filters and optimizing by sandwiching.
\newblock In \emph{Advances in Neural Information Processing Systems}, pages
  464--473, 2018.

\bibitem[Jiang et~al.(2019)Jiang, Li, Lin, Ruan, and
  Woodruff]{jiang2019learning}
Tanqiu Jiang, Yi~Li, Honghao Lin, Yisong Ruan, and David~P Woodruff.
\newblock Learning-augmented data stream algorithms.
\newblock In \emph{International Conference on Learning Representations}, 2019.

\bibitem[Hsu et~al.(2019)Hsu, Indyk, Katabi, and Vakilian]{hsu2019learning}
Chen-Yu Hsu, Piotr Indyk, Dina Katabi, and Ali Vakilian.
\newblock Learning-based frequency estimation algorithms.
\newblock In \emph{International Conference on Learning Representations}, 2019.

\bibitem[Foster et~al.(2018)Foster, Rakhlin, and Sridharan]{foster2018online}
Dylan~J Foster, Alexander Rakhlin, and Karthik Sridharan.
\newblock Online learning: Sufficient statistics and the burkholder method.
\newblock In \emph{Conference On Learning Theory}, pages 3028--3064, 2018.

\bibitem[Kim et~al.(2018)Kim, Jiang, Rana, Kannan, Oh, and
  Viswanath]{kim2018communication}
Hyeji Kim, Yihan Jiang, Ranvir Rana, Sreeram Kannan, Sewoong Oh, and Pramod
  Viswanath.
\newblock Communication algorithms via deep learning.
\newblock In \emph{6th International Conference on Learning Representations,
  ICLR 2018}, 2018.

\bibitem[Kim et~al.(2020)Kim, Jiang, Kannan, Oh, and
  Viswanath]{kim2020deepcode}
Hyeji Kim, Yihan Jiang, Sreeram Kannan, Sewoong Oh, and Pramod Viswanath.
\newblock Deepcode: Feedback codes via deep learning.
\newblock \emph{IEEE Journal on Selected Areas in Information Theory},
  1\penalty0 (1):\penalty0 194--206, 2020.

\bibitem[Mitzenmacher(2020)]{mitzenmacher2020scheduling}
Michael Mitzenmacher.
\newblock Scheduling with predictions and the price of misprediction.
\newblock In \emph{11th Innovations in Theoretical Computer Science Conference
  (ITCS 2020)}. Schloss Dagstuhl-Leibniz-Zentrum f{\"u}r Informatik, 2020.

\bibitem[Indyk et~al.(2019)Indyk, Vakilian, and Yuan]{indyk2019learning}
Piotr Indyk, Ali Vakilian, and Yang Yuan.
\newblock Learning-based low-rank approximations.
\newblock In \emph{Advances in Neural Information Processing Systems}, pages
  7402--7412, 2019.

\bibitem[Hutter et~al.(2019)Hutter, Kotthoff, and
  Vanschoren]{hutter2019automated}
Frank Hutter, Lars Kotthoff, and Joaquin Vanschoren.
\newblock \emph{Automated machine learning: methods, systems, challenges}.
\newblock Springer Nature, 2019.

\bibitem[He et~al.(2021)He, Zhao, and Chu]{he2021automl}
Xin He, Kaiyong Zhao, and Xiaowen Chu.
\newblock Automl: A survey of the state-of-the-art.
\newblock \emph{Knowledge-Based Systems}, 212:\penalty0 106622, 2021.

\bibitem[Monga et~al.(2019)Monga, Li, and Eldar]{monga2019algorithm}
Vishal Monga, Yuelong Li, and Yonina~C Eldar.
\newblock Algorithm unrolling: Interpretable, efficient deep learning for
  signal and image processing.
\newblock \emph{arXiv preprint arXiv:1912.10557}, 2019.

\bibitem[Shlezinger et~al.(2020)Shlezinger, Whang, Eldar, and
  Dimakis]{shlezinger2020model}
Nir Shlezinger, Jay Whang, Yonina~C Eldar, and Alexandros~G Dimakis.
\newblock Model-based deep learning.
\newblock \emph{arXiv preprint arXiv:2012.08405}, 2020.

\bibitem[Li et~al.(2018)Li, Chen, and
  Koltun]{LiChenQifengKoltun2018_Combinatorial}
Zhuwen Li, Qifeng Chen, and Vladlen Koltun.
\newblock Combinatorial optimization with graph convolutional networks and
  guided tree search.
\newblock In \emph{Advances in Neural Information Processing Systems}, pages
  539--548, 2018.

\bibitem[Dai et~al.(2018)Dai, Khalil, Zhang, Dilkina, and
  Song]{DaiKhalilZhangDilkinaSong2018_learning}
Hanjun Dai, Elias~B. Khalil, Yuyu Zhang, Bistra Dilkina, and Le~Song.
\newblock Learning {Combinatorial} {Optimization} {Algorithms} over {Graphs}.
\newblock \emph{arXiv:1704.01665}, February 2018.

\bibitem[Bertsimas and Stellato(2019)]{bertsimas2019online}
Dimitris Bertsimas and Bartolomeo Stellato.
\newblock Online mixed-integer optimization in milliseconds.
\newblock \emph{arXiv preprint arXiv:1907.02206}, 2019.

\bibitem[Bertsimas and Stellato(2021)]{bertsimas2021voice}
Dimitris Bertsimas and Bartolomeo Stellato.
\newblock The voice of optimization.
\newblock \emph{Machine Learning}, 110\penalty0 (2):\penalty0 249--277, 2021.

\bibitem[Cauligi et~al.(2020)Cauligi, Culbertson, Stellato, Bertsimas,
  Schwager, and Pavone]{cauligi2020learning}
Abhishek Cauligi, Preston Culbertson, Bartolomeo Stellato, Dimitris Bertsimas,
  Mac Schwager, and Marco Pavone.
\newblock Learning mixed-integer convex optimization strategies for robot
  planning and control.
\newblock In \emph{2020 59th IEEE Conference on Decision and Control (CDC)},
  pages 1698--1705. IEEE, 2020.

\bibitem[Bengio et~al.(2018)Bengio, Lodi, and
  Prouvost]{BengioLodiProuvost2018_Machine}
Yoshua Bengio, Andrea Lodi, and Antoine Prouvost.
\newblock Machine learning for combinatorial optimization: a methodological
  tour d'horizon.
\newblock \emph{arXiv:1811.06128}, 2018.

\bibitem[Lv et~al.(2017)Lv, Jiang, and Li]{lv2017learning}
Kaifeng Lv, Shunhua Jiang, and Jian Li.
\newblock Learning gradient descent: Better generalization and longer horizons.
\newblock \emph{arXiv preprint arXiv:1703.03633}, 2017.

\bibitem[Wichrowska et~al.(2017)Wichrowska, Maheswaranathan, Hoffman,
  Colmenarejo, Denil, de~Freitas, and Sohl-Dickstein]{wichrowska2017learned}
Olga Wichrowska, Niru Maheswaranathan, Matthew~W. Hoffman, Sergio~Gomez
  Colmenarejo, Misha Denil, Nando de~Freitas, and Jascha Sohl-Dickstein.
\newblock Learned optimizers that scale and generalize, 2017.

\bibitem[Metz et~al.(2019{\natexlab{a}})Metz, Maheswaranathan, Nixon, Freeman,
  and Sohl-Dickstein]{pmlr-v97-metz19a}
Luke Metz, Niru Maheswaranathan, Jeremy Nixon, Daniel Freeman, and Jascha
  Sohl-Dickstein.
\newblock Understanding and correcting pathologies in the training of learned
  optimizers.
\newblock In \emph{International Conference on Machine Learning}, pages
  4556--4565, 2019{\natexlab{a}}.

\bibitem[Li and Malik(2017{\natexlab{b}})]{li2016learning}
Ke~Li and Jitendra Malik.
\newblock Learning to optimize.
\newblock In \emph{International Conference on Learning Representations
  (ICLR)}, 2017{\natexlab{b}}.

\bibitem[Bello et~al.(2017)Bello, Zoph, Vasudevan, and Le]{pmlr-v70-bello17a}
Irwan Bello, Barret Zoph, Vijay Vasudevan, and Quoc~V. Le.
\newblock Neural optimizer search with reinforcement learning.
\newblock In Doina Precup and Yee~Whye Teh, editors, \emph{International
  Conference on Machine Learning}, volume~70 of \emph{Proceedings of Machine
  Learning Research}, pages 459--468, International Convention Centre, Sydney,
  Australia, 06--11 Aug 2017. PMLR.
\newblock URL \url{http://proceedings.mlr.press/v70/bello17a.html}.

\bibitem[Jiang et~al.(2018)Jiang, Chen, Shi, Dai, and Zhao]{jiang2018learning}
Haoming Jiang, Zhehui Chen, Yuyang Shi, Bo~Dai, and Tuo Zhao.
\newblock Learning to defense by learning to attack, 2018.

\bibitem[Xiong and Hsieh(2020)]{xiong2020improved}
Yuanhao Xiong and Cho-Jui Hsieh.
\newblock Improved adversarial training via learned optimizer, 2020.

\bibitem[Leclerc and Madry(2020)]{leclerc2020two}
Guillaume Leclerc and Aleksander Madry.
\newblock The two regimes of deep network training.
\newblock \emph{arXiv preprint arXiv:2002.10376}, 2020.

\bibitem[Lewkowycz et~al.(2020)Lewkowycz, Bahri, Dyer, Sohl-Dickstein, and
  Gur-Ari]{lewkowycz2020large}
Aitor Lewkowycz, Yasaman Bahri, Ethan Dyer, Jascha Sohl-Dickstein, and Guy
  Gur-Ari.
\newblock The large learning rate phase of deep learning: the catapult
  mechanism.
\newblock \emph{arXiv preprint arXiv:2003.02218}, 2020.

\bibitem[Wu et~al.(2020)Wu, Zou, Braverman, and Gu]{wu2020direction}
Jingfeng Wu, Difan Zou, Vladimir Braverman, and Quanquan Gu.
\newblock Direction matters: On the implicit regularization effect of
  stochastic gradient descent with moderate learning rate.
\newblock \emph{arXiv preprint arXiv:2011.02538}, 2020.

\bibitem[Chen et~al.(2020{\natexlab{c}})Chen, Dai, Li, Gao, and
  Song]{chen2020learning2stop}
Xinshi Chen, Hanjun Dai, Yu~Li, Xin Gao, and Le~Song.
\newblock Learning to stop while learning to predict.
\newblock \emph{arXiv preprint arXiv:2006.05082}, 2020{\natexlab{c}}.

\bibitem[Kingma and Welling(2013)]{kingma2013auto}
Diederik~P Kingma and Max Welling.
\newblock Auto-encoding variational bayes.
\newblock \emph{arXiv preprint arXiv:1312.6114}, 2013.

\bibitem[Wierstra et~al.(2008)Wierstra, Schaul, Peters, and
  Schmidhuber]{wierstra2008natural}
Daan Wierstra, Tom Schaul, Jan Peters, and Juergen Schmidhuber.
\newblock Natural evolution strategies.
\newblock In \emph{2008 IEEE Congress on Evolutionary Computation (IEEE World
  Congress on Computational Intelligence)}, pages 3381--3387. IEEE, 2008.

\bibitem[Metz et~al.(2019{\natexlab{b}})Metz, Maheswaranathan, Shlens,
  Sohl-Dickstein, and Cubuk]{metz2019using}
Luke Metz, Niru Maheswaranathan, Jonathon Shlens, Jascha Sohl-Dickstein, and
  Ekin~D Cubuk.
\newblock Using learned optimizers to make models robust to input noise.
\newblock \emph{arXiv preprint arXiv:1906.03367}, 2019{\natexlab{b}}.

\bibitem[Parikh and Boyd(2014)]{parikh2014proximal}
Neal Parikh and Stephen Boyd.
\newblock Proximal algorithms.
\newblock \emph{Foundations and Trends in optimization}, 1\penalty0
  (3):\penalty0 127--239, 2014.

\bibitem[Chen et~al.(2020{\natexlab{d}})Chen, Zhang, Zhou, Chang, Liu, Amini,
  and Wang]{chen2020training}
Tianlong Chen, Weiyi Zhang, Jingyang Zhou, Shiyu Chang, Sijia Liu, Lisa Amini,
  and Zhangyang Wang.
\newblock Training stronger baselines for learning to optimize.
\newblock \emph{arXiv preprint arXiv:2010.09089}, 2020{\natexlab{d}}.

\bibitem[Almeida et~al.(2021)Almeida, Winter, Tang, and
  Zaremba]{almeida2021generalizable}
Diogo Almeida, Clemens Winter, Jie Tang, and Wojciech Zaremba.
\newblock A generalizable approach to learning optimizers.
\newblock \emph{arXiv preprint arXiv:2106.00958}, 2021.

\bibitem[Real et~al.(2020)Real, Liang, So, and Le]{real2020automl}
Esteban Real, Chen Liang, David So, and Quoc Le.
\newblock Automl-zero: evolving machine learning algorithms from scratch.
\newblock In \emph{International Conference on Machine Learning}, pages
  8007--8019. PMLR, 2020.

\bibitem[Snoek et~al.(2014)Snoek, Swersky, Zemel, and Adams]{snoek2014input}
Jasper Snoek, Kevin Swersky, Rich Zemel, and Ryan Adams.
\newblock Input warping for bayesian optimization of non-stationary functions.
\newblock In \emph{International Conference on Machine Learning}, pages
  1674--1682, 2014.

\bibitem[Bergstra et~al.(2011)Bergstra, Bardenet, Bengio, and
  K{\'e}gl]{bergstra2011algorithms}
James Bergstra, R{\'e}mi Bardenet, Yoshua Bengio, and Bal{\'a}zs K{\'e}gl.
\newblock Algorithms for hyper-parameter optimization.
\newblock \emph{Advances in neural information processing systems},
  24:\penalty0 2546--2554, 2011.

\bibitem[Moal and Bates(2010)]{moal2010swarmdock}
Iain~H Moal and Paul~A Bates.
\newblock Swarmdock and the use of normal modes in protein-protein docking.
\newblock \emph{International journal of molecular sciences}, 11\penalty0
  (10):\penalty0 3623--3648, 2010.

\bibitem[Goodfellow et~al.(2014)Goodfellow, Pouget-Abadie, Mirza, Xu,
  Warde-Farley, Ozair, Courville, and Bengio]{goodfellow2014generative}
Ian Goodfellow, Jean Pouget-Abadie, Mehdi Mirza, Bing Xu, David Warde-Farley,
  Sherjil Ozair, Aaron Courville, and Yoshua Bengio.
\newblock Generative adversarial nets.
\newblock In \emph{Advances in neural information processing systems}, pages
  2672--2680, 2014.

\bibitem[Madry et~al.(2018)Madry, Makelov, Schmidt, Tsipras, and
  Vladu]{madry2018towards}
Aleksander Madry, Aleksandar Makelov, Ludwig Schmidt, Dimitris Tsipras, and
  Adrian Vladu.
\newblock Towards deep learning models resistant to adversarial attacks.
\newblock In \emph{International Conference on Learning Representations}, 2018.

\bibitem[Wu et~al.(2018)Wu, Wang, Wang, and Jin]{wu2018towards}
Zhenyu Wu, Zhangyang Wang, Zhaowen Wang, and Hailin Jin.
\newblock Towards privacy-preserving visual recognition via adversarial
  training: A pilot study.
\newblock In \emph{Proceedings of the European Conference on Computer Vision
  (ECCV)}, pages 606--624, 2018.

\bibitem[Ruan et~al.(2019)Ruan, Xiong, Reddi, Kumar, and
  Hsieh]{ruan2019learning}
Yangjun Ruan, Yuanhao Xiong, Sashank Reddi, Sanjiv Kumar, and Cho-Jui Hsieh.
\newblock Learning to learn by zeroth-order oracle.
\newblock \emph{arXiv preprint arXiv:1910.09464}, 2019.

\bibitem[Shen et~al.(2021)Shen, Chen, Heaton, Chen, Liu, Yin, and
  Wang]{minimax}
Jiayi Shen, Xiaohan Chen, Howard Heaton, Tianlong Chen, Jialin Liu, Wotao Yin,
  and Zhangyang Wang.
\newblock Learning a minimax optimizer: A pilot study.
\newblock In \emph{International Conference on Learning Representations
  (ICLR)}, 2021.

\bibitem[Ravi and Larochelle(2016)]{ravi2016optimization}
Sachin Ravi and Hugo Larochelle.
\newblock Optimization as a model for few-shot learning.
\newblock In \emph{International Conference on Learning Representations
  (ICLR)}, 2016.

\bibitem[Li et~al.(2017)Li, Zhou, Chen, and Li]{li2017meta}
Zhenguo Li, Fengwei Zhou, Fei Chen, and Hang Li.
\newblock Meta-sgd: Learning to learn quickly for few-shot learning.
\newblock \emph{arXiv preprint arXiv:1707.09835}, 2017.

\bibitem[Venkatakrishnan et~al.(2013)Venkatakrishnan, Bouman, and
  Wohlberg]{venkatakrishnan2013plug}
Singanallur~V Venkatakrishnan, Charles~A Bouman, and Brendt Wohlberg.
\newblock Plug-and-play priors for model based reconstruction.
\newblock In \emph{2013 IEEE Global Conference on Signal and Information
  Processing}, pages 945--948. IEEE, 2013.

\bibitem[Rudin et~al.(1992)Rudin, Osher, and Fatemi]{rudin1992nonlinear}
Leonid~I Rudin, Stanley Osher, and Emad Fatemi.
\newblock Nonlinear total variation based noise removal algorithms.
\newblock \emph{Physica D: nonlinear phenomena}, 60\penalty0 (1-4):\penalty0
  259--268, 1992.

\bibitem[Ryu et~al.(2019)Ryu, Liu, Wang, Chen, Wang, and Yin]{ryu2019plug}
Ernest Ryu, Jialin Liu, Sicheng Wang, Xiaohan Chen, Zhangyang Wang, and Wotao
  Yin.
\newblock Plug-and-play methods provably converge with properly trained
  denoisers.
\newblock In \emph{International Conference on Machine Learning (ICML)}, pages
  5546--5557, 2019.

\bibitem[Dabov et~al.(2007)Dabov, Foi, Katkovnik, and
  Egiazarian]{dabov2007image}
Kostadin Dabov, Alessandro Foi, Vladimir Katkovnik, and Karen Egiazarian.
\newblock Image denoising by sparse 3-d transform-domain collaborative
  filtering.
\newblock \emph{IEEE Transactions on image processing}, 16\penalty0
  (8):\penalty0 2080--2095, 2007.

\bibitem[Aharon et~al.(2006)Aharon, Elad, and Bruckstein]{aharon2006k}
Michal Aharon, Michael Elad, and Alfred Bruckstein.
\newblock K-svd: An algorithm for designing overcomplete dictionaries for
  sparse representation.
\newblock \emph{IEEE Transactions on signal processing}, 54\penalty0
  (11):\penalty0 4311--4322, 2006.

\bibitem[Buades et~al.(2005)Buades, Coll, and Morel]{buades2005nonlocal}
Antoni Buades, Bartomeu Coll, and J-M Morel.
\newblock A non-local algorithm for image denoising.
\newblock In \emph{2005 IEEE Computer Society Conference on Computer Vision and
  Pattern Recognition (CVPR'05)}, volume~2, pages 60--65. IEEE, 2005.

\bibitem[Heide et~al.(2014)Heide, Steinberger, Tsai, Rouf, Paj{{a}}k, Reddy,
  Gallo, Liu, Heidrich, Egiazarian, et~al.]{heide2014flexisp}
Felix Heide, Markus Steinberger, Yun-Ta Tsai, Mushfiqur Rouf, Dawid Paj{{a}}k,
  Dikpal Reddy, Orazio Gallo, Jing Liu, Wolfgang Heidrich, Karen Egiazarian,
  et~al.
\newblock Flexisp: A flexible camera image processing framework.
\newblock \emph{ACM Transactions on Graphics (TOG)}, 33\penalty0 (6):\penalty0
  1--13, 2014.

\bibitem[Metzler et~al.(2015)Metzler, Maleki, and Baraniuk]{metzler2015bm3d}
Christopher~A Metzler, Arian Maleki, and Richard~G Baraniuk.
\newblock Bm3d-amp: A new image recovery algorithm based on bm3d denoising.
\newblock In \emph{2015 IEEE International Conference on Image Processing
  (ICIP)}, pages 3116--3120. IEEE, 2015.

\bibitem[Rond et~al.(2016)Rond, Giryes, and Elad]{rond2016poisson}
Arie Rond, Raja Giryes, and Michael Elad.
\newblock Poisson inverse problems by the plug-and-play scheme.
\newblock \emph{Journal of Visual Communication and Image Representation},
  41:\penalty0 96--108, 2016.

\bibitem[Brifman et~al.(2016)Brifman, Romano, and Elad]{brifman2016turning}
Alon Brifman, Yaniv Romano, and Michael Elad.
\newblock Turning a denoiser into a super-resolver using plug and play priors.
\newblock In \emph{2016 IEEE International Conference on Image Processing
  (ICIP)}, pages 1404--1408. IEEE, 2016.

\bibitem[Wang and Chan(2017)]{wang2017parameter}
Xiran Wang and Stanley~H Chan.
\newblock Parameter-free plug-and-play admm for image restoration.
\newblock In \emph{2017 IEEE International Conference on Acoustics, Speech and
  Signal Processing (ICASSP)}, pages 1323--1327. IEEE, 2017.

\bibitem[Ono(2017)]{ono2017primal}
Shunsuke Ono.
\newblock Primal-dual plug-and-play image restoration.
\newblock \emph{IEEE Signal Processing Letters}, 24\penalty0 (8):\penalty0
  1108--1112, 2017.

\bibitem[Kamilov et~al.(2017)Kamilov, Mansour, and Wohlberg]{kamilov2017plug}
Ulugbek~S Kamilov, Hassan Mansour, and Brendt Wohlberg.
\newblock A plug-and-play priors approach for solving nonlinear imaging inverse
  problems.
\newblock \emph{IEEE Signal Processing Letters}, 24\penalty0 (12):\penalty0
  1872--1876, 2017.

\bibitem[He et~al.(2018)He, Yang, Wang, Zeng, Bian, Zhang, Sun, Xu, and
  Ma]{he2018optimizing}
Ji~He, Yan Yang, Yongbo Wang, Dong Zeng, Zhaoying Bian, Hao Zhang, Jian Sun,
  Zongben Xu, and Jianhua Ma.
\newblock Optimizing a parameterized plug-and-play admm for iterative low-dose
  ct reconstruction.
\newblock \emph{IEEE transactions on medical imaging}, 38\penalty0
  (2):\penalty0 371--382, 2018.

\bibitem[Gupta et~al.(2018)Gupta, Jin, Nguyen, McCann, and Unser]{gupta2018cnn}
Harshit Gupta, Kyong~Hwan Jin, Ha~Q Nguyen, Michael~T McCann, and Michael
  Unser.
\newblock Cnn-based projected gradient descent for consistent ct image
  reconstruction.
\newblock \emph{IEEE transactions on medical imaging}, 37\penalty0
  (6):\penalty0 1440--1453, 2018.

\bibitem[Yang and Sun(2018)]{yang2018proximal}
Dong Yang and Jian Sun.
\newblock Proximal dehaze-net: A prior learning-based deep network for single
  image dehazing.
\newblock In \emph{Proceedings of the European Conference on Computer Vision
  (ECCV)}, pages 702--717, 2018.

\bibitem[Ye et~al.(2018)Ye, Srivastava, Thibault, Sauer, and
  Bouman]{ye2018deep}
Dong~Hye Ye, Somesh Srivastava, Jean-Baptiste Thibault, Ken Sauer, and Charles
  Bouman.
\newblock Deep residual learning for model-based iterative ct reconstruction
  using plug-and-play framework.
\newblock In \emph{2018 IEEE International Conference on Acoustics, Speech and
  Signal Processing (ICASSP)}, pages 6668--6672. IEEE, 2018.

\bibitem[Lyu et~al.(2019)Lyu, Ruan, Hoffman, Neph, McNitt-Gray, and
  Sheng]{lyu2019iterative}
Qihui Lyu, Dan Ruan, John Hoffman, Ryan Neph, Michael McNitt-Gray, and
  Ke~Sheng.
\newblock Iterative reconstruction for low dose ct using plug-and-play
  alternating direction method of multipliers (admm) framework.
\newblock In \emph{Medical Imaging 2019: Image Processing}, volume 10949, page
  1094906. International Society for Optics and Photonics, 2019.

\bibitem[Zhang et~al.(2019{\natexlab{a}})Zhang, Zuo, and Zhang]{zhang2019deep}
Kai Zhang, Wangmeng Zuo, and Lei Zhang.
\newblock Deep plug-and-play super-resolution for arbitrary blur kernels.
\newblock In \emph{Proceedings of the IEEE Conference on Computer Vision and
  Pattern Recognition}, pages 1671--1681, 2019{\natexlab{a}}.

\bibitem[Yuan et~al.(2020)Yuan, Liu, Suo, and Dai]{yuan2020plug}
Xin Yuan, Yang Liu, Jinli Suo, and Qionghai Dai.
\newblock Plug-and-play algorithms for large-scale snapshot compressive
  imaging.
\newblock In \emph{Proceedings of the IEEE/CVF Conference on Computer Vision
  and Pattern Recognition}, pages 1447--1457, 2020.

\bibitem[Ahmad et~al.(2020)Ahmad, Bouman, Buzzard, Chan, Liu, Reehorst, and
  Schniter]{ahmad2020plug}
Rizwan Ahmad, Charles~A Bouman, Gregery~T Buzzard, Stanley Chan, Sizhuo Liu,
  Edward~T Reehorst, and Philip Schniter.
\newblock Plug-and-play methods for magnetic resonance imaging: Using denoisers
  for image recovery.
\newblock \emph{IEEE Signal Processing Magazine}, 37\penalty0 (1):\penalty0
  105--116, 2020.

\bibitem[Mataev et~al.(2019)Mataev, Milanfar, and Elad]{mataev2019deepred}
Gary Mataev, Peyman Milanfar, and Michael Elad.
\newblock Deepred: Deep image prior powered by red.
\newblock In \emph{Proceedings of the IEEE/CVF International Conference on
  Computer Vision Workshops}, pages 0--0, 2019.

\bibitem[Song et~al.(2020)Song, Sun, Liu, Wang, and Kamilov]{song2020new}
Guangxiao Song, Yu~Sun, Jiaming Liu, Zhijie Wang, and Ulugbek~S Kamilov.
\newblock A new recurrent plug-and-play prior based on the multiple
  self-similarity network.
\newblock \emph{IEEE Signal Processing Letters}, 27:\penalty0 451--455, 2020.

\bibitem[Sreehari et~al.(2016)Sreehari, Venkatakrishnan, Wohlberg, Buzzard,
  Drummy, Simmons, and Bouman]{sreehari2016plug}
Suhas Sreehari, S~Venkat Venkatakrishnan, Brendt Wohlberg, Gregery~T Buzzard,
  Lawrence~F Drummy, Jeffrey~P Simmons, and Charles~A Bouman.
\newblock Plug-and-play priors for bright field electron tomography and sparse
  interpolation.
\newblock \emph{IEEE Transactions on Computational Imaging}, 2\penalty0
  (4):\penalty0 408--423, 2016.

\bibitem[Chan et~al.(2016)Chan, Wang, and Elgendy]{chan2016plug}
Stanley~H Chan, Xiran Wang, and Omar~A Elgendy.
\newblock Plug-and-play admm for image restoration: Fixed-point convergence and
  applications.
\newblock \emph{IEEE Transactions on Computational Imaging}, 3\penalty0
  (1):\penalty0 84--98, 2016.

\bibitem[Teodoro et~al.(2017)Teodoro, Bioucas-Dias, and
  Figueiredo]{teodoro2017scene}
Afonso~M Teodoro, Jos{\'e}~M Bioucas-Dias, and M{\'a}rio~AT Figueiredo.
\newblock Scene-adapted plug-and-play algorithm with convergence guarantees.
\newblock In \emph{2017 IEEE 27th International Workshop on Machine Learning
  for Signal Processing (MLSP)}, pages 1--6. IEEE, 2017.

\bibitem[Teodoro et~al.(2019)Teodoro, Bioucas-Dias, and
  Figueiredo]{teodoro2019image}
Afonso~M Teodoro, Jos{\'e}~M Bioucas-Dias, and M{\'a}rio~AT Figueiredo.
\newblock Image restoration and reconstruction using targeted plug-and-play
  priors.
\newblock \emph{IEEE Transactions on Computational Imaging}, 5\penalty0
  (4):\penalty0 675--686, 2019.

\bibitem[Buzzard et~al.(2018)Buzzard, Chan, Sreehari, and
  Bouman]{buzzard2018plug}
Gregery~T Buzzard, Stanley~H Chan, Suhas Sreehari, and Charles~A Bouman.
\newblock Plug-and-play unplugged: Optimization-free reconstruction using
  consensus equilibrium.
\newblock \emph{SIAM Journal on Imaging Sciences}, 11\penalty0 (3):\penalty0
  2001--2020, 2018.

\bibitem[Dong et~al.(2018)Dong, Wang, Yin, Shi, Wu, and Lu]{dong2018denoising}
Weisheng Dong, Peiyao Wang, Wotao Yin, Guangming Shi, Fangfang Wu, and Xiaotong
  Lu.
\newblock Denoising prior driven deep neural network for image restoration.
\newblock \emph{IEEE transactions on pattern analysis and machine
  intelligence}, 41\penalty0 (10):\penalty0 2305--2318, 2018.

\bibitem[Tirer and Giryes(2018)]{tirer2018image}
Tom Tirer and Raja Giryes.
\newblock Image restoration by iterative denoising and backward projections.
\newblock \emph{IEEE Transactions on Image Processing}, 28\penalty0
  (3):\penalty0 1220--1234, 2018.

\bibitem[Chan(2019)]{chan2019performance}
Stanley~H Chan.
\newblock Performance analysis of plug-and-play admm: A graph signal processing
  perspective.
\newblock \emph{IEEE Transactions on Computational Imaging}, 5\penalty0
  (2):\penalty0 274--286, 2019.

\bibitem[Sun et~al.(2019)Sun, Wohlberg, and Kamilov]{sun2019online}
Yu~Sun, Brendt Wohlberg, and Ulugbek~S Kamilov.
\newblock An online plug-and-play algorithm for regularized image
  reconstruction.
\newblock \emph{IEEE Transactions on Computational Imaging}, 5\penalty0
  (3):\penalty0 395--408, 2019.

\bibitem[Gavaskar and Chaudhury(2020)]{gavaskar2020plug}
Ruturaj~G Gavaskar and Kunal~N Chaudhury.
\newblock Plug-and-play ista converges with kernel denoisers.
\newblock \emph{IEEE Signal Processing Letters}, 27:\penalty0 610--614, 2020.

\bibitem[Xu et~al.(2020)Xu, Sun, Liu, Wohlberg, and Kamilov]{xu2020provable}
Xiaojian Xu, Yu~Sun, Jiaming Liu, Brendt Wohlberg, and Ulugbek~S Kamilov.
\newblock Provable convergence of plug-and-play priors with mmse denoisers.
\newblock \emph{IEEE Signal Processing Letters}, 27:\penalty0 1280--1284, 2020.

\bibitem[Sun et~al.(2020{\natexlab{a}})Sun, Wu, Wohlberg, and
  Kamilov]{sun2020scalable}
Yu~Sun, Zihui Wu, Brendt Wohlberg, and Ulugbek~S Kamilov.
\newblock Scalable plug-and-play admm with convergence guarantees.
\newblock \emph{arXiv preprint arXiv:2006.03224}, 2020{\natexlab{a}}.

\bibitem[Sun et~al.(2020{\natexlab{b}})Sun, Liu, Sun, Wohlberg, and
  Kamilov]{sun2020async}
Yu~Sun, Jiaming Liu, Yiran Sun, Brendt Wohlberg, and Ulugbek~S Kamilov.
\newblock Async-red: A provably convergent asynchronous block parallel
  stochastic method using deep denoising priors.
\newblock \emph{arXiv preprint arXiv:2010.01446}, 2020{\natexlab{b}}.

\bibitem[Meinhardt et~al.(2017)Meinhardt, Moller, Hazirbas, and
  Cremers]{meinhardt2017learning}
Tim Meinhardt, Michael Moller, Caner Hazirbas, and Daniel Cremers.
\newblock Learning proximal operators: Using denoising networks for
  regularizing inverse imaging problems.
\newblock In \emph{Proceedings of the IEEE International Conference on Computer
  Vision}, pages 1781--1790, 2017.

\bibitem[Rick~Chang et~al.(2017)Rick~Chang, Li, Poczos, Vijaya~Kumar, and
  Sankaranarayanan]{rick2017one}
JH~Rick~Chang, Chun-Liang Li, Barnabas Poczos, BVK Vijaya~Kumar, and Aswin~C
  Sankaranarayanan.
\newblock One network to solve them all--solving linear inverse problems using
  deep projection models.
\newblock In \emph{Proceedings of the IEEE International Conference on Computer
  Vision}, pages 5888--5897, 2017.

\bibitem[Zhang et~al.(2017)Zhang, Zuo, Gu, and Zhang]{zhang2017learning}
Kai Zhang, Wangmeng Zuo, Shuhang Gu, and Lei Zhang.
\newblock Learning deep cnn denoiser prior for image restoration.
\newblock In \emph{Proceedings of the IEEE conference on computer vision and
  pattern recognition}, pages 3929--3938, 2017.

\bibitem[Bigdeli et~al.(2017)Bigdeli, Zwicker, Favaro, and
  Jin]{bigdeli2017deep}
Siavash~Arjomand Bigdeli, Matthias Zwicker, Paolo Favaro, and Meiguang Jin.
\newblock Deep mean-shift priors for image restoration.
\newblock In \emph{Advances in Neural Information Processing Systems}, pages
  763--772, 2017.

\bibitem[Bigdeli. and Zwicker.(2018)]{visapp18}
Siavash~Arjomand Bigdeli. and Matthias Zwicker.
\newblock Image restoration using autoencoding priors.
\newblock In \emph{Proceedings of the 13th International Joint Conference on
  Computer Vision, Imaging and Computer Graphics Theory and Applications -
  Volume 5: VISAPP}, pages 33--44. INSTICC, SciTePress, 2018.
\newblock ISBN 978-989-758-290-5.
\newblock \doi{10.5220/0006532100330044}.

\bibitem[Lunz et~al.(2018)Lunz, {\"O}ktem, and
  Sch{\"o}nlieb]{lunz2018adversarial}
Sebastian Lunz, Ozan {\"O}ktem, and Carola-Bibiane Sch{\"o}nlieb.
\newblock Adversarial regularizers in inverse problems.
\newblock In \emph{Advances in Neural Information Processing Systems}, pages
  8507--8516, 2018.

\bibitem[Wei et~al.(2020)Wei, Aviles-Rivero, Liang, Fu, Schnlieb, and
  Huang]{wei2020tuning}
Kaixuan Wei, Angelica Aviles-Rivero, Jingwei Liang, Ying Fu, Carola-Bibiane
  Schnlieb, and Hua Huang.
\newblock Tuning-free plug-and-play proximal algorithm for inverse imaging
  problems.
\newblock \emph{arXiv preprint arXiv:2002.09611}, 2020.

\bibitem[Terris et~al.(2021)Terris, Repetti, Pesquet, and
  Wiaux]{terris2021enhanced}
Matthieu Terris, Audrey Repetti, Jean-Christophe Pesquet, and Yves Wiaux.
\newblock Enhanced convergent pnp algorithms for image restoration.
\newblock In \emph{IEEE ICIP 2021 Conference Proceedings}. IEEE, 2021.

\bibitem[Cohen et~al.(2020)Cohen, Elad, and Milanfar]{cohen2020regularization}
Regev Cohen, Michael Elad, and Peyman Milanfar.
\newblock Regularization by denoising via fixed-point projection (red-pro).
\newblock \emph{arXiv preprint arXiv:2008.00226}, 2020.

\bibitem[Heaton et~al.(2020{\natexlab{b}})Heaton, Fung, Lin, Osher, and
  Yin]{heaton2020projecting}
Howard Heaton, Samy~Wu Fung, Alex~Tong Lin, Stanley Osher, and Wotao Yin.
\newblock Projecting to manifolds via unsupervised learning.
\newblock \emph{arXiv preprint arXiv:2008.02200}, 2020{\natexlab{b}}.

\bibitem[Sprechmann et~al.(2013)Sprechmann, Litman, Ben~Yakar, Bronstein, and
  Sapiro]{sprechmann2013supervised}
Pablo Sprechmann, Roee Litman, Tal Ben~Yakar, Alexander~M Bronstein, and
  Guillermo Sapiro.
\newblock Supervised sparse analysis and synthesis operators.
\newblock \emph{Advances in Neural Information Processing Systems},
  26:\penalty0 908--916, 2013.

\bibitem[Moreau and Bruna(2017{\natexlab{a}})]{Moreau_Bruna_2017}
Thomas Moreau and Joan Bruna.
\newblock Understanding the learned iterative soft thresholding algorithm with
  matrix factorization.
\newblock In \emph{International Conference on Learning Representations
  (ICLR)}, 2017{\natexlab{a}}.
\newblock URL \url{http://arxiv.org/abs/1706.01338}.

\bibitem[Perdios et~al.(2017)Perdios, Besson, Rossinelli, and
  Thiran]{Perdios_Besson_Rossinelli_Thiran_2017}
Dimitris Perdios, Adrien Besson, Philippe Rossinelli, and Jean-Philippe Thiran.
\newblock Learning the weight matrix for sparsity averaging in compressive
  imaging.
\newblock In \emph{2017 IEEE International Conference on Image Processing
  (ICIP)}, pages 3056--3060. IEEE, 2017.

\bibitem[Giryes et~al.(2018)Giryes, Eldar, Bronstein, and
  Sapiro]{Giryes_Eldar_Bronstein_Sapiro_2018}
Raja Giryes, Yonina~C. Eldar, Alex~M. Bronstein, and Guillermo Sapiro.
\newblock Tradeoffs between convergence speed and reconstruction accuracy in
  inverse problems.
\newblock \emph{IEEE Transactions on Signal Processing}, 66\penalty0
  (7):\penalty0 1676–1690, Apr 2018.
\newblock ISSN 1941-0476.
\newblock \doi{10.1109/TSP.2018.2791945}.

\bibitem[Ablin et~al.(2019)Ablin, Moreau, Massias, and
  Gramfort]{ablin2019learning}
Pierre Ablin, Thomas Moreau, Mathurin Massias, and Alexandre Gramfort.
\newblock Learning step sizes for unfolded sparse coding.
\newblock In \emph{Advances in Neural Information Processing Systems}, pages
  13100--13110, 2019.

\bibitem[Hara et~al.(2019)Hara, Chen, Washio, Wazawa, and
  Nagai]{Hara_Chen_Washio_Wazawa_Nagai_2019}
Satoshi Hara, Weichih Chen, Takashi Washio, Tetsuichi Wazawa, and Takeharu
  Nagai.
\newblock Spod-net: Fast recovery of microscopic images using learned ista.
\newblock In \emph{Asian Conference on Machine Learning}, page 694–709, Oct
  2019.

\bibitem[Cowen et~al.(2019)Cowen, Saridena, and Choromanska]{cowen2019lsalsa}
Benjamin Cowen, Apoorva~Nandini Saridena, and Anna Choromanska.
\newblock Lsalsa: accelerated source separation via learned sparse coding.
\newblock \emph{Machine Learning}, 108\penalty0 (8-9):\penalty0 1307--1327,
  2019.

\bibitem[Wu et~al.(2019{\natexlab{a}})Wu, Dimakis, Sanghavi, Yu, Holtmann-Rice,
  Storcheus, Rostamizadeh, and
  Kumar]{Wu_Dimakis_Sanghavi_Yu_Holtmann-Rice_Storcheus_Rostamizadeh_Kumar_2019}
Shanshan Wu, Alexandros~G. Dimakis, Sujay Sanghavi, Felix~X. Yu, Daniel
  Holtmann-Rice, Dmitry Storcheus, Afshin Rostamizadeh, and Sanjiv Kumar.
\newblock Learning a compressed sensing measurement matrix via gradient
  unrolling.
\newblock \emph{arXiv:1806.10175 [cs, math, stat]}, Jul 2019{\natexlab{a}}.
\newblock URL \url{http://arxiv.org/abs/1806.10175}.
\newblock arXiv: 1806.10175.

\bibitem[Sprechmann et~al.(2012)Sprechmann, Bronstein, and
  Sapiro]{Bronstein_Sprechmann_Sapiro_2012}
Pablo Sprechmann, Alex Bronstein, and Guillermo Sapiro.
\newblock Learning efficient structured sparse models.
\newblock In \emph{International Conference on Machine Learning (ICML)}, pages
  615--622, 2012.

\bibitem[Papyan et~al.(2017)Papyan, Romano, and Elad]{papyan2017convolutional}
Vardan Papyan, Yaniv Romano, and Michael Elad.
\newblock Convolutional neural networks analyzed via convolutional sparse
  coding.
\newblock \emph{The Journal of Machine Learning Research}, 18\penalty0
  (1):\penalty0 2887--2938, 2017.

\bibitem[Aberdam et~al.(2019)Aberdam, Sulam, and Elad]{aberdam2019multi}
Aviad Aberdam, Jeremias Sulam, and Michael Elad.
\newblock Multi-layer sparse coding: The holistic way.
\newblock \emph{SIAM Journal on Mathematics of Data Science}, 1\penalty0
  (1):\penalty0 46--77, 2019.

\bibitem[Sulam et~al.(2019)Sulam, Aberdam, Beck, and Elad]{sulam2019multi}
Jeremias Sulam, Aviad Aberdam, Amir Beck, and Michael Elad.
\newblock On multi-layer basis pursuit, efficient algorithms and convolutional
  neural networks.
\newblock \emph{IEEE transactions on pattern analysis and machine
  intelligence}, 42\penalty0 (8):\penalty0 1968--1980, 2019.

\bibitem[Cherkaoui et~al.(2020)Cherkaoui, Sulam, and
  Moreau]{cherkaoui2020learning}
Hamza Cherkaoui, Jeremias Sulam, and Thomas Moreau.
\newblock Learning to solve tv regularised problems with unrolled algorithms.
\newblock \emph{Advances in Neural Information Processing Systems}, 33, 2020.

\bibitem[Mal{\'e}zieux et~al.(2021)Mal{\'e}zieux, Moreau, and
  Kowalski]{malezieux2021dictionary}
Beno{\^\i}t Mal{\'e}zieux, Thomas Moreau, and Matthieu Kowalski.
\newblock Dictionary and prior learning with unrolled algorithms for
  unsupervised inverse problems.
\newblock \emph{arXiv preprint arXiv:2106.06338}, 2021.

\bibitem[Wang et~al.(2016{\natexlab{b}})Wang, Ling, and
  Huang]{wang2016learning}
Zhangyang Wang, Qing Ling, and Thomas~S Huang.
\newblock Learning deep $\ell$ 0 encoders.
\newblock In \emph{Thirtieth AAAI Conference on Artificial Intelligence},
  2016{\natexlab{b}}.

\bibitem[Xin et~al.(2016)Xin, Wang, Gao, Wipf, and Wang]{xin2016maximal}
Bo~Xin, Yizhou Wang, Wen Gao, David Wipf, and Baoyuan Wang.
\newblock Maximal sparsity with deep networks?
\newblock In \emph{Advances in Neural Information Processing Systems}, pages
  4340--4348, 2016.

\bibitem[Wang et~al.(2016{\natexlab{c}})Wang, Yang, Chang, Ling, and
  Huang]{WangYangChangLingHuang2016_learning}
Zhangyang Wang, Yingzhen Yang, Shiyu Chang, Qing Ling, and Thomas~S. Huang.
\newblock Learning a {Deep} $l_\infty$ {Encoder} for {Hashing}.
\newblock In \emph{Proceedings of the {Twenty}-{Fifth} {International} {Joint}
  {Conference} on {Artificial} {Intelligence}}, {IJCAI}'16, pages 2174--2180.
  AAAI Press, 2016{\natexlab{c}}.

\bibitem[Studer et~al.(2014)Studer, Goldstein, Yin, and
  Baraniuk]{studer2014democratic}
Christoph Studer, Tom Goldstein, Wotao Yin, and Richard~G Baraniuk.
\newblock Democratic representations.
\newblock \emph{arXiv preprint arXiv:1401.3420}, 2014.

\bibitem[Hershey et~al.(2014)Hershey, Roux, and
  Weninger]{HersheyRouxWeninger2014_deep}
John~R. Hershey, Jonathan~Le Roux, and Felix Weninger.
\newblock Deep {Unfolding}: {Model}-{Based} {Inspiration} of {Novel} {Deep}
  {Architectures}.
\newblock \emph{arXiv:1409.2574}, 2014.

\bibitem[{Sprechmann} et~al.(2014){Sprechmann}, {Bronstein}, and
  {Sapiro}]{sprechmann2014supervised}
P.~{Sprechmann}, A.~M. {Bronstein}, and G.~{Sapiro}.
\newblock Supervised non-euclidean sparse nmf via bilevel optimization with
  applications to speech enhancement.
\newblock In \emph{2014 4th Joint Workshop on Hands-free Speech Communication
  and Microphone Arrays (HSCMA)}, pages 11--15, 2014.

\bibitem[Sprechmann et~al.(2015)Sprechmann, Bronstein, and
  Sapiro]{SprechmannBronsteinSapiro2015_learning}
P.~Sprechmann, A.~M. Bronstein, and G.~Sapiro.
\newblock Learning {Efficient} {Sparse} and {Low} {Rank} {Models}.
\newblock \emph{IEEE Transactions on Pattern Analysis and Machine
  Intelligence}, 37\penalty0 (9):\penalty0 1821--1833, September 2015.

\bibitem[Yakar et~al.(2013)Yakar, Litman, Sprechmann, Bronstein, and
  Sapiro]{yakar2013bilevel}
Tal~Ben Yakar, Roee Litman, Pablo Sprechmann, Alexander~M Bronstein, and
  Guillermo Sapiro.
\newblock Bilevel sparse models for polyphonic music transcription.
\newblock In \emph{ISMIR}, pages 65--70, 2013.

\bibitem[Zheng et~al.(2015)Zheng, Jayasumana, Romera-Paredes, Vineet, Su, Du,
  Huang, and Torr]{zheng2015conditional}
Shuai Zheng, Sadeep Jayasumana, Bernardino Romera-Paredes, Vibhav Vineet,
  Zhizhong Su, Dalong Du, Chang Huang, and Philip~HS Torr.
\newblock Conditional random fields as recurrent neural networks.
\newblock In \emph{Proceedings of the IEEE international conference on computer
  vision}, pages 1529--1537, 2015.

\bibitem[Chen and Pock(2017)]{Chen_Pock_2017}
Yunjin Chen and Thomas Pock.
\newblock Trainable nonlinear reaction diffusion: A flexible framework for fast
  and effective image restoration.
\newblock \emph{IEEE Transactions on Pattern Analysis and Machine
  Intelligence}, 39\penalty0 (6):\penalty0 1256–1272, Jun 2017.
\newblock ISSN 0162-8828, 2160-9292.
\newblock \doi{10.1109/TPAMI.2016.2596743}.

\bibitem[Long et~al.(2018)Long, Lu, Ma, and Dong]{long_pde-net_2018}
Zichao Long, Yiping Lu, Xianzhong Ma, and Bin Dong.
\newblock Pde-net: Learning pdes from data.
\newblock In \emph{International Conference on Machine Learning (ICML)}, pages
  3208--3216, 2018.

\bibitem[Long et~al.(2019{\natexlab{a}})Long, Lu, and Dong]{long_pde-net_2019}
Zichao Long, Yiping Lu, and Bin Dong.
\newblock {PDE}-{Net} 2.0: {Learning} {PDEs} from data with a numeric-symbolic
  hybrid deep network.
\newblock \emph{Journal of Computational Physics}, 399:\penalty0 108925,
  December 2019{\natexlab{a}}.
\newblock ISSN 0021-9991.
\newblock \doi{10.1016/j.jcp.2019.108925}.

\bibitem[Greenfeld et~al.(2019)Greenfeld, Galun, Kimmel, Yavneh, and
  Basri]{greenfeld_learning_2019}
Daniel Greenfeld, Meirav Galun, Ron Kimmel, Irad Yavneh, and Ronen Basri.
\newblock Learning to {Optimize} {Multigrid} {PDE} {Solvers}.
\newblock \emph{arXiv:1902.10248 [cs, math]}, August 2019.
\newblock URL \url{http://arxiv.org/abs/1902.10248}.
\newblock arXiv: 1902.10248.

\bibitem[Wiewel et~al.(2019)Wiewel, Becher, and Thuerey]{wiewel2019latent}
Steffen Wiewel, Moritz Becher, and Nils Thuerey.
\newblock Latent space physics: Towards learning the temporal evolution of
  fluid flow.
\newblock In \emph{Computer Graphics Forum}, volume~38, pages 71--82. Wiley
  Online Library, 2019.

\bibitem[Chen et~al.(2020{\natexlab{e}})Chen, Zhang, Reisinger, and
  Song]{Chen_Zhang_Reisinger_Song_2020}
Xinshi Chen, Yufei Zhang, Christoph Reisinger, and Le~Song.
\newblock Understanding deep architectures with reasoning layer.
\newblock \emph{arXiv:2006.13401 [cs, stat]}, Jun 2020{\natexlab{e}}.
\newblock URL \url{http://arxiv.org/abs/2006.13401}.
\newblock arXiv: 2006.13401.

\bibitem[Wang et~al.(2016{\natexlab{d}})Wang, Chang, Zhou, Wang, and
  Huang]{Wang_Chang_Zhou_Wang_Huang_2016}
Zhangyang Wang, Shiyu Chang, Jiayu Zhou, Meng Wang, and Thomas~S Huang.
\newblock Learning a task-specific deep architecture for clustering.
\newblock In \emph{Proceedings of the 2016 SIAM International Conference on
  Data Mining}, pages 369--377. SIAM, 2016{\natexlab{d}}.

\bibitem[Liu et~al.(2019{\natexlab{b}})Liu, Sun, Wang, Liu, and
  Zha]{liu2019frank}
Dong Liu, Ke~Sun, Zhangyang Wang, Runsheng Liu, and Zheng-Jun Zha.
\newblock Frank-wolfe network: An interpretable deep structure for non-sparse
  coding.
\newblock \emph{IEEE Transactions on Circuits and Systems for Video
  Technology}, 2019{\natexlab{b}}.

\bibitem[Pauwels et~al.(2021)Pauwels, Tsiligianni, and
  Deligiannis]{pauwels2021hcgm}
Ruben Pauwels, Evaggelia Tsiligianni, and Nikos Deligiannis.
\newblock Hcgm-net: A deep unfolding network for financial index tracking.
\newblock In \emph{ICASSP 2021-2021 IEEE International Conference on Acoustics,
  Speech and Signal Processing (ICASSP)}, pages 3910--3914. IEEE, 2021.

\bibitem[Bertocchi et~al.(2020)Bertocchi, Chouzenoux, Corbineau, Pesquet, and
  Prato]{Bertocchi_Chouzenoux_Corbineau_Pesquet_Prato_2020}
Carla Bertocchi, Emilie Chouzenoux, Marie-Caroline Corbineau, Jean-Christophe
  Pesquet, and Marco Prato.
\newblock Deep unfolding of a proximal interior point method for image
  restoration.
\newblock \emph{arXiv:1812.04276 [cs, math]}, Jan 2020.
\newblock URL \url{http://arxiv.org/abs/1812.04276}.
\newblock arXiv: 1812.04276.

\bibitem[Diamond et~al.(2018)Diamond, Sitzmann, Heide, and
  Wetzstein]{Diamond_Sitzmann_Heide_Wetzstein_2018}
Steven Diamond, Vincent Sitzmann, Felix Heide, and Gordon Wetzstein.
\newblock Unrolled optimization with deep priors.
\newblock \emph{arXiv:1705.08041 [cs]}, Dec 2018.
\newblock URL \url{http://arxiv.org/abs/1705.08041}.
\newblock arXiv: 1705.08041.

\bibitem[Takabe and Wadayama(2020)]{Takabe_Wadayama_2020}
Satoshi Takabe and Tadashi Wadayama.
\newblock Theoretical interpretation of learned step size in deep-unfolded
  gradient descent.
\newblock \emph{arXiv:2001.05142 [cs, math, stat]}, Jan 2020.
\newblock URL \url{http://arxiv.org/abs/2001.05142}.
\newblock arXiv: 2001.05142.

\bibitem[Adler and {\"O}ktem(2017)]{adler2017solving}
Jonas Adler and Ozan {\"O}ktem.
\newblock Solving ill-posed inverse problems using iterative deep neural
  networks.
\newblock \emph{Inverse Problems}, 33\penalty0 (12):\penalty0 124007, 2017.

\bibitem[Adler and {\"O}ktem(2018)]{adler2018learned}
Jonas Adler and Ozan {\"O}ktem.
\newblock Learned primal-dual reconstruction.
\newblock \emph{IEEE transactions on medical imaging}, 37\penalty0
  (6):\penalty0 1322--1332, 2018.

\bibitem[Domke(2012)]{domke2012generic}
Justin Domke.
\newblock Generic methods for optimization-based modeling.
\newblock In \emph{Artificial Intelligence and Statistics}, pages 318--326,
  2012.

\bibitem[Putzky and Welling(2017)]{Putzky_Welling_2017}
Patrick Putzky and Max Welling.
\newblock Recurrent inference machines for solving inverse problems.
\newblock \emph{arXiv:1706.04008 [cs]}, Jun 2017.
\newblock URL \url{http://arxiv.org/abs/1706.04008}.
\newblock arXiv: 1706.04008.

\bibitem[Mukherjee et~al.(2020)Mukherjee, Dittmer, Shumaylov, Lunz, {\"O}ktem,
  and Sch{\"o}nlieb]{mukherjee2020learned}
Subhadip Mukherjee, S{\"o}ren Dittmer, Zakhar Shumaylov, Sebastian Lunz, Ozan
  {\"O}ktem, and Carola-Bibiane Sch{\"o}nlieb.
\newblock Learned convex regularizers for inverse problems.
\newblock \emph{arXiv preprint arXiv:2008.02839}, 2020.

\bibitem[Wadayama and Takabe(2019)]{Wadayama_Takabe_2019}
Tadashi Wadayama and Satoshi Takabe.
\newblock Deep learning-aided trainable projected gradient decoding for ldpc
  codes.
\newblock \emph{arXiv:1901.04630 [cs, math]}, Jan 2019.
\newblock URL \url{http://arxiv.org/abs/1901.04630}.
\newblock arXiv: 1901.04630.

\bibitem[Kofler et~al.(2020)Kofler, Haltmeier, Schaeffter, Kachelrie{\ss},
  Dewey, Wald, and Kolbitsch]{kofler2020neural}
Andreas Kofler, Markus Haltmeier, Tobias Schaeffter, Marc Kachelrie{\ss}, Marc
  Dewey, Christian Wald, and Christoph Kolbitsch.
\newblock Neural networks-based regularization for large-scale medical image
  reconstruction.
\newblock \emph{Physics in Medicine \& Biology}, 65\penalty0 (13):\penalty0
  135003, 2020.

\bibitem[Blumensath and Davies(2008)]{blumensath2008iterative}
Thomas Blumensath and Mike~E Davies.
\newblock Iterative thresholding for sparse approximations.
\newblock \emph{Journal of Fourier analysis and Applications}, 14\penalty0
  (5-6):\penalty0 629--654, 2008.

\bibitem[Takabe et~al.(2020)Takabe, Wadayama, and Eldar]{takabe2020complex}
Satoshi Takabe, Tadashi Wadayama, and Yonina~C Eldar.
\newblock Complex trainable ista for linear and nonlinear inverse problems.
\newblock In \emph{IEEE International Conference on Acoustics, Speech and
  Signal Processing (ICASSP)}, pages 5020--5024, 2020.

\bibitem[Ito et~al.(2019)Ito, Takabe, and Wadayama]{Ito_Takabe_Wadayama_2019}
Daisuke Ito, Satoshi Takabe, and Tadashi Wadayama.
\newblock Trainable ista for sparse signal recovery.
\newblock \emph{IEEE Transactions on Signal Processing}, 67\penalty0
  (12):\penalty0 3113–3125, Jun 2019.
\newblock ISSN 1053-587X, 1941-0476.
\newblock \doi{10.1109/TSP.2019.2912879}.

\bibitem[Yao et~al.(2019)Yao, Dang, Zhang, and Wu]{Yao_Dang_Zhang_Wu_2019}
Mengcheng Yao, Jian Dang, Zaichen Zhang, and Liang Wu.
\newblock Sure-tista: A signal recovery network for compressed sensing.
\newblock In \emph{ICASSP 2019 - 2019 IEEE International Conference on
  Acoustics, Speech and Signal Processing (ICASSP)}, page 3832–3836. IEEE,
  May 2019.
\newblock ISBN 978-1-4799-8131-1.
\newblock \doi{10.1109/ICASSP.2019.8683182}.
\newblock URL \url{https://ieeexplore.ieee.org/document/8683182/}.

\bibitem[Aberdam et~al.(2020)Aberdam, Golts, and Elad]{Aberdam_Golts_Elad_2020}
Aviad Aberdam, Alona Golts, and Michael Elad.
\newblock Ada-lista: Learned solvers adaptive to varying models.
\newblock \emph{arXiv:2001.08456 [cs, stat]}, Feb 2020.
\newblock URL \url{http://arxiv.org/abs/2001.08456}.
\newblock arXiv: 2001.08456.

\bibitem[Behrens et~al.(2021)Behrens, Sauder, and
  Jung]{Behrens_Sauder_Jung_2020}
Freya Behrens, Jonathan Sauder, and Peter Jung.
\newblock Neurally augmented alista.
\newblock In \emph{International Conference on Learning Representations}, 2021.
\newblock URL \url{https://openreview.net/forum?id=q_S44KLQ_Aa}.

\bibitem[Tolooshams et~al.(2018)Tolooshams, Dey, and
  Ba]{Tolooshams_Dey_Ba_2018}
Bahareh Tolooshams, Sourav Dey, and Demba Ba.
\newblock Scalable convolutional dictionary learning with constrained recurrent
  sparse auto-encoders.
\newblock In \emph{2018 IEEE 28th International Workshop on Machine Learning
  for Signal Processing (MLSP)}, pages 1--6. IEEE, 2018.

\bibitem[Beck and Teboulle(2009)]{beck2009fast}
Amir Beck and Marc Teboulle.
\newblock A fast iterative shrinkage-thresholding algorithm for linear inverse
  problems.
\newblock \emph{SIAM journal on imaging sciences}, 2\penalty0 (1):\penalty0
  183--202, 2009.

\bibitem[Borgerding and Schniter(2016)]{Borgerding_Schniter_2016}
Mark Borgerding and Philip Schniter.
\newblock Onsager-corrected deep learning for sparse linear inverse problems.
\newblock In \emph{2016 IEEE Global Conference on Signal and Information
  Processing (GlobalSIP)}, page 227–231, Dec 2016.
\newblock \doi{10.1109/GlobalSIP.2016.7905837}.

\bibitem[Metzler et~al.(2017)Metzler, Mousavi, and
  Baraniuk]{metzler2017learned}
Chris Metzler, Ali Mousavi, and Richard Baraniuk.
\newblock Learned d-amp: Principled neural network based compressive image
  recovery.
\newblock In \emph{Advances in Neural Information Processing Systems}, pages
  1772--1783, 2017.

\bibitem[He et~al.(2020)He, Wen, Jin, and Li]{he2019model}
Hengtao He, Chao-Kai Wen, Shi Jin, and Geoffrey~Ye Li.
\newblock Model-driven deep learning for mimo detection.
\newblock \emph{IEEE Transactions on Signal Processing}, 68:\penalty0
  1702--1715, 2020.

\bibitem[Ma and Ping(2017)]{ma2017orthogonal}
Junjie Ma and Li~Ping.
\newblock Orthogonal amp.
\newblock \emph{IEEE Access}, 5:\penalty0 2020--2033, 2017.

\bibitem[Kobler et~al.(2017)Kobler, Klatzer, Hammernik, and
  Pock]{kobler2017variational}
Erich Kobler, Teresa Klatzer, Kerstin Hammernik, and Thomas Pock.
\newblock Variational networks: connecting variational methods and deep
  learning.
\newblock In \emph{German conference on pattern recognition}, pages 281--293.
  Springer, 2017.

\bibitem[Liu et~al.(2019{\natexlab{c}})Liu, Cheng, He, Fan, Lin, and
  Luo]{liu2019convergence}
Risheng Liu, Shichao Cheng, Yi~He, Xin Fan, Zhouchen Lin, and Zhongxuan Luo.
\newblock On the convergence of learning-based iterative methods for nonconvex
  inverse problems.
\newblock \emph{IEEE transactions on pattern analysis and machine
  intelligence}, 2019{\natexlab{c}}.

\bibitem[Sun et~al.(2016)Sun, Li, Xu, et~al.]{yang_Sun_Li_Xu_2016}
Jian Sun, Huibin Li, Zongben Xu, et~al.
\newblock Deep admm-net for compressive sensing mri.
\newblock In \emph{Advances in neural information processing systems}, pages
  10--18, 2016.

\bibitem[Xie et~al.(2019)Xie, Wu, Liu, Zhong, and Lin]{xie2019differentiable}
Xingyu Xie, Jianlong Wu, Guangcan Liu, Zhisheng Zhong, and Zhouchen Lin.
\newblock Differentiable linearized admm.
\newblock In \emph{International Conference on Machine Learning}, pages
  6902--6911, 2019.

\bibitem[Cheng et~al.(2019)Cheng, Wang, Ying, and Liang]{cheng2019model}
Jing Cheng, Haifeng Wang, Leslie Ying, and Dong Liang.
\newblock Model learning: Primal dual networks for fast mr imaging.
\newblock In \emph{International Conference on Medical Image Computing and
  Computer-Assisted Intervention}, pages 21--29. Springer, 2019.

\bibitem[Jiu and Pustelnik(2020)]{jiu2020deep}
Mingyuan Jiu and Nelly Pustelnik.
\newblock A deep primal-dual proximal network for image restoration.
\newblock \emph{arXiv preprint arXiv:2007.00959}, 2020.

\bibitem[Xiong and De~la Torre(2013)]{Xiong_De_la_Torre_2013}
Xuehan Xiong and Fernando De~la Torre.
\newblock Supervised descent method and its applications to face alignment.
\newblock In \emph{2013 IEEE Conference on Computer Vision and Pattern
  Recognition}, page 532–539. IEEE, Jun 2013.
\newblock ISBN 978-0-7695-4989-7.
\newblock \doi{10.1109/CVPR.2013.75}.
\newblock URL \url{http://ieeexplore.ieee.org/document/6618919/}.

\bibitem[Li et~al.()Li, Liu, and Yin]{lilearning}
Maojia Li, Jialin Liu, and Wotao Yin.
\newblock Learning to combine quasi-newton methods.

\bibitem[Long et~al.(2019{\natexlab{b}})Long, Lu, and Dong]{Long_Lu_Dong_2019}
Zichao Long, Yiping Lu, and Bin Dong.
\newblock Pde-net 2.0: Learning pdes from data with a numeric-symbolic hybrid
  deep network.
\newblock \emph{Journal of Computational Physics}, 399:\penalty0 108925, Dec
  2019{\natexlab{b}}.
\newblock ISSN 0021-9991.
\newblock \doi{10.1016/j.jcp.2019.108925}.

\bibitem[Zhang et~al.(2018)Zhang, Lu, Liu, and Dong]{Zhang_Lu_Liu_Dong_2018}
Xiaoshuai Zhang, Yiping Lu, Jiaying Liu, and Bin Dong.
\newblock Dynamically unfolding recurrent restorer: A moving endpoint control
  method for image restoration.
\newblock \emph{arXiv:1805.07709 [cs]}, Oct 2018.
\newblock URL \url{http://arxiv.org/abs/1805.07709}.
\newblock arXiv: 1805.07709.

\bibitem[Aggarwal et~al.(2018)Aggarwal, Mani, and Jacob]{aggarwal2018modl}
Hemant~K Aggarwal, Merry~P Mani, and Mathews Jacob.
\newblock Modl: Model-based deep learning architecture for inverse problems.
\newblock \emph{IEEE transactions on medical imaging}, 38\penalty0
  (2):\penalty0 394--405, 2018.

\bibitem[Khatib et~al.(2020)Khatib, Simon, and Elad]{khatib2020learned}
Rajaei Khatib, Dror Simon, and Michael Elad.
\newblock Learned greedy method (lgm): A novel neural architecture for sparse
  coding and beyond.
\newblock \emph{arXiv preprint arXiv:2010.07069}, 2020.

\bibitem[Huang et~al.(2020)Huang, Preuhs, Manhart, Lauritsch, and
  Maier]{huang2020data}
Yixing Huang, Alexander Preuhs, Michael Manhart, Guenter Lauritsch, and Andreas
  Maier.
\newblock Data consistent ct reconstruction from insufficient data with learned
  prior images.
\newblock \emph{arXiv preprint arXiv:2005.10034}, 2020.

\bibitem[Moeller et~al.(2019)Moeller, Mollenhoff, and
  Cremers]{moeller2019controlling}
Michael Moeller, Thomas Mollenhoff, and Daniel Cremers.
\newblock Controlling neural networks via energy dissipation.
\newblock In \emph{Proceedings of the IEEE International Conference on Computer
  Vision}, pages 3256--3265, 2019.

\bibitem[Hammernik et~al.(2017)Hammernik, W{\"u}rfl, Pock, and
  Maier]{hammernik2017deep}
Kerstin Hammernik, Tobias W{\"u}rfl, Thomas Pock, and Andreas Maier.
\newblock A deep learning architecture for limited-angle computed tomography
  reconstruction.
\newblock In \emph{Bildverarbeitung f{\"u}r die Medizin 2017}, pages 92--97.
  Springer, 2017.

\bibitem[Kobler et~al.(2020)Kobler, Effland, Kunisch, and
  Pock]{kobler2020total}
Erich Kobler, Alexander Effland, Karl Kunisch, and Thomas Pock.
\newblock Total deep variation: A stable regularizer for inverse problems.
\newblock \emph{arXiv preprint arXiv:2006.08789}, 2020.

\bibitem[Li et~al.(2020{\natexlab{b}})Li, Schwab, Antholzer, and
  Haltmeier]{li2020nett}
Housen Li, Johannes Schwab, Stephan Antholzer, and Markus Haltmeier.
\newblock Nett: Solving inverse problems with deep neural networks.
\newblock \emph{Inverse Problems}, 2020{\natexlab{b}}.

\bibitem[Egidio et~al.(2020)Egidio, Hansson, and Wahlberg]{egidio2020learning}
Lucas~N Egidio, Anders Hansson, and Bo~Wahlberg.
\newblock Learning the step-size policy for the limited-memory
  broyden-fletcher-goldfarb-shanno algorithm.
\newblock \emph{arXiv preprint arXiv:2010.01311}, 2020.

\bibitem[Shu et~al.(2020)Shu, Zhu, Zhao, Meng, and Xu]{shu2020meta}
Jun Shu, Yanwen Zhu, Qian Zhao, Deyu Meng, and Zongben Xu.
\newblock Meta-lr-schedule-net: Learned lr schedules that scale and generalize.
\newblock \emph{arXiv preprint arXiv:2007.14546}, 2020.

\bibitem[You et~al.(2020)You, Chen, Wang, and Shen]{you2020l2}
Yuning You, Tianlong Chen, Zhangyang Wang, and Yang Shen.
\newblock L2-gcn: Layer-wise and learned efficient training of graph
  convolutional networks.
\newblock In \emph{Proceedings of the IEEE/CVF Conference on Computer Vision
  and Pattern Recognition}, pages 2127--2135, 2020.

\bibitem[Wang et~al.(2019)Wang, Sun, and Xu]{wang2019hyperadam}
Shipeng Wang, Jian Sun, and Zongben Xu.
\newblock Hyperadam: A learnable task-adaptive adam for network training.
\newblock In \emph{Proceedings of the AAAI Conference on Artificial
  Intelligence}, volume~33, pages 5297--5304, 2019.

\bibitem[Chen et~al.(2020{\natexlab{f}})Chen, Chen, Chen, Yuan, Gong, Chen, and
  Wang]{chen2020self}
Xuxi Chen, Wuyang Chen, Tianlong Chen, Ye~Yuan, Chen Gong, Kewei Chen, and
  Zhangyang Wang.
\newblock Self-pu: Self boosted and calibrated positive-unlabeled training.
\newblock \emph{International Conference on Machine Learning (ICML)},
  2020{\natexlab{f}}.

\bibitem[Bora et~al.(2017)Bora, Jalal, Price, and Dimakis]{bora2017compressed}
Ashish Bora, Ajil Jalal, Eric Price, and Alexandros~G Dimakis.
\newblock Compressed sensing using generative models.
\newblock In \emph{International Conference on Machine Learning}, pages
  537--546, 2017.

\bibitem[Van~Veen et~al.(2018)Van~Veen, Jalal, Soltanolkotabi, Price,
  Vishwanath, and Dimakis]{van2018compressed}
Dave Van~Veen, Ajil Jalal, Mahdi Soltanolkotabi, Eric Price, Sriram Vishwanath,
  and Alexandros~G Dimakis.
\newblock Compressed sensing with deep image prior and learned regularization.
\newblock \emph{arXiv preprint arXiv:1806.06438}, 2018.

\bibitem[Sreter and Giryes(2018)]{sreter2018learned}
Hillel Sreter and Raja Giryes.
\newblock Learned convolutional sparse coding.
\newblock In \emph{2018 IEEE International Conference on Acoustics, Speech and
  Signal Processing (ICASSP)}, pages 2191--2195. IEEE, 2018.

\bibitem[Mardani et~al.(2019)Mardani, Sun, Papyan, Vasanawala, Pauly, and
  Donoho]{mardani2019degrees}
Morteza Mardani, Qingyun Sun, Vardan Papyan, Shreyas Vasanawala, John Pauly,
  and David Donoho.
\newblock Degrees of freedom analysis of unrolled neural networks.
\newblock \emph{arXiv preprint arXiv:1906.03742}, 2019.

\bibitem[Liu et~al.(2018{\natexlab{a}})Liu, Cheng, Ma, Fan, Luo,
  et~al.]{liu2018bridging}
Risheng Liu, Shichao Cheng, Long Ma, Xin Fan, Zhongxuan Luo, et~al.
\newblock A bridging framework for model optimization and deep propagation.
\newblock \emph{Advances in Neural Information Processing Systems},
  31:\penalty0 4318--4327, 2018{\natexlab{a}}.

\bibitem[Li et~al.(2020{\natexlab{c}})Li, Tofighi, Geng, Monga, and
  Eldar]{li2020efficient}
Yuelong Li, Mohammad Tofighi, Junyi Geng, Vishal Monga, and Yonina~C Eldar.
\newblock Efficient and interpretable deep blind image deblurring via algorithm
  unrolling.
\newblock \emph{IEEE Transactions on Computational Imaging}, 6:\penalty0
  666--681, 2020{\natexlab{c}}.

\bibitem[Mardani et~al.(2018)Mardani, Sun, Donoho, Papyan, Monajemi,
  Vasanawala, and Pauly]{mardani2018neural}
Morteza Mardani, Qingyun Sun, David Donoho, Vardan Papyan, Hatef Monajemi,
  Shreyas Vasanawala, and John Pauly.
\newblock Neural proximal gradient descent for compressive imaging.
\newblock \emph{Advances in Neural Information Processing Systems},
  31:\penalty0 9573--9583, 2018.

\bibitem[Fu et~al.(2019)Fu, Zha, Wu, Ding, and Paisley]{fu2019jpeg}
Xueyang Fu, Zheng-Jun Zha, Feng Wu, Xinghao Ding, and John Paisley.
\newblock Jpeg artifacts reduction via deep convolutional sparse coding.
\newblock In \emph{Proceedings of the IEEE International Conference on Computer
  Vision}, pages 2501--2510, 2019.

\bibitem[Wang et~al.(2020{\natexlab{b}})Wang, Xie, Zhao, and
  Meng]{wang2020model}
Hong Wang, Qi~Xie, Qian Zhao, and Deyu Meng.
\newblock A model-driven deep neural network for single image rain removal.
\newblock In \emph{Proceedings of the IEEE/CVF Conference on Computer Vision
  and Pattern Recognition}, pages 3103--3112, 2020{\natexlab{b}}.

\bibitem[Tamir et~al.(2019)Tamir, Stella, and Lustig]{tamir2019unsupervised}
Jonathan~I Tamir, X~Yu Stella, and Michael Lustig.
\newblock Unsupervised deep basis pursuit: Learning reconstruction without
  ground-truth data.
\newblock In \emph{Proc. Intl. Soc. Mag. Reson. Med}, volume~27, page 0660,
  2019.

\bibitem[Dardikman-Yoffe and Eldar(2020)]{dardikman2020learned}
Gili Dardikman-Yoffe and Yonina~C Eldar.
\newblock Learned sparcom: Unfolded deep super-resolution microscopy.
\newblock \emph{arXiv preprint arXiv:2004.09270}, 2020.

\bibitem[Solomon et~al.(2019)Solomon, Cohen, Zhang, Yang, He, Luo, van Sloun,
  and Eldar]{solomon2019deep}
Oren Solomon, Regev Cohen, Yi~Zhang, Yi~Yang, Qiong He, Jianwen Luo, Ruud~JG
  van Sloun, and Yonina~C Eldar.
\newblock Deep unfolded robust pca with application to clutter suppression in
  ultrasound.
\newblock \emph{IEEE transactions on medical imaging}, 39\penalty0
  (4):\penalty0 1051--1063, 2019.

\bibitem[Sun et~al.(2017)Sun, Chen, Shi, Hong, Fu, and
  Sidiropoulos]{sun2017learning}
Haoran Sun, Xiangyi Chen, Qingjiang Shi, Mingyi Hong, Xiao Fu, and Nikos~D
  Sidiropoulos.
\newblock Learning to optimize: Training deep neural networks for wireless
  resource management.
\newblock In \emph{2017 IEEE 18th International Workshop on Signal Processing
  Advances in Wireless Communications (SPAWC)}, pages 1--6. IEEE, 2017.

\bibitem[Chowdhury et~al.(2020)Chowdhury, Verma, Rao, Swami, and
  Segarra]{chowdhury2020unfolding}
Arindam Chowdhury, Gunjan Verma, Chirag Rao, Ananthram Swami, and Santiago
  Segarra.
\newblock Unfolding wmmse using graph neural networks for efficient power
  allocation.
\newblock \emph{arXiv preprint arXiv:2009.10812}, 2020.

\bibitem[Wang et~al.(2018{\natexlab{a}})Wang, Yang, and Ma]{wang2018velocity}
Wenlong Wang, Fangshu Yang, and Jianwei Ma.
\newblock Velocity model building with a modified fully convolutional network.
\newblock In \emph{SEG Technical Program Expanded Abstracts 2018}, pages
  2086--2090. Society of Exploration Geophysicists, 2018{\natexlab{a}}.

\bibitem[Yang and Ma(2019)]{yang2019deep}
Fangshu Yang and Jianwei Ma.
\newblock Deep-learning inversion: A next-generation seismic velocity model
  building method.
\newblock \emph{Geophysics}, 84\penalty0 (4):\penalty0 R583--R599, 2019.

\bibitem[Yang et~al.(2020{\natexlab{a}})Yang, Pham, Gupta, Unser, and
  Ma]{yang2020deep}
Fangshu Yang, Thanh-an Pham, Harshit Gupta, Michael Unser, and Jianwei Ma.
\newblock Deep-learning projector for optical diffraction tomography.
\newblock \emph{Optics Express}, 28\penalty0 (3):\penalty0 3905--3921,
  2020{\natexlab{a}}.

\bibitem[Zhang et~al.(2020)Zhang, Yang, and Ma]{zhang2020can}
Hao Zhang, Xiuyan Yang, and Jianwei Ma.
\newblock Can learning from natural image denoising be used for seismic data
  interpolation?
\newblock \emph{Geophysics}, 85\penalty0 (4):\penalty0 WA115--WA136, 2020.

\bibitem[Peng et~al.(2016)Peng, Xiao, Feng, Yau, and Yi]{peng2016deep}
Xi~Peng, Shijie Xiao, Jiashi Feng, Wei-Yun Yau, and Zhang Yi.
\newblock Deep subspace clustering with sparsity prior.
\newblock In \emph{IJCAI}, pages 1925--1931, 2016.

\bibitem[Peng et~al.(2018)Peng, Tsang, Zhou, and Zhu]{peng2018k}
Xi~Peng, Ivor~W Tsang, Joey~Tianyi Zhou, and Hongyuan Zhu.
\newblock k-meansnet: When k-means meets differentiable programming.
\newblock \emph{arXiv preprint arXiv:1808.07292}, 2018.

\bibitem[I{\c{s}}{\i}l et~al.(2019)I{\c{s}}{\i}l, Oktem, and
  Ko{\c{c}}]{icsil2019deep}
{\c{C}}a{\u{g}}atay I{\c{s}}{\i}l, Figen~S Oktem, and Aykut Ko{\c{c}}.
\newblock Deep iterative reconstruction for phase retrieval.
\newblock \emph{Applied Optics}, 58\penalty0 (20):\penalty0 5422--5431, 2019.

\bibitem[Wang et~al.(2018{\natexlab{b}})Wang, Roux, Wang, and
  Hershey]{wang2018end}
Zhong-Qiu Wang, Jonathan~Le Roux, DeLiang Wang, and John~R Hershey.
\newblock End-to-end speech separation with unfolded iterative phase
  reconstruction.
\newblock \emph{arXiv preprint arXiv:1804.10204}, 2018{\natexlab{b}}.

\bibitem[Lohit et~al.(2019)Lohit, Liu, Mansour, and
  Boufounos]{lohit2019unrolled}
Suhas Lohit, Dehong Liu, Hassan Mansour, and Petros~T Boufounos.
\newblock Unrolled projected gradient descent for multi-spectral image fusion.
\newblock In \emph{ICASSP 2019-2019 IEEE International Conference on Acoustics,
  Speech and Signal Processing (ICASSP)}, pages 7725--7729. IEEE, 2019.

\bibitem[Zhang et~al.(2019{\natexlab{b}})Zhang, Wang, and
  Giannakis]{zhang2019real}
Liang Zhang, Gang Wang, and Georgios~B Giannakis.
\newblock Real-time power system state estimation and forecasting via deep
  unrolled neural networks.
\newblock \emph{IEEE Transactions on Signal Processing}, 67\penalty0
  (15):\penalty0 4069--4077, 2019{\natexlab{b}}.

\bibitem[Shrivastava et~al.(2020)Shrivastava, Chen, Chen, Lan, Aluru, Liu, and
  Song]{Shrivastava2020GLAD:}
Harsh Shrivastava, Xinshi Chen, Binghong Chen, Guanghui Lan, Srinivas Aluru,
  Han Liu, and Le~Song.
\newblock Glad: Learning sparse graph recovery.
\newblock In \emph{International Conference on Learning Representations}, 2020.
\newblock URL \url{https://openreview.net/forum?id=BkxpMTEtPB}.

\bibitem[Anil et~al.(2019)Anil, Lucas, and Grosse]{anil2019sorting}
Cem Anil, James Lucas, and Roger Grosse.
\newblock Sorting out lipschitz function approximation.
\newblock In \emph{International Conference on Machine Learning}, pages
  291--301. PMLR, 2019.

\bibitem[Wu et~al.(2019{\natexlab{b}})Wu, Guo, Li, and Zhang]{wu2019sparse}
Kailun Wu, Yiwen Guo, Ziang Li, and Changshui Zhang.
\newblock Sparse coding with gated learned ista.
\newblock In \emph{International Conference on Learning Representations},
  2019{\natexlab{b}}.

\bibitem[Yang et~al.(2020{\natexlab{b}})Yang, Gu, Chen, Ma, and
  So]{yang2020learning}
Chengzhu Yang, Yuantao Gu, Badong Chen, Hongbing Ma, and Hing~Cheung So.
\newblock Learning proximal operator methods for nonconvex sparse recovery with
  theoretical guarantee.
\newblock \emph{IEEE Transactions on Signal Processing}, 2020{\natexlab{b}}.

\bibitem[Zarka et~al.(2019)Zarka, Thiry, Angles, and Mallat]{zarka2019deep}
John Zarka, Louis Thiry, Tom{\'a}s Angles, and St{\'e}phane Mallat.
\newblock Deep network classification by scattering and homotopy dictionary
  learning.
\newblock \emph{arXiv preprint arXiv:1910.03561}, 2019.

\bibitem[Chen et~al.(2020{\natexlab{g}})Chen, Liu, Ye, and
  Zhang]{chen2020learnable}
Yunmei Chen, Hongcheng Liu, Xiaojing Ye, and Qingchao Zhang.
\newblock Learnable descent algorithm for nonsmooth nonconvex image
  reconstruction.
\newblock \emph{arXiv preprint arXiv:2007.11245}, 2020{\natexlab{g}}.

\bibitem[Liu et~al.(2018{\natexlab{b}})Liu, Ma, Wang, and
  Zhang]{liu2018learning}
Risheng Liu, Long Ma, Yiyang Wang, and Lei Zhang.
\newblock Learning converged propagations with deep prior ensemble for image
  enhancement.
\newblock \emph{IEEE Transactions on Image Processing}, 28\penalty0
  (3):\penalty0 1528--1543, 2018{\natexlab{b}}.

\bibitem[Metz et~al.(2020)Metz, Maheswaranathan, Freeman, Poole, and
  Sohl-Dickstein]{metz2020tasks}
Luke Metz, Niru Maheswaranathan, C~Daniel Freeman, Ben Poole, and Jascha
  Sohl-Dickstein.
\newblock Tasks, stability, architecture, and compute: Training more effective
  learned optimizers, and using them to train themselves.
\newblock \emph{arXiv preprint arXiv:2009.11243}, 2020.

\bibitem[Behboodi et~al.(2020)Behboodi, Rauhut, and
  Schnoor]{behboodi2020generalization}
Arash Behboodi, Holger Rauhut, and Ekkehard Schnoor.
\newblock Generalization bounds for deep thresholding networks.
\newblock \emph{arXiv preprint arXiv:2010.15658}, 2020.

\bibitem[Van~Luong et~al.(2020)Van~Luong, Joukovsky, and
  Deligiannis]{van2020interpretable}
Huynh Van~Luong, Boris Joukovsky, and Nikos Deligiannis.
\newblock Interpretable deep recurrent neural networks via unfolding reweighted
  $\ell_1$-$\ell_1$ minimization: Architecture design and generalization
  analysis.
\newblock \emph{arXiv preprint arXiv:2003.08334}, 2020.

\bibitem[Schwab et~al.(2019)Schwab, Antholzer, and Haltmeier]{schwab2019deep}
Johannes Schwab, Stephan Antholzer, and Markus Haltmeier.
\newblock Deep null space learning for inverse problems: convergence analysis
  and rates.
\newblock \emph{Inverse Problems}, 35\penalty0 (2):\penalty0 025008, 2019.

\bibitem[Martin et~al.(2001)Martin, Fowlkes, Tal, and
  Malik]{martin2001database}
David Martin, Charless Fowlkes, Doron Tal, and Jitendra Malik.
\newblock A database of human segmented natural images and its application to
  evaluating segmentation algorithms and measuring ecological statistics.
\newblock In \emph{Proceedings Eighth IEEE International Conference on Computer
  Vision. ICCV 2001}, volume~2, pages 416--423. IEEE, 2001.

\bibitem[Xu and Yin(2013)]{xu2013block}
Yangyang Xu and Wotao Yin.
\newblock A block coordinate descent method for regularized multiconvex
  optimization with applications to nonnegative tensor factorization and
  completion.
\newblock \emph{SIAM Journal on imaging sciences}, 6\penalty0 (3):\penalty0
  1758--1789, 2013.

\bibitem[Moreau and Bruna(2017{\natexlab{b}})]{MoreauBruna2017_understanding}
Thomas Moreau and Joan Bruna.
\newblock Understanding {Trainable} {Sparse} {Coding} with {Matrix}
  {Factorization}.
\newblock 2017{\natexlab{b}}.

\bibitem[Chen et~al.(2019{\natexlab{b}})Chen, Sun, and Yin]{chen2019run}
Yifan Chen, Yuejiao Sun, and Wotao Yin.
\newblock Run-and-inspect method for nonconvex optimization and global
  optimality bounds for r-local minimizers.
\newblock \emph{Mathematical Programming}, 176\penalty0 (1):\penalty0 39--67,
  2019{\natexlab{b}}.

\bibitem[Meng et~al.(2021)Meng, Chen, Jiang, and Wang]{meng2021a}
Tianjian Meng, Xiaohan Chen, Yifan Jiang, and Zhangyang Wang.
\newblock A design space study for lista and beyond.
\newblock In \emph{International Conference on Learning Representations}, 2021.
\newblock URL \url{https://openreview.net/forum?id=GMgHyUPrXa}.

\end{thebibliography}

\end{document}